\documentclass[11pt]{article}

\usepackage[utf8]{inputenc}
\usepackage{amsmath}
\usepackage{amssymb}
\usepackage{amsfonts}
\usepackage{amsthm}
\usepackage{thmtools}
\usepackage{enumerate}
\usepackage[top=3.5cm, bottom=3.5cm, left=3.5cm, right=3.5cm]{geometry}
\usepackage{graphicx}
\usepackage{titlesec}
\usepackage{hyperref}
\usepackage{rotating}
\usepackage{xcolor}
\usepackage{setspace}
\usepackage[T1]{fontenc}
\usepackage{tikz}
\usepackage{comment}
\usepackage[all]{xy}
\usepackage{stmaryrd}
\usepackage[nameinlink,capitalize,noabbrev]{cleveref}
\usepackage{mathrsfs}
\usepackage{imakeidx}
\usepackage{authblk}

\titleformat{\section}
{\filcenter\large\bf} 
{\thesection.\ } 
{1pt} 
{} %

\titleformat{\subsection}[hang]
{\filcenter\bf}
{\thesubsection.}
{1pt}
{}

\titleformat{\subsubsection}[hang]
{\filcenter\bf}
{\thesubsubsection.}
{1pt}
{}

\declaretheoremstyle[bodyfont=\normalfont]{normalbody}

\declaretheorem[numberwithin=section,name=Theorem]{theorem}
\declaretheorem[sibling=theorem,style=normalbody,name=Definition]{definition}
\declaretheorem[sibling=theorem,name=Corollary]{corollary}
\declaretheorem[sibling=theorem,name=Lemma]{lemma}
\declaretheorem[sibling=theorem,name=Proposition]{proposition}

\declaretheorem[sibling=theorem,style=normalbody,name=Example]{example}
\declaretheorem[sibling=theorem,style=normalbody,name=Remark]{remark}

\newcommand{\Z}{\mathbb{Z}}
\newcommand{\N}{\mathbb{N}}
\newcommand{\Q}{\mathbb{Q}}
\newcommand{\R}{\mathbb{R}}
\newcommand{\C}{\mathbb{C}}

\renewcommand{\P}{\mathbb{P}}
\newcommand{\G}{\mathbb{G}}

\newcommand{\Hom}{\operatorname{Hom}}
\newcommand{\Aut}{\operatorname{Aut}}
\newcommand{\Gal}{\operatorname{Gal}}

\newcommand{\GL}{\mathrm{GL}}
\newcommand{\PGL}{\mathrm{PGL}}
\newcommand{\Div}{\mathrm{Div}}
\newcommand{\divr}{\mathrm{div}}

\newcommand{\id}{\mathrm{id}}

\newcommand{\Spec}{\mathrm{Spec}}
\newcommand{\A}{\mathbb{A}}
\newcommand{\PPDiv}{\mathrm{PPDiv}}
\newcommand{\D}{\mathfrak{D}}
\newcommand{\E}{\mathfrak{E}}
\renewcommand{\S}{\mathfrak{S}}
\newcommand{\cadiv}{\mathrm{CaDiv}}
\newcommand{\Gm}{\mathbb{G}_{\mathrm{m}}}
\newcommand{\cone}{\mathrm{cone}}

\newcommand{\SAut}{\mathrm{SAut}}
\newcommand{\face}{\mathrm{face}}

\newcommand{\sep}[1]{#1^{\mathrm{sep}}}
\newcommand{\Pol}{\mathrm{Pol}}

\newcommand{\Loc}{\mathrm{Loc}}
\renewcommand{\H}{\mathrm{H}}

\newcommand{\Addresses}{{
  \bigskip
  \footnotesize

  G.~Martínez-Núñez, \textsc{Departamento de Matemáticas, Facultad de Ciencias, Universidad de Chile}\par\nopagebreak \textit{e-mail address}, G.~Martínez-Núñez: \texttt{gary.martinez@ug.uchile.cl}

}}

\DeclareFontFamily{U}{wncy}{}
\DeclareFontShape{U}{wncy}{m}{n}{<->wncyr10}{}
\DeclareSymbolFont{mcy}{U}{wncy}{m}{n}
\DeclareMathSymbol{\Sh}{\mathord}{mcy}{"58}

\DeclareFontFamily{U}{wncy}{}
\DeclareFontShape{U}{wncy}{m}{n}{<->wncyr10}{}
\DeclareSymbolFont{mcy}{U}{wncy}{m}{n}
\DeclareMathSymbol{\Ch}{\mathord}{mcy}{"51}

\title{\textsc{Divisorial fans and algebraic torus actions over arbitrary fields}}

\author[]{\textsc{Gary Mart\'inez-N\'unez}}

\begin{document}

\maketitle

\begin{abstract}
We provide an algebro-geometric combinatorial description of normal varieties endowed with an effective action of an algebraic torus over arbitrary fields. This description is achieved in terms of divisorial fans endowed with a Galois semilinear action. This work completes the classification of normal $T$-varieties over fields.
\end{abstract}

		\mbox{\hspace{1.6em} \textbf{ Keywords:} normal varieties, torus actions, Galois descent.}

\vspace{1em}
		\mbox{\hspace{1.6em} \textbf{ MSC codes (2020):} 14L24, 14L30.} 
		\mbox{\hspace{1.3em} }

\tableofcontents

%
%
%
%

\section{Introduction}
\noindent

In 2006, an algebro-geometric combinatorial description for affine normal $T$-varieties over algebraically closed fields of characteristic zero emerged. These varieties are described in terms of \textit{proper polyhedral divisors}, as introduced by Altmann and Hausen \cite{AH06}. A generalization to complexity-one affine normal $T$-varieties over arbitrary fields was achieved by Langlois in \cite{Lan15}, and a generalization to any field of characteristic zero was accomplished by Gillard in \cite{Gil22a} and \cite{Gil22b} for arbitrary complexity. In \cite{MN25}, a generalization to arbitrary fields and any complexity was established. A proper polyhedral divisor is a formal sum
\[
\D := \sum \Delta_{D} \otimes D,
\]
where the $\Delta_{D}$'s are polyhedra in $N_{\Q}$ with tail cone $\omega$, with $N$ a lattice, and the $D$'s are prime divisors in $\cadiv_{}(Y)$, where $Y$ is a normal semi-projective variety over $k$. Given a proper polyhedral divisor $\D$ (hereafter referred to as a pp-divisor), we can associate with it an affine normal variety $X(\D)$ endowed with an effective action of $T := \Spec(k[M])$. Moreover, it can be shown that such a variety is normal and that every affine normal $T$-variety arises in this way.

\begin{theorem}\cite[Theorem 1.2]{MN25}\label{theoremmainMN25}
Let $k$ be a field and $T$ be a split torus over $k$.
\begin{enumerate}[i)]
\item The scheme $X(\D)$ is a normal $k$-variety with an effective action of $T := \Spec(k[M])$.
\item Let $X$ be an affine normal $k$-variety with an effective $T$-action. Then, there exists a pp-divisor $\D$ such that $X\cong X(\D)$ as $T$-varieties.
\end{enumerate}
\end{theorem}

In the non-affine case, over an algebraically closed field of characteristic zero, Altmann, Hausen, and S\"u\ss\ \cite{AHS08} extended the well-known results in toric geometry, stating that normal toric varieties are determined by fans on the cocharacter lattice of the torus (see, for instance, \cite{Ful93} or \cite{CLS11}). The objects extending fans in toric geometry are referred to as \textit{divisorial fans}. A divisorial fan is a finite set of pp-divisors defined over the same normal semi-projective variety that is stable under intersections, and the intersection of any two of them is a \emph{face} of both (cf. \cref{definition divisorial fan}). The face relation is equivalent to the condition that the corresponding equivariant morphism of varieties induced by the face relation is an open embedding, as in toric geometry.

By Sumihiro Theorem \cite{Sum74}, a normal $T$-variety, where $T$ is a split torus, admits a $T$-stable affine open covering. For each $T$-stable affine open subvariety, by \cite[Theorems 3.1 and 3.4]{AH06}, there exists a pp-divisor on a normal semi-projective variety; however, this normal semi-projective variety may not necessarily be the same for all the pp-divisors. Nonetheless, Altmann, Hausen, and S\"u\ss\ \cite{AHS08} proved that all these pp-divisors can be modified to reside on the same normal semi-projective variety. Under this modification, the pp-divisors fit together into a divisorial fan. They further demonstrated that, from such a divisorial fan, the normal $T$-variety can be fully recovered.

\begin{theorem}\cite[Theorem 5.6]{AHS08}\label{theoremmainalthausus}
Let $k$ be an algebraically closed field of characteristic zero. Up to equivariant isomorphism, every normal variety endowed with an effective algebraic torus action arises from a divisorial fan on a normal semi-projective variety $Y$.
\end{theorem}

This result generalizes those from Mumford \cite{KKMSD73}, Timashev \cite{Tim97}, Pinkham \cite{Pin77} and Flenner and Zaidenberg \cite{FZ03} for complexity-one normal normal $T$-varieties.

Once this theory was developed, numerous contributions emerged. For example, Hausen and S\"u\ss\ \cite{HS10} studied the Cox ring of normal $T$-varieties; Petersen and S\"u\ss\ \cite{PS11} focused on $T$-invariant divisors; Ilten and S\"u\ss\ \cite{IS11} explored polarized $T$-varieties of complexity one; Liendo and S\"u\ss\ \cite{LS13} analyzed the resolution of singularities for $T$-varieties. Many other results can be found in the survey \cite{AIPSV12}.

When $k$ is no longer algebraically closed, the combinatorial framework needs to be enriched for non-split algebraic tori over $k$, similar to toric geometry \cite{Hur11}. An appropriate language to describe the non-split context is through \textit{Galois semilinear (equivariant) actions} or \textit{Galois semilinear actions} (see: \cite[Sections 5.1, 5.2 and 5.3]{MN25}), depending on whether the variety is endowed with an action of an algebraic group or not; and we aim to descend it with the action.

From a \textit{Galois semilinear action $\varrho$ on a pp-divisor} $\D$, we obtain a Galois semilinear equivariant action $X(\varrho)$ on $X(\D)$ and, therefore, a Galois descent datum over $X(\D)$. Nevertheless, not every equivariant Galois descent datum over $X(\D)$ induces a Galois semilinear action on the pp-divisor $\D$. However, there always exists another pp-divisor $\D'$ encoding the same affine normal variety $X(\D)$ such that any equivariant Galois descent datum over $X(\D')$ induces a Galois semilinear action on $\D'$, for instance a \emph{minimal} pp-divisor (\cite[Definition 8.7]{AH06} or \cite[Section 4.4.3]{MN25}). Thus, we have the following result.

\begin{theorem}\cite[Theorem 1.3]{MN25}\label{maintheoremofpaperaffine}
Let $k$ be a field and $L/k$ be a finite Galois extension with Galois group $\Gamma_{L}$.
\begin{enumerate}[a)]
\item\label{theorem main affine minimal part i} Let $(\D_{L},\varrho)$ be a pair consisting of a pp-divisor $\D_{L}$ and a Galois semilinear action on $\D$. Then, $X(\D_{L},\varrho)$ is an affine normal variety endowed with an effective action of an algebraic torus $T$ over $k$ such that $T$ splits over $L$ and $X(\D_{L},\varrho)_{L}\cong X(\D_{L})$ as $T_{\D_{L}}$-varieties over $L$.
\item\label{theorem main affine minimal part ii} Let $X$ be an affine normal variety over $k$ endowed with an effective $T$-action such that $T_{L}$ is split. Then, there exists a pair $(\D_{L},\varrho)$ as above such that $X \cong X(\D_{L},\varrho)$ as $T$-varieties.
\end{enumerate}
\end{theorem}

\subsection*{Main results}
\noindent

Several works seeking generalizations of \cite{AH06} exist, as mentioned above, but not for \cite{AHS08}. In this work, we prove that normal $T$-varieties over arbitrary fields with split torus actions, as in \cref{theoremmainalthausus}, also arise from divisorial fans over a normal semi-projective variety over the ground field.

\begin{theorem}\label{maintheoremofpaperspli}
Let $k$ be a field and $T$ be a split torus over $k$. Up to equivariant isomorphism, every normal $T$-variety arises from a divisorial fan on a normal semi-projective variety $Y$ over $k$.
\end{theorem}

On the one hand, a \textit{Galois semilinear action on a divisorial fan} $(\S,Y)$ induces a Galois semilinear equivariant action on $X(\S)$, therefore, a Galois descent datum over $X(\S)$, the normal $T$-variety encoded by the divisorial fan. On the other hand, not every Galois descent datum over $X(\S)$ induces a Galois semilinear action on the divisorial fan $(\S,Y)$. However, if the Galois group is finite, there always exists a divisorial fan $(\S',Y')$ such that $(\S',Y')$ admits a compatible Galois semilinear action and $X(\S)\cong X(\S')$ for the given Galois descent datum. Not every Galois descent datum over a variety is effective, the respective variety has to admit a Galois stable quasi-projective open covering in order to do so. Hence, we need to impose an additional condition. Thus, we prove the following result, which is the main theorem of this work and is also a generalization of \cite{Hur11}.

\begin{theorem}\label{maintheoremofpaper}
Let $k$ be a field and $L/k$ be a finite Galois extension with Galois group $\Gamma_{L}$.
\begin{enumerate}[a)]
\item \label{maintheoremofpaper part a} Let $T$ be a split algebraic torus over $L$ and $X$ be a normal $T$-variety over $L$. If there exists a divisorial fan $(\S_{L},Y_{L})$ for $X$ admitting a Galois semilinear action such that the subvariety $X(\S(\D,\Gamma_{L}))$ is quasi-projective for every $\D\in\S_{L}$, then there exists an algebraic torus $T'$ over $k$ and a normal $T'$-variety $X'$ over $k$ such that $X'_{L}\cong X$ as $T$-varieties over $L$.
\item \label{maintheoremofpaper part b} Let $T$ be an algebraic torus over $k$ that splits over $L$. Let $X$ be a normal variety endowed with an effective $T$-action over $k$. Then, there exists a divisorial fan $(\S_{L},Y_{L})$ admitting a Galois semilinear action such that the subvariety $X(\S(\D,\Gamma_{L}))$ is quasi-projective for every $\D\in\S_{L}$, and $X_{L}\cong X(\S_{L})$ as $T_{L}$-varieties.
\end{enumerate}
\end{theorem}

In \cref{maintheoremofpaper}, $\S(\D,\Gamma_{L})$ stands for the \textit{sub divisorial fan generated by $\D$ and $\Gamma_{L}$}, cf. \cref{definition semilinear orbit divisorial sub fan}. Moreover, the variety $X(\S(\D,\Gamma_{L}))$ corresponds to the Galois orbit of the affine variety $X(\D)$.

\paragraph{Conventions and settings.} Throughout this paper, $k$ denotes an arbitrary field. We fix an algebraic closure $\bar{k}$, and consequently, we define $\sep{k} \subset \bar{k}$ as the separable closure of $k$. Additionally, $k \subset L \subset \sep{k}$ denotes a finite Galois extension, and we denote its Galois group by $\Gamma_L$.

A \emph{variety} over $k$ (or simply a $k$-variety) is a geometrically integral separated scheme of finite type over $k$. Furthermore, $T$ denotes an algebraic torus over $k$. By extension, a \emph{$T$-variety} over $k$ is a variety equipped with an effective action of $T$.

For any extension $K/k$, we denote $X_{K} := X \times_{k} K$ the corresponding base change.

From convex geometry, $N$ denotes a lattice and $M:=\Hom_{\Z}(N,\Z)$ its dual. These come with a natural perfect pairing $\langle\bullet,\bullet\rangle:M\times N\to \Z$. Moreover, we denote $N_{\Q}:=N\otimes_{\Z}\Q$ and $M_{\Q}:=M\otimes_{\Z}\Q$ for the respective $\Q$-vector spaces. Notice that $M_{\Q}\cong \Hom_{\Q}(N_{\Q},\Q)$ canonically and the perfect pairing extends to a perfect pairing $\langle\bullet,\bullet\rangle:M_{\Q}\times N_{\Q}\to \Q$. By a \emph{cone} $\omega$ in $N_{\Q}$, we mean a convex polyhedral cone.

\paragraph{Layout of the paper.} In \cref{Section affine case minimal}, an overview of the affine case is given. The basic concepts of the theory of pp-divisors (see, for instance, \cite{AH06} and \cite{MN25}) are provided, followed by a discussion of the combinatorial description of affine normal $T$-varieties in the split case; the section concludes with the case of affine normal $T$-varieties with non-split torus actions.

The main goal of \cref{section divisorial fans} is to give a proof of \cref{maintheoremofpaperspli}. To this end, it begins by providing definitions and conventions throughout this section, leading to the definition of a \emph{divisorial fan}. Hence, in \cref{Section Open embeddings}, open embeddings are described in combinatorial terms, following the parallelism with toric geometry. Once open embeddings are well understood, in \cref{Section Divisorial fan}, the notion of \emph{face} in the world of pp-divisors is introduced, followed by that of \emph{divisorial fan}. In the same section, it is proved that \emph{faces} correspond to $T$-stable affine open subvarieties, again as in toric geometry. In \cref{section: from divisorial fans to normal T-varieties}, the construction of a normal $T$-variety from a divisorial fan is given. In general, it is a \emph{prevariety}, and those divisorial fans encoding varieties are studied in \cref{section: separated divisorial fans}. Finally, \cref{maintheoremofpaperspli} is proved in \cref{section: from normal T-varieties to divisorial fans}.

\cref{Appendix localization} is devoted to the study of the \emph{localization} of the category of pp-divisors in order to obtain an equivalence of categories with the category of affine normal $T$-varieties. This construction is a fundamental step for the proof of \cref{maintheoremofpaper}. Once a localization is achieved, in \cref{Section general pp-divisor} we generalize \cref{maintheoremofpaperaffine} for any pp-divisor, i.e., for non-minimal pp-divisors.

In \cref{section: normal varieties and non-split tori}, the main objective is to give a proof of \cref{maintheoremofpaper}. It begins in \cref{sectioneqaut}, where \emph{semilinear morphisms of divisorial fans} are introduced, followed by \cref{Section semilinear aut div fan}, where the group of semilinear automorphisms of a given divisorial fan is explored. This allows us to study the action of finite groups on divisorial fans in \cref{Section semilinear actions of finite groups}. Finally, in \cref{Section proof main theorem}, \cref{maintheoremofpaper} is proved.

In \cref{Section Complexity one and applications}, a simplification of \cref{maintheoremofpaper} is explored in the context of complexity-one normal $T$-varieties, and some applications are given.

\subsection*{Acknowledgements}
I would like to thank Giancarlo Lucchini-Arteche, \'Alvaro Liendo, Adrien Dubouloz, Ronan Terpereau, and Michel Brion for many enriching conversations. This work was partially supported by ANID via Beca de Doctorado Nacional 2021 Folio 21211482 and ECOS-ANID Project ECOS230044.

\section{The affine case}\label{Section affine case minimal}
\noindent

In this section, we briefly recall some results and constructions from \cite{AH06} and \cite{MN25}, focusing on those relevant to the present work. Moreover, new tools are provided in this work; for instance, \cref{lemma polyhedron face}.

\subsection{Proper polyhedral divisors}
\noindent

As stated in the introduction, a proper polyhedral divisor is a formal sum of prime divisors of a normal variety whose coefficients are polyhedra. In the present section, a brief exposition on convex geometry is provided, followed by one on proper polyhedral divisors and a description of the category of proper polyhedral divisors.

\paragraph{Convex Geometry.} Let $N$ be a lattice and $\omega\subset N_{\Q}$ be a pointed cone, i.e., it contains no lines. A \emph{polyhedron} is a convex set $\Delta\subset N_{\Q}$ that can be written as the \emph{Minkowski sum} of $\Pi\subset N_{\Q}$, the convex hull of a finite set of points, and a cone $\omega\subset N_{\Q}$, i.e., for every $n\in\Delta$ there exist $p\in \Pi$ and $w\in\omega$ such that $n=p+w$. For example, the set $\Delta:=\{(a,b)\in \Q^{2}\mid a+b\geq 1\}$ can be written as the Minkowski sum of $\Pi$, the convex hull of $\{(1,0),(0,1)\}$, and the cone $\omega:=\cone((1,0),(0,1))$.

\begin{center}
\begin{tikzpicture}
\fill[gray] (0,0.7) -- (0,2) -- (2,2) -- (2,0) -- (0.7,0);
\node[scale=0.6] at (0.7,-0.2) {$(1,0)$};
\node[scale=0.6] at (-0.3,0.7) {$(0,1)$};
\draw[->] (0,0) -- (2,0);
\draw[->] (0,0) -- (0,2);
\node[scale=0.9] at (1,1) {$\Delta$};

\node[scale=0.9] at (3,1) {$=$} ;

\draw[-] (4,0.7) -- (4.7,0);
\node[scale=0.9] at (4.35,1) {$\Pi$};
\node[scale=0.6] at (4.7,-0.2) {$(1,0)$};
\node[scale=0.6] at (3.7,0.7) {$(0,1)$};

\node[scale=0.9] at (6,1) {$+$} ;

\fill[gray] (7,0) -- (7,2) -- (9,2) -- (9,0) -- (7.7,0);
\draw[->] (7,0) -- (9,0);
\draw[->] (7,0) -- (7,2);
\node[scale=0.9] at (7.85,1) {$\omega$};

\end{tikzpicture}
\end{center}

Every polyhedron $\Delta\subset N_{\Q}$ has a \emph{Minkowski decomposition} $\Delta=\Pi+\omega$, and $\omega$ is called the \emph{tail} of $\Delta$. For a polyhedron, its tail is unique in the decomposition. A polyhedron of tail $\omega$ is called an $\omega$-polyhedron. The Minkowski sum of two $\omega$-polyhedra is also an $\omega$-polyhedron. If $\Pi\subset N_{\Q}$ is generated by elements in $N$, the polyhedron $\Delta$ is called \emph{integral}, and the Minkowski sum of two integral $\omega$-polyhedra is also an integral $\omega$-polyhedron. Thus, given a cone $\omega\subset N_{\Q}$, the set of all $\omega$-polyhedra $\Pol_{\omega}^{+}(N_{\Q})$, endowed with the Minkowski sum, has a monoidal structure, with $\omega$ as a neutral element, and the set of all integral $\omega$-polyhedra $\Pol_{\omega}^{+}(N)$ is a submono\"id with the same neutral element. Let us denote by $\Pol_{\omega}(N_{\Q})$ (respectively $\Pol_{\omega}(N)$) the Grothendieck group associated to $\Pol_{\omega}^{+}(N_{\Q})$ (respectively $\Pol_{\omega}^{+}(N)$).

\begin{definition}
The \textit{face of $\Delta\in\mathrm{Pol}_{\omega}^{+}(N)$ defined by} $m\in\omega^{\vee}\cap M$ is
\[\mathrm{face}(\Delta,m):=\left\{v\in\Delta \mid \langle m,v \rangle\leq \langle m,v'\rangle \textrm{ for all }v'\in\Delta\right\}\in\mathrm{Pol}_{\omega\cap m^{\perp}}^{+}(N).\]
\end{definition}

The following lemma is introduced to simplify the proofs in \cref{Section Open embeddings}. In particular, it is used in the proof of \cref{Proposition 4.3 AHS08}.

\begin{lemma}\label{lemma polyhedron face}
Let $\Delta$ and $\Delta'$ be two polyhedra in $\mathrm{Pol}_{\omega}(N)$ and $\mathrm{Pol}_{\omega'}(N)$, respectively. Let $m\in \omega^{\vee}\cap M$. Then, the following holds:
\begin{enumerate}[a)]
\item\label{part a lemma polyhedron face} if $\Delta'\subset \Delta$ and $m\in(\omega')^{\vee}$, then $\Delta'\cap\mathrm{face}\left(\Delta,m\right)\subset \mathrm{face}\left(\Delta',m\right)$,
\item\label{part b lemma polyhedron face} for any pair $\Delta$ and $\Delta'$ such that $\Delta\cup\Delta'$ is a polyhedron and $m\in(\omega')^{\vee}$, we have that
\[\mathrm{face}\left(\Delta\cup\Delta',m\right)\subset \mathrm{face}\left(\Delta,m\right)\cup \mathrm{face}\left(\Delta',m\right),\]
\item\label{part c lemma polyhedron face} if $m'\in\omega^{\vee}$ and $\face\left(\Delta,m\right)\cap\face\left(\Delta,m'\right)\neq\emptyset$, then
\[\face\left(\face\left(\Delta,m\right),m'\right)=\face\left(\Delta,m\right)\cap\face\left(\Delta,m'\right)=\face\left(\face\left(\Delta,m'\right),m\right).\]
\end{enumerate}
\end{lemma}

\begin{proof}
Let $v\in\Delta'\cap\mathrm{face}\left(\Delta,m\right)$. Then $\langle m,v \rangle\leq \langle m,v'\rangle$ for all $v'\in\Delta$. In particular, $\langle m,v \rangle\leq \langle m,v'\rangle$ for all $v'\in\Delta'$. Thus, $v\in\mathrm{face}\left( \Delta',m \right)$. This proves \eqref{part a lemma polyhedron face}.

By the first part, we have that
\[\Delta\cap\mathrm{face}\left(\Delta\cup\Delta',m\right)\subset \mathrm{face}\left(\Delta,m\right)\]
and
\[\Delta'\cap\mathrm{face}\left(\Delta\cup\Delta',m\right)\subset \mathrm{face}\left(\Delta',m\right).\]
Then,
\[\mathrm{face}\left(\Delta\cup\Delta',m\right)\subset \mathrm{face}\left(\Delta,m\right)\cup \mathrm{face}\left(\Delta',m\right).\]

Let us now prove the last part. For $m\in \omega^{\vee}$, we denote by $\mathrm{eval}_{m}:N_{\Q}\to \Q$ the map defined by $\mathrm{eval}_{m}(v):=\langle m,v\rangle$, for all $v\in N_{\Q}$. Given that $m$ and $m'$ are in $\omega^{\vee}$, we have that $\min\left\{\mathrm{eval}_{m}(\Delta)\right\}$ and $\min\left\{\mathrm{eval}_{m'}(\Delta)\right\}$ exist. Let us denote $a=\min\left\{\mathrm{eval}_{m}(\Delta)\right\}$ and $b=\min\left\{\mathrm{eval}_{m'}(\Delta)\right\}$. Hence, the following equalities hold:
\[\face\left(\Delta,m\right)=\{v\in\Delta\mid \mathrm{eval}_{m}(v)=a\}\quad \mathrm{and}\quad \face\left(\Delta,m'\right)=\{v\in\Delta\mid \mathrm{eval}_{m'}(v)=b\}.\]

On the one hand,  we have that $\face\left(\Delta,m\right)\cap \face\left(\Delta,m'\right)\neq \emptyset$ and
\[\face\left(\Delta,m\right)\cap \face\left(\Delta,m'\right)=\{v\in\Delta\mid \mathrm{eval}_{m}(v)=a\textrm{ and }\mathrm{eval}_{m'}(v)=b \}.\]
On the other hand, since $b\leq\mathrm{eval}_{m'}(v')$ for all $v'\in\Delta$ it follows that
\begin{align*}
\face\left(\face\left(\Delta,m\right),m'\right) &=\{v\in\face\left(\Delta,m\right) \mid  \mathrm{eval}_{m'}(v)\leq \mathrm{eval}_{m'}(v') \textrm{ for all }v'\in\face\left(\Delta,m\right) \} \\
&=\{v\in\face\left(\Delta,m\right) \mid b\leq \mathrm{eval}_{m'}(v)\leq b \} \\
&=\{v\in\Delta\mid \mathrm{eval}_{m}(v)=a\textrm{ and }\mathrm{eval}_{m'}(v)=b \}.
\end{align*}
Hence, we have that $\face\left(\face\left(\Delta,m\right),m'\right)=\face\left(\Delta,m\right)\cap \face\left(\Delta,m'\right)$.
\end{proof}

\paragraph{Proper polyhedral divisors.} Let $Y$ be a normal variety over $k$. Given that $\mathrm{Pol}_{\omega}(N_{\Q})$ and $\mathrm{Pol}_{\omega}(N)$ are abelian groups, we can take the tensor products
\[\mathrm{Pol}_{\omega}(N_{\Q})\otimes_{\Z}\cadiv(Y)\quad\mathrm{ and }\quad\mathrm{Pol}_{\omega}(N_{})\otimes_{\Z}\cadiv(Y),\]
the group of \textit{rational (resp. integral) polyhedral Cartier divisors}, denoted by $\cadiv_{\Q}(Y,\omega)$ and $\cadiv_{}(Y,\omega)$ respectively, and
\[\mathrm{Pol}_{\omega}(N_{\Q})\otimes_{\Z}\Div(Y)\quad\mathrm{ and }\quad\mathrm{Pol}_{\omega}(N_{})\otimes_{\Z}\Div(Y),\]
the group of \textit{rational (resp. integral) Weil divisors}, denoted by $\mathrm{Div}_{\Q}(Y,\omega)$ and $\mathrm{Div}_{}(Y,\omega)$ respectively.

Recall that, for a normal variety $Y$ over $k$, there is a canonical embedding
\[\cadiv(Y)\to \mathrm{Div}(Y),\]
which allows us to consider $\cadiv(Y)\subset\mathrm{Div}(Y)$ and, therefore, $\cadiv_{\Q}(Y,\omega)\subset\mathrm{Div}_{\Q}(Y,\omega)$ and $\cadiv_{}(Y,\omega)\subset\mathrm{Div}_{}(Y,\omega)$. In particular, we can require $D\in\cadiv(Y)$ to be prime. With this in mind, note that we can always write an element in any of these groups as $\D=\sum_D \Delta_{D}\otimes D$, where the sum runs through the prime divisors $D$ of $Y$ and the $\Delta_{D}$'s are elements in $\mathrm{Pol}_{\omega}(N_{})$ or $\mathrm{Pol}_{\omega}(N_{\Q})$.

By a \textit{polyhedral divisor with respect to} $\omega\in N_{\Q}$ we mean an element of the group $\mathrm{Div}_{\Q}(Y,\omega)$ and $\omega$ is called the \emph{tail} of the pp-divisor. They are also called \emph{$\omega$-tailed polyhedral divisor}.

\begin{definition}\label{definitionppdiv}
Let $Y$ be a normal $k$-variety, $N$ be a lattice and $\omega\subset N_{\Q}$ a pointed cone. A polyhedral divisor $\D=\sum_D \Delta_{D}\otimes D\in\mathrm{CaDiv}_{\Q}(Y,\omega)$ is called \textit{proper} if
\begin{enumerate}
\item all the $D\in\Div(Y)$ are prime divisors and the $\Delta_{D}$ are in $\mathrm{Pol}_{\omega}^{+}(N_{\Q})$;
\item for every $m\in\mathrm{relint}(\omega^{\vee})\cap M$, the evaluation
\[\D(m):=\sum \min_{v\in\Delta_{D}}\{\langle m,v\rangle\}D\in\mathrm{CaDiv}_{\Q}(Y)\] is a big divisor on $Y$, i.e., for some $n\in\N$ there exists a section $f\in \H^{0}(Y,\mathscr{O}(n\D(m)))$ such that $Y_{f}$ is affine;
\item for every $m\in\omega^{\vee}\cap M$, the evaluation $\D(m)\in\cadiv_{\Q}(Y)$ is semiample, i.e., it admits a basepoint-free multiple. In other words, for some $n\in\N$ the sets $Y_{f}$ cover $Y$, where $f\in \H^{0}(Y,\mathscr{O}(n\D(m)))$.
\end{enumerate}
The cone $\omega$ is called the \emph{tail cone} of $\D$. Given a pointed cone $\omega\subset N_{\Q}$, the mono\"id of the $\omega$-tailed proper polyhedral divisors (pp-divisors for short) is denoted by $\PPDiv_{\Q}(Y,\omega)$. The mono\"id is partially ordered as follows: if $\D=\sum_D \Delta_{D}\otimes D$ and $\D'=\sum_D \Delta_{D}'\otimes D$, then $\D'\leq \D$ if and only if $\Delta_{D}\subset \Delta_{D}'$ for every $D$.
\end{definition}

\begin{center}
\begin{tikzpicture}
\fill[gray] (1/3,0) -- (1/2,2) -- (2,2) -- (2,0) -- (1/3,0);
\draw[->] (0,0) -- (2,0);
\draw[->] (0,0) -- (0,2);
\node[scale=0.9] at (1,0.6) {$\Delta_{0}$};

\fill[gray] (11/4,0) -- (35/12,2) -- (5,2) -- (5,0);
\node[scale=0.6] at (3.7,-0.2) {$(1,0)$};
\node[scale=0.6] at (2.7,0.7) {$(0,1)$};
\draw[->] (3,0) -- (5,0);
\draw[->] (3,0) -- (3,2);
\node[scale=0.9] at (3.6,0.6) {$\Delta_{1}$};

\fill[gray] (6,0) -- (6,1) --(73/12,2) -- (8,0) -- (6.7,0);
\node[scale=0.6] at (6.7,-0.2) {$(1,0)$};
\node[scale=0.6] at (5.7,0.7) {$(0,1)$};
\draw[->] (6,0) -- (8,0);
\draw[->] (6,0) -- (6,2);
\node[scale=0.9] at (6.6,0.6) {$\Delta_{\infty}$};
\end{tikzpicture}
\end{center}

\paragraph{The category of pp-divisors.} There exists the category of pp-divisors $\mathfrak{PPDiv}(k)$ (see \cite[Section 8]{AH06} and \cite[Section 3]{MN25}). The morphisms in this category are defined by three pieces of data. Among them, one is called a \textit{plurifunction}. A \textit{plurifunction} with respect to the lattice $N$ and a normal semi-projective variety $Y$ is an element
\[\mathfrak{f}=\sum v_{i}\otimes f_{i}\in k(Y,N)^{*}:=N\otimes_{\Z}k(Y)^{*}.\]
For $m\in M=\Hom(N,\Z)$, the \textit{evaluation} of $\mathfrak{f}$ is given by
\[\mathfrak{f}(m):=\prod f_{i}^{\langle m,v_{i} \rangle}\in k(Y)^{*}.\]
Besides, we can associate with it a \textit{polyhedral principal divisor} with respect to a cone $\omega\subset N_{\Q}$, which is given by
\[\mathrm{div}(\mathfrak{f}):=\sum(v_{i}+\omega)\otimes\mathrm{div}(f_{i})\in\cadiv(Y,\omega).\]

Notice that a morphism of lattices $F:N\to N'$ induces a morphism between the groups $F_{*}:k(N,Y)^{*}\to k(N',Y)^{*}$ given by
\[F_{*}\left(\sum v_{i}\otimes f_{i}\right)=\sum F(v_{i})\otimes f_{i}.\]

A morphism $\psi:Y\to Y'$ also induces a morphism between the respective groups of plurifunctions $\psi^{*}:k(N,Y')^{*}\to k(N,Y)^{*}$ given by
\[\psi^{*}\left(\sum v_{i}\otimes f_{i}\right)=\sum v_{i}\otimes \psi^{*}(f_{i}).\]

The mono\"id $\PPDiv_{\Q}(Y,\omega)$ is partially ordered, defined as follows: $\D'\leq \D$ if and only if $\Delta_{D}\subset \Delta_{D}'$ for every $D$.

\begin{definition}\label{def83}\cite[Definition 8.3]{AH06}
Let $Y$ and $Y'$ be normal $k$-varieties, $N$ and $N'$ be lattices and $\omega\subset N$ and $\omega'\subset N'$ be pointed cones. Let us consider
\[\D=\sum \Delta_{i}\otimes D_{i}\in\PPDiv_{\Q}(Y,\omega)\, \textrm{ and }\, \D'=\sum \Delta_{i}'\otimes D_{i}'\in\PPDiv_{\Q}(Y',\omega')\]
two pp-divisors.
\begin{enumerate}[a)]
\item For morphisms $\psi:Y\to Y'$ such that none of the supports $\mathrm{Supp}(D_{i}')$ contains $\psi(Y)$, we define the (not necessarily proper) \textit{polyhedral pullback} as
\[\psi^{*}(\D'):=\sum \Delta_{i}'\otimes \psi^{*}(D_{i}')\in\cadiv_{\Q}(Y,\omega').\]
\item\label{def83b} For linear maps $F:N\to N'$ with $F(\omega)\subset \omega'$, we define the (not necessarily proper) \textit{polyhedral pushforward} as
\[F_{*}(\D):=\sum (F(\Delta_{i})+\omega')\otimes D_{i}\in \cadiv_{\Q}(Y,\omega').\]
\item \label{defmappdiv} A map $\D\to \D'$ is a triple $(\psi,F,\mathfrak{f})$ with a \emph{dominant} morphism $\psi:Y\to Y'$, $F$ a linear map as in \ref{def83b}), and a plurifunction $\mathfrak{f}\in k(Y,N')^{*}$ such that
\[\psi^{*}(\D')\leq F_{*}(\D)+\divr(\mathfrak{f}).\]
\end{enumerate}
\end{definition}

The triple $(\id,\id_{N},1)$ is the identity map $\D\to\D$ for a pp-divisor. The composition of two morphisms $(\psi,F,\mathfrak{f})$ and $(\psi',F',\mathfrak{f}')$ is defined as follows:
\[(\psi',F',\mathfrak{f}')\circ(\psi,F,\mathfrak{f})=(\psi'\circ\psi,F'\circ F,F_{*}'(\mathfrak{f})\cdot \psi^{*}(\mathfrak{f}')).\]
 
 \subsection{Split torus actions}
\noindent

A pp-divisor gives rise to an affine normal variety endowed with an effective action of a split torus. Below, we summarize some of the key constructions and results regarding these objects.

\paragraph{Combinatorial description.} Let $Y$ be a normal semi-projective variety. Let $\D$ be a pp-divisor in $\PPDiv_{\Q}(Y,\omega)$. For any pair of elements $m$ and $m'$ in $\omega^{\vee}\cap M$, it holds that $\D(m)+\D(m')\leq \D(m+m')$ by definition. Hence, this condition induces a product
\begin{align*}
\H^{0}(Y,\mathscr{O}_{Y}(\D(m)))\times \H^{0}(Y,\mathscr{O}_{Y}(\D(m'))) &\to \H^{0}(Y,\mathscr{O}_{Y}(\D(m+m'))) \\
(f,f')&\mapsto f\cdot f',
\end{align*}
as subsets of $k(Y)$. This endows the $\H^{0}(Y,\mathscr{O}_{Y})$-module
\[
A[Y,\D]:=\bigoplus_{m\in\omega^{\vee}\cap M}\H^{0}(Y,\mathscr{O}_{Y}(\D(m)))
\]
with a structure of $M$-graded $k$-algebra. Moreover, the scheme $X(\D):=\Spec(A[Y,\D])$ is an affine normal variety endowed with an effective action of $T := \Spec(k[M])$. Conversely, for an affine normal $T$-variety $X$ over $k$, where $T$ is a split torus over $k$, there exists a pp-divisor $\D\in\PPDiv_{\Q}(Y,\omega)$ such that $X\cong X(\D)$ as $T$-varieties. These are summarized in the following result, corresponding to \cite[Theorems 3.1 and 3.4]{AH06} and \cite[Theorem 1.2]{MN25}.

\begin{proposition}\label{proposition: theoremmainMN25}
Let $T$ be a split torus over $k$.
\begin{enumerate}[i)]
\item The scheme $X(\D)$ is a normal $k$-variety with an effective action of $T := \Spec(k[M])$.
\item Let $X$ be an affine normal $k$-variety with an effective $T$-action. Then, there exists a pp-divisor $\D$ such that $X\cong X(\D)$ as $T$-varieties.
\end{enumerate}
\end{proposition}

\paragraph{Functoriality.} By \cref{proposition: theoremmainMN25}, from any object $\D$ in $\mathfrak{PPDiv}(k)$ we obtain an affine normal $T$-variety $X(\D)$, with $T$ a split torus over $k$. Let $\D$ and $\D'$ be two objects in $\mathfrak{PPDiv}(k)$ and $(\psi,F,\mathfrak{f}):\D'\to\D$ be a morphism of pp-divisors over $k$. This morphism induces a morphism of modules between the $\H^{0}(Y,\mathscr{O}_{Y})$-module $\H^{0}(Y,\mathscr{O}(\D(m)))$ and the $\H^{0}(Y',\mathscr{O}_{Y'})$-module $\H^{0}(Y',\mathscr{O}(\D'(F^{*}(m))))$ given by
\begin{align*}
\H^{0}(Y,\mathscr{O}(\D(m))) &\to \H^{0}(Y',\mathscr{O}(\D'(F^{*}(m)))), \\
h &\mapsto \mathfrak{f}(m)\psi^{*}(h).
\end{align*}
Hence, all these morphisms fit together into a graded morphism of $k$-algebras
\[A[Y,\D]=\bigoplus_{m\in \omega^{\vee}\cap M}\H^{0}(Y,\mathscr{O}(\D(m)))\to \bigoplus_{m\in \omega'^{\vee}\cap M'}\H^{0}(Y',\mathscr{O}(\D'(m)))=A[Y',\D'],\]
which induces an equivariant morphism
\[X(\psi,F,\mathfrak{f}):=(\varphi,f):X(\D)\to X(\D'),\]
where $\varphi:T'\to T$ is determined by $F:N'\to N$.

In what follows, $\mathcal{E}(k)$ denotes the category of affine normal varieties endowed with an effective action of a split torus over $k$.

\begin{proposition}\label{Proposition X is faithful}
The assignment $\D\mapsto X(\D)$ defines a faithful covariant functor $X:\mathfrak{PPDiv}(k)\to\mathcal{E}(k)$.
\end{proposition}

\subsection{Non-split torus actions}
\noindent

As in toric geometry \cite{Hur11}, combinatorial data alone does not describe non-split torus actions. The missing information is encoded by \emph{Galois semilinear actions}, or \emph{Galois semilinear equivariant actions}, which, under mild hypotheses, can be translated into the algebro-geometric combinatorial framework.

The following results can be found in \cite[Section 5.2]{MN25} and \cite[Section 6.1]{MN25} (based on \cite[Section 14.20]{GW20} and \cite{Bor20}).

\paragraph{Semilinear morphism.} Recall that $L/k$ denotes a finite Galois extension with Galois group $\Gamma_{L}$. Let $X_{1}$ and $X_{2}$ be varieties over $L$ and $\gamma\in\Gamma_{L}$. A \emph{semilinear morphism with respect to $\Gamma_{L}$} is a morphism of schemes $h_{\gamma}:X_{1}\to X_{2}$ such that the following diagram commutes
\[\xymatrix{ X_{1} \ar[r]^{h_{\gamma}} \ar[d] & X_{2} \ar[d] \\ \Spec(L)\ar[r]_{\gamma^{\natural}} & \Spec(L), }\]
where $\gamma^{\natural}:=\Spec(\gamma^{-1})$. It is called a \emph{semilinear isomorphism} if $h_{\gamma}$ is a scheme isomorphism.

\begin{remark}
The choice of $\gamma^{\natural}=\Spec(\gamma^{-1})$ instead of $\Spec(\gamma)$ is to maintain \emph{covariance} in the following sense: for $\gamma$ and $\gamma'$ in $\Gamma_{L}$, and $h_{\gamma}$ and $g_{\gamma'}$ two semilinear morphisms, the composition $h_{\gamma}\circ g_{\gamma'}$ is a semilinear morphism with respect to $\gamma\cdot \gamma'$. More details can be found in \cite[Section 5.2]{MN25}.
\end{remark}

The group of semilinear automorphisms of a variety $X$ is denoted by $\SAut(X)$, which clearly contains $\Aut(X)$ as a subgroup. Moreover, they fit into the following exact sequence
\begin{equation}\label{equation: aut saut}
1 \to \Aut(X) \to \SAut(X)\to \Gamma_{L}.
\end{equation}

In the context of algebraic groups, a semilinear morphism $\varphi_{\gamma}:G_{1}\to G_{2}$ that respects the group structure is called a \emph{semilinear group homomorphism with respect to $\Gamma_{L}$}. The group of semilinear group automorphisms of an algebraic group $G$ is denoted by $\SAut_{\mathrm{gr}}(G)$, and it also contains $\Aut_{\mathrm{gr}}(G)$ as a subgroup. Moreover, they satisfy an exact sequence as in \eqref{equation: aut saut}.

There is also an equivariant version of this morphism. A \emph{semilinear equivariant morphism with respect to $\Gamma_{L}$} is a pair $(\varphi_{\gamma}, h_{\gamma})$ consisting of a semilinear group homomorphism $\varphi_{\gamma}:G_{1}\to G_{2}$ and a semilinear morphism $h_{\gamma}:X_{1}\to X_{2}$ such that the following diagram commutes
\[
\xymatrix{ G_{1}\times X_{1} \ar[r]^{\mu_{1}} \ar[d]_{(\varphi_{\gamma},h_{\gamma})} & X_{1} \ar[d]^{h_{\gamma}} \\
G_{2}\times X_{2} \ar[r]_{\mu_{2}}& X_{2}}
\]
and $\mu_{1}$ and $\mu_{2}$ are the respective actions of $G_{1}$ on $X_{1}$ and $G_{2}$ on $X_{2}$.

For an algebraic group $G$ acting on a variety $X$, denote by $\SAut_{G}(X)$ the group of semilinear automorphisms, which naturally contains $\Aut_{G}(X)$. We set $\SAut(G; X)$ as the subgroup of $\SAut_{\mathrm{gr}}(G)\times \SAut_{G}(X)$ defined as the preimage of the diagonal inclusion $\Gamma_{L}\to\Gamma_{L}\times\Gamma_{L}$. This group fits into the following exact sequence
\begin{equation}\label{equation: autGautX sautGX}
1 \to \Aut_{\mathrm{gr}}(G)\times\Aut(X) \to \SAut(G;X)\to \Gamma_{L},
\end{equation}
where the map on the right is defined as $(\varphi_{\gamma},h_{\gamma})\mapsto \gamma$.

A \emph{Galois semilinear equivariant action on $X$} is a group homomorphism
\[\varrho:\Gamma_{L}\to \SAut(G;X),\]
where $\varrho$ is a section of the exact sequence \eqref{equation: autGautX sautGX}. Moreover, in the particular case when $G$ is trivial, $\varrho$ is simply called a \emph{Galois semilinear action}.

A Galois semilinear equivariant action $\varrho:\Gamma_{L}\to \SAut(G;X)$ comes with natural projections inducing Galois semilinear actions $\varrho_{G}=\varphi$ on $G$ and $\varrho_{X}=h$ on $X$.

\paragraph{Galois descent through semilinear actions.} Let $S$ be a scheme over $k$. The induced scheme $S_{L}$ over $L$ comes with a canonical Galois semilinear action $\varrho_{\mathrm{can}}:\Gamma_{L}\to \SAut(S_{L})$ defined as $\varrho_{\mathrm{can}}(\gamma):=\id_{S}\otimes\Spec(\gamma^{-1}).$ This defines a functor
\begin{align*}
\Phi: \mathbf{Sch}_{k} &\to \mathbf{Sch}(L/k) \\
S &\mapsto (S_{L},\varrho_{\mathrm{can}}),
\end{align*}
from the category of schemes over $k$ to the category of schemes over $L$ endowed with a Galois semilinear action, which is faithful and covariant. Not every pair $(Z,\varrho)$ lies in the essential image of the functor. Those pairs are called \emph{effective}. For example, a pair $(Z,\varrho)$ is always effective whenever $Z$ is quasi-projective. This holds in particular for algebraic groups. In what follows, we address the equivariant version of this problem for varieties. Before doing so, we introduce the following definition.

\begin{definition}
Let $Z$ be a scheme over $L$. A $k$-form of $Z$ is a pair $(S,\nu)$ consisting of a scheme $S$ over $k$ and an isomorphism $\nu:S_{L}\to Z$ of schemes over $L$.
\end{definition}

Let $G'$ be an algebraic group over $k$ and $S$ be a $G'$-variety over $k$. The $G'_{L}$-variety $S_{L}$ comes with a canonical Galois semilinear equivariant action $\varrho_{\mathrm{can}}:\Gamma_{L}\to \SAut(G'_{L};S_{L})$ defined as
\[\varrho_{\mathrm{can}}(\gamma):=(\id_{G}\otimes\Spec(\gamma^{-1}),\id_{S}\otimes\Spec(\gamma^{-1})),\]
for each $\gamma$ in $\Gamma_{L}$.

Meanwhile, for an algebraic group $G$ over $L$, any semilinear action on it defines an effective pair; this does not hold for an arbitrary $G$-variety over $L$. The following result gives a complete answer to this problem.

\begin{proposition}\cite[Proposition 6.14]{MN25}\label{proposition equivalence normal G-varieties}
Let $G$ be a connected algebraic group over $L$ and $\varphi$ be a Galois semilinear action. Let $G'$ be the algebraic group over $k$ associated with $(G,\varphi)$. Then, the functor $\Phi: \mathbf{Sch}_{k} \to \mathbf{Sch}(L/k)$, given by $S\mapsto (S_{L},\varrho_{\mathrm{can}})$, induces an equivalence of categories between the category of normal varieties with effective $G'$-actions and the category of normal varieties over $L$ with effective $G$-actions endowed with Galois semilinear equivariant actions, which are covered by $G$-stable and $\Gamma_{L}$-stable quasi-projective subvarieties and for which $\varrho_{G}$ is isomorphic to $\varphi$.
\end{proposition}

\paragraph{Semilinear morphism of pp-divisors.} Let $\D$ and $\D'$ be two objects in $\mathfrak{PPDiv}(L)$ that are defined over normal varieties $Y$ and $Y'$ respectively. Denote by $\omega_{\D}$ and $\omega_{\D'}$ their respective tail cones. Let $\gamma\in \Gamma_{L}$. A \emph{semilinear morphism of pp-divisors with respect to $\gamma$} is a triple $(\psi_{\gamma},F,\mathfrak{f}):\D\to\D'$, where $\psi_{\gamma}$ is a semilinear dominant morphism, $F:N\to N'$ is a morphism of lattices with $F(\omega_{\D'})\subset \omega_{\D}$, and $\mathfrak{f}\in L(N',Y)^{*}$ is a plurifunction such that
\[ \psi_{\gamma}^{*}(\D')\leq F_{*}(\D)+\mathrm{div}(\mathfrak{f}). \]
The composition of two semilinear morphisms of pp-divisors is a semilinear morphism of pp-divisors. Denote by $\mathfrak{PPDiv}(L/k)$ the category of pp-divisors over $L$ whose morphisms are semilinear morphisms of pp-divisors.

A semilinear morphism of pp-divisors $(\psi_{\gamma},F,\mathfrak{f}):\D\to\D'$ induces a semilinear equivariant automorphism $X(\psi_{\gamma},F,\mathfrak{f}):X(\D)\to X(\D')$. Denote by $\mathcal{E}(L/k)$ the category of affine normal varieties endowed with an effective split torus action whose morphisms are dominant semilinear equivariant morphisms. This assignment $\D\mapsto X(\D)$ defines a functor. Moreover, it holds that

\begin{proposition}\cite[Proposition 5.7]{MN25}\label{proposition: semilinear functoriality}
The functor $X:\mathfrak{PPDiv}(L/k)\to \mathcal{E}(L/k)$ is faithful and covariant.
\end{proposition}

However, this functor turns out not to be an equivalence of categories. Indeed, let $Y$ be a normal projective surface and $\pi: \tilde{Y}\to Y$ be a blow-up such that $\pi$ is not an isomorphism. Let $\D$ be a pp-divisor on $Y$ whose support is outside the exceptional divisor. Hence, $\pi^{*}\D$ is a pp-divisor on $\tilde{Y}$. The semilinear morphism of pp-divisors with respect to the identity in $\Gamma_{L}$ given by $(\pi,\mathrm{id}_{N},\mathfrak{1}):\pi^{*}\D\to \D$ corresponds to the identity $\mathrm{id}:X(\pi^{*}\D)\to X(\D)$. However, $(\pi,\mathrm{id}_{N},\mathfrak{1}):\pi^{*}\D\to \D$ has no inverse while the identity has one. Thus, $X:\mathfrak{PPDiv}(L/k)\to \mathcal{E}(L/k)$ is not an equivalence of categories due to a lack of morphisms. In general, for a pp-divisor $\D$ in $\mathfrak{PPDiv}(L/k)$, the set
\[S(\D):=\{(\psi_{\gamma},F,\mathfrak{f})\in \mathrm{Mor}_{\mathfrak{PPDiv}(L/k)}(\D,\D) \mid X(\psi_{\gamma},F,\mathfrak{f})\textrm{ is an automorphism} \}\]
may not be a group. Nevertheless, there exist pp-divisors such that $S(\D)$ is a group. For instance, if $\D$ is a \emph{minimal} pp-divisor (\cite[Definition 8.7]{AH06} or \cite[Section 4.4.3]{MN25}), then $S(\D)$ is a group and it is denoted by $\SAut(\D)$. In such a case, $\Aut(\D)$ exists and is a subgroup of $\SAut(\D)$; both fit into the exact sequence
\begin{equation}\label{equation: autD sautD}
1 \to \Aut(\D) \to \SAut(\D) \to \Gamma_{L},
\end{equation}
where $(\psi_{\gamma},F,\mathfrak{f})\mapsto \gamma$.

Let $\D$ be a pp-divisor such that $S(\D)$ is a group. A \emph{Galois semilinear action on $\D$} is a group homomorphism $\varrho:\Gamma_{L}\to \SAut(\D)$ such that $\varrho$ is a section of \eqref{equation: autD sautD}.

\paragraph{Combinatorial description.} Let $\D$ be a pp-divisor in $\mathfrak{PPDiv}(L/k)$ such that $\SAut(\D)$ exists. Let $\varrho:\Gamma_{L}\to \SAut(\D)$ be a Galois semilinear action. By \cref{proposition: semilinear functoriality}, it induces a Galois semilinear equivariant action $X(\varrho):\Gamma_{L}\to\SAut(T,X(\D))$ through the functor $X:\mathfrak{PPDiv}(L/k)\to\mathcal{E}(L/k)$. Hence, by \cref{proposition equivalence normal G-varieties}, the pair $(X(\D),X(\varrho))$ corresponds to an affine normal variety endowed with an effective action of an algebraic torus over $k$ that splits over $L$; this is a variety over $k$.

Conversely, let $X$ be an affine normal variety endowed with an effective action of an algebraic torus $T$ over $k$ that splits over $L$. Hence, $(X_{L},\varrho_{\mathrm{can}})$ is an affine normal $T_{L}$-variety over $L$. By \cref{proposition: theoremmainMN25}, there exists a pp-divisor $\D$ over $L$ such that $X_{L}\cong X(\D)$. Nevertheless, the Galois semilinear equivariant action might not translate to information over the pp-divisor. However, if such a pp-divisor is minimal (see \cite[Section 4.4.3]{MN25}), it does by \cite[Theorem 5.11]{MN25}. By \cite[Proposition 4.30]{MN25}, there always exists a minimal pp-divisor encoding a given affine normal $T$-variety. We thus have the following result.

\begin{proposition}\cite[Theorem 1.3]{MN25}\label{proposition: maintheoremofpaper}
Let $k$ be a field and $L/k$ be a finite Galois extension with Galois group $\Gamma_{L}$.
\begin{enumerate}[a)]
\item\label{proposition main affine minimal part i} Let $(\D_{L},\varrho)$ be a pair consisting of a pp-divisor and a Galois semilinear action on it. Then, $X(\D_{L},\varrho)$ is an affine normal variety endowed with an effective action of an algebraic torus $T$ over $k$ such that $T$ splits over $L$ and $X(\D_{L},\varrho)_{L}\cong X(\D_{L})$ as $T_{\D_{L}}$-varieties over $L$.
\item\label{proposition main affine minimal part ii} Let $X$ be an affine normal variety over $k$ endowed with an effective $T$-action such that $T_{L}$ is split. Then, there exists a pair $(\D_{L},\varrho)$ as above such that $X \cong X(\D_{L},\varrho)$ as $T$-varieties.
\end{enumerate}
\end{proposition}

This is precisely the crucial point that makes the description of non-affine normal varieties endowed with an effective action of a non-split algebraic torus different and more difficult, since there is no notion of \emph{minimality} for \emph{divisorial fans} (cf. \cref{definition divisorial fan}). In order to solve this problem, we first deal with the affine case for non-minimal pp-divisors in \cref{Appendix localization} and then address the general problem in \cref{section: normal varieties and non-split tori}.

\section{Normal varieties and split tori}\label{section divisorial fans}
\noindent

When $k=\sep{k}$ and $\mathrm{char}(k)=0$, every normal $T$ variety arises from a divisorial fan $(\S,Y)$ by \cref{theoremmainalthausus}. The goal of this section is to prove that this holds too when $k$ is any field and $T$ is a split $k$-torus. 
 
 \begin{theorem}\label{theoremahssplitcase}
Let $T$ be a split $k$-torus. Up to isomorphism, every normal $k$-variety with an effective $T$-action arises from a divisorial fan $(\S,Y)$ over $k$.
\end{theorem}

The first part of this section, we recall some definitions on pp-divisors. Hence, to give a characterization for when the morphism $(\id_{Y}, \id_{N}, \mathfrak{1}): \D’ \to \D$ induces a $T$-equivariant open embedding.

\subsection{More on pp-divisors.}
\noindent

A pp-divisor $\D\in\PPDiv_{\Q}(Y,\omega)$ defines a map $\mathfrak{h}_{\D}:\omega^{\vee}\to \cadiv_{\Q}(Y)$ given by $\mathfrak{h}_{\D}(m):=\D(m)$. This map satisfies certain properties summarized in the following definition.

\begin{definition}
    Let $Y$ be a normal $k$-variety; let $M$ be a lattice, and let $\omega^{\vee}\subset M_{\Q}$ be a cone of full dimension. We say that a map \[\mathfrak{h}:\omega^{\vee}\to \cadiv_{\Q}(Y)\] is 
    \begin{enumerate}[i)]
        \item \textit{convex} if $\mathfrak{h}(m)+\mathfrak{h}(m')\leq \mathfrak{h}(m+m')$ holds for any two elements $m,m'\in\omega^{\vee}$,
        \item \textit{piecewise linear} if there is a quasifan $\Lambda$ in $M_{\Q}$ having $\omega^{\vee}$ as its support such that $\mathfrak{h}$ is linear on the cones of $\Lambda$,
        \item \textit{strictly semiample} if $\mathfrak{h}(m)$ is semiample for all $m\in\omega^{\vee}$ and big for all $m\in\mathrm{relint}(\omega^{\vee})$.
    \end{enumerate}
    The set of all convex, piecewise linear and strictly semiample maps $\mathfrak{h}:\omega^{\vee}\to \cadiv_{\Q}(Y)$ is denoted by $\mathrm{CPL}_{\Q}(Y,\omega)$. 
\end{definition}

To each $\D\in\PPDiv_{\Q}(Y,\omega)$ we can associate a convex, piecewise linear and strictly semiample map $\mathfrak{h}_{\D}\in\mathrm{CPL}_{\Q}(Y,\omega)$. Thus, we have a natural map 
\begin{align*}
        \PPDiv_{\Q}(Y,\omega) &\to\mathrm{CPL}_{\Q}(Y,\omega), \\ \D &\mapsto \mathfrak{h}_{\D}.
    \end{align*}

The following proposition, which corresponds to \cite[Proposition 2.11]{AH06} and \cite[Proposition 3.6]{MN25}, helps us to prove one of the main theorems of \cref{Section Open embeddings}.

\begin{proposition}\label{proposition ppdiv cplmaps split}
Let $Y$ be a normal $k$-variety, $N$ be a lattice, and $\omega\subset N_{\Q}$ be a pointed cone. Then, the set $\mathrm{CPL}_{\Q}(Y,\omega)$ is a mono\"id and the canonical map
    $\PPDiv_{\Q}(Y,\omega) \to\mathrm{CPL}_{\Q}(Y,\omega)$ given by $ \D \mapsto \mathfrak{h}_{\D}$ is an isomorphism. Moreover, the integral polyhedral divisors correspond to maps $\mathfrak{h}:\omega^{\vee}\to \cadiv_{\Q}(Y)$ such that $\mathfrak{h}(\omega^{\vee}\cap M)\subset\cadiv(Y)$.
\end{proposition}

The following definitions concern \cref{Section Divisorial fan}, where the definition of \emph{divisorial fan} is given.

\begin{definition}
Let $Y$ be a normal variety over $k$. Let $\D\in\PPDiv_{\Q}(Y,\omega)$ and $\D'\in\PPDiv_{\Q}(Y,\omega')$ be pp-divisors.
    \begin{enumerate}[i)]
        \item We define the intersection of $\D$ and $\D'$ as 
        \[\D\cap\D':=\sum (\Delta_{D}\cap\Delta_{D}')\otimes D.\]
        \item For $y\in Y$, we define the \textit{fiber polyhedron at} $y$ as 
        \[\Delta_{y}:=\sum_{y\in D}\Delta_{D}.\]
    \end{enumerate}
\end{definition}

Let $m\in\omega^{\vee}\cap M$ and $f\in \H^{0}(Y,\mathscr{O}(\D(m)))$. Denote $Z(f)=\mathrm{Supp}(\mathrm{div}(f)+\D(m))$ its zero set and $Y_{f}=Y\smallsetminus Z(f)$ its principal open subset. If $m\in\mathrm{relint}(\omega^{\vee})$, $Y_{f}$ is an affine open subvariety of $Y$. Then, the Cartier divisors $D$ appearing in $\D$ induce Cartier divisors $D|_{Y_{f}}$ of $Y_{f}$ once we restrict them. Thus, we can define the following pp-divisor in $\PPDiv_{\Q}(Y_{f},\omega\cap m^{\perp})$ by taking the faces defined by $m$ of the polyhedra $\Delta_{D}$:
\[\D':=\sum\textrm{face}(\Delta_{D},m)\otimes D|_{Y_{f}} .\]

 We would like to have both pp-divisors over the same base to define a nice concept of a \textit{face} (cf. \cref{definitionface}) in the context of pp-divisors. Then, we will accept pp-divisors with $\emptyset$ coefficients satisfying the following conditions:
\[\emptyset+\Delta:=\emptyset\textrm{ and }0\cdot\emptyset:=\omega.\] 
Thus, we will always assume that $\bigcup_{\Delta_{D}=\emptyset}\mathrm{Supp}(D)$ is the support of a semiample divisor, and by a pp-divisor $\D\in\PPDiv_{\Q}(Y,\omega)$ we mean $\D_{|\mathrm{Loc(\D)}}\in\PPDiv_{\Q}(\mathrm{Loc(\D)},\omega)$, where $\mathrm{Loc}(\D)$ is defined here below.

\begin{definition}
    Let $k$ be a field. Let $Y$ be a normal $k$-variety. Let $\D\in\PPDiv_{\Q}(Y,\omega)$ and $m\in\mathrm{relint}(\omega^{\vee})\cap M$.
    \begin{enumerate}[i)]
        \item We define the \textit{locus of} $\D$ as 
        \[\mathrm{Loc}(\D):=Y\smallsetminus \bigcup_{\Delta_{i}=\emptyset} \mathrm{Supp}(D_{i}).\]
        \item The \textit{localization} of $\D$ by $f\in \H^{0}(Y,\mathscr{O}(\D(m)))$ is 
        \[\D_{f}:=\sum\textrm{face}(\Delta_{D},m)\otimes D|_{Y_{f}}=\emptyset\otimes(\mathrm{div}(f)+\D(m))+\sum\mathrm{face}(\Delta_{D},m)\otimes D.\]
    \end{enumerate}
\end{definition}

\subsection{Open embeddings}\label{Section Open embeddings}
\noindent

Let $Y$ be a normal semi-projective variety over $k$. Recall that, for any pair of pp-divisors $\D=\sum_D \Delta_{D}\otimes D$ and $\D’=\sum_D \Delta’_{D}\otimes D$ on $Y$ whose tail cones are defined over the same lattice $N$, we have that $\D \leq \D’$ if and only if $\Delta_{D}' \subset \Delta_{D}$ for every $D$. This relation encodes a morphism of pp-divisors $(\id_{Y}, \id_{N}, \mathfrak{1}): \D’ \to \D$, which induces a $T$-equivariant morphism $X(\id_{Y}, \id_{N}, \mathfrak{1}): X(\D’) \to X(\D)$. In the toric case, for cones $\omega’$ and $\omega$, such a morphism corresponds to $X_{\omega’} \to X_{\omega}$ induced by the inclusion $\omega’ \subset \omega$, and we know that $X_{\omega’} \to X_{\omega}$ is an open embedding if and only if $\omega’ \preceq \omega$, i.e. $\omega’$ is a face of $\omega$. In such a context, the $T$-equivariant open embeddings correspond to localizations. However, for general $T$-varieties, faces do not come from a single localization.

The main results of this section are \cref{Proposition 3.4 AHS08}, \cref{Proposition 4.3 AHS08} and \cref{Corollary 4.4 AHS08}.  The following lemma corresponds to \cite[Lemma 3.2]{AHS08}, whose proof adapts to this frame work.
\begin{lemma}\label{Lemma 3.2 AHS08}
Let $Y$ be a normal variety over $k$, and $\omega$ be a cone in a lattice $N$. Let $\D$ be a pp-divisor in $\PPDiv_{\Q}(Y,\omega)$. Let $m\in\omega^{\vee}\cap M$ and $f\in A_{m}:=\H^{0}(Y,\mathscr{O}_{Y}(\D(m)))$. Then, for every $m'\in M$, there exists $k\gg 0$ such that $m'+km\in \omega^{\vee}$. Moreover, if $m'\in \left(\omega\cap m^{\perp}\right)^{\vee}\cap M$, then 
\[\D_{f}(m')=\D(m'+km)|_{Y_{f}}-\D(km)|_{Y_{f}}.\]
\end{lemma}

\begin{proposition}\label{Proposition 3.3 AHS08}
Let $\D$ be an object in $\mathfrak{PPDiv}(k)$ and $\D_{f}$ be a localization by an element $f$ of $\H^{0}(Y,\mathscr{O}_{Y}(\D(m)))$. Let $\iota:Y_{f}\to Y$ be the canonical open embedding. Then, $\D_{f}$ is a pp-divisor in $\PPDiv_{\Q}(Y_{f},\omega\cap m^{\perp})$ and $X(\D_{f})\cong X(\D)_{f}$. Also, the morphism of pp-divisors $(\iota, \id_{N}, \mathfrak{1}):\D_{f}\to\D$ induces a $T$-equivariant open embedding $X(\iota, \id_{N}, \mathfrak{1}):X(\D_{f})\to X(\D)$.

\end{proposition}

\begin{proof}
Set $\tau:=\omega\cap m^{\perp}$. By \cref{Lemma 3.2 AHS08}, we have that $\D_{f}:M\cap\tau^{\vee}\to\cadiv_{\Q}(Y_{f},\tau)$ is a convex piecewise linear strictly semiample map. Then, by \cref{proposition ppdiv cplmaps split}, it arises by a pp-divisor which coincides with the polyhedral divisor $\D_{f}$. Thus, $\D_{f}$ is a pp-divisor.

The morphism $(\iota, \id_{N}, \mathfrak{1}):\D_{f}\to\D$ induces a morphism of algebras
\begin{align*}
	A[Y,\D]=\bigoplus_{m'\in\omega^{\vee}\cap M}\H^{0}\left( Y,\mathscr{O}_{Y}(\D(m'))\right)&\to\bigoplus_{m'\in\tau^{\vee}\cap M}\H^{0}\left( Y_{f},\mathscr{O}_{Y_{f}}(\D_{f}(m'))\right)=A[Y_{f},\D_{f}] \\
	g &\mapsto \iota^{*}g,
\end{align*}
which is injective since $\iota^{*}g$ is the restriction of g to $Y_{f}$. Given that $f$ has no zeros in $Y_{f}$, we have that $\iota^{*}f$ is invertible in $A[Y_{f},\D_{f}]$. Hence, the morphism $\iota^{*}:A[Y,\D]\to A[Y_{f},\D_{f}]$ factorizes through the localization
\[\xymatrix{ A[Y,\D] \ar[r]^{\iota^{*}} \ar[rd]_{} & A[Y_{f},\D_{f}] \\ & A[Y,\D]_{f} \ar[u]_{\alpha} .}\]
Since $\iota^{*}$ is injective, the morphism $\alpha$ is injective. Given that $f$ is homogeneous, the morphism $\alpha:A[Y,\D]_{f}\to A[Y_{f},\D_{f}]$ is a morphism of graded algebras. Let us prove now that $\alpha$ is surjective. Let $h\in \H^{0}(Y_{f},\mathscr{O}_{Y_{f}}(\D_{f}(m')))$, then $\mathrm{div}(h)+\D_{f}(m')\geq 0$ and $h\in k(Y)$. By \cref{Lemma 3.2 AHS08},
\[\mathrm{div}(h)+\D_{f}(m')=\mathrm{div}(h)+\D(m'+km)|_{Y_{f}}-\D(km)|_{Y_{f}}\geq 0.\]
Since $Y_{f}=Y\smallsetminus Z(f)$ and $Z(f)=\mathrm{Supp}(\mathrm{div}(f)+\D(m))$, for some $l\in \N$, we have that
\[\mathrm{div}(h)+\D(m'+km)-\D(km)+l\mathrm{div}(f)+l\D(km)\geq 0.\]
Otherwise stated,
\[\mathrm{div}(hf^{l})+\D(m'+km)+\D((lk-k)m)\geq 0.\] 
Given that $\D(m'+lkm)\geq \D(m'+km)+\D((lk-k)m)$, we have that 
\[\mathrm{div}(hf^{l})+\D(m'+lkm)\geq 0\]
and, therefore, $hf^{l}\in \H^{0}(Y,\D(m'+lkm))$. Hence, in the localization, $\alpha(hf^{l}/f^{l})=h$. This proves that $\alpha$ is surjective and, therefore, $A[Y,\D]_{f}\cong A[Y_{f},\D_{f}]$. This implies that $X(\D_{f})\cong X(\D)_{f}$ and $X(\iota, \id_{N}, \mathfrak{1}):X(\D_{f})\to X(\D)$ is an open embedding.
\end{proof}

\begin{proposition}\label{Proposition 3.4 AHS08}
Let $k$ be a field. Let $\D$ and $\D'$ be two pp-divisors over the same semi-projective variety over $k$ such that $\D\leq \D'$. Let $(\id_{Y}, \id_{N}, \mathfrak{1}):\D'\to\D$ be the corresponding morphism of pp-divisors. Then, the morphism $X(\id_{Y}, \id_{N}, \mathfrak{1}):X(\D')\to X(\D)$ is an open embedding if and only if there exist $m_{1},\dots,m_{r}\in\omega^{\vee}\cap M$ and $f_{i}\in \H^{0}(Y,\mathscr{O}_{Y}(\D(m_{i})))$ such that $\bigcup \Loc(\D)_{f_{i}}=\Loc(\D')$ and $\face(\Delta_{D}',m_{i})=\face(\Delta_{D},m_{i})$ for every $D\not\subset Z(f_{i})$ and $i\in\{1,\dots,r\}$.
\end{proposition}

\begin{proof}
Let us suppose that $\Loc(\D)=Y$ and denote $Y':=\Loc(\D')$. Also, we denote by $\Psi:A[Y,\D]\to A[Y',\D']$ the respective morphism of graded algebras. 

Let us suppose that $X(\id_{Y}, \id_{N}, \mathfrak{1}):X(\D')\to X(\D)$ is an open embedding, then $\Psi$ is injective. The open embedding condition is equivalent to the following one: there exist $f_{1},\dots,f_{r}\in A[Y,\D]$ such that $1\in (\Psi(f_{1}),\dots,\Psi(f_{r}))$ and $A[Y,\D]_{f_{i}}\cong A[Y',\D']_{\Psi(f_{i})}$ for every $i\in\{1,\dots,r\}$. Moreover, each $f_{i}$ can be considered homogeneous of degree $m_{i}\in\omega^{\vee}\cap M$ and, therefore, $\Psi(f_{i})$ is also homogeneous of degree $m_{i}$. Denote $\D'_{\Psi(f_{i})}$ and $\D_{f_{i}}$ the respective localizations. By \cref{Proposition 3.3 AHS08}, we have that $\D'_{\Psi(f_{i})}$ and $\D_{f_{i}}$ are pp-divisors over $Y'_{\Psi(f_{i})}$ and $Y_{f_{i}}$ respectively and 
\[A\left[Y_{f_{i}},\D_{f_{i}}\right]\cong A\left[Y,\D\right]_{f_{i}}\cong A\left[Y',\D'\right]_{\Psi(f_{i})}\cong A\left[Y'_{\Psi(f_{i})},\D'_{\Psi(f_{i})}\right].\]
Such isomorphism arises from $(\iota_{i},\id_{N},\mathfrak{1}):\D'_{\Psi(f_{i})}\to\D_{f_{i}}$, where $\iota_{i}:Y'_{\Psi(f_{i})}\to Y_{f_{i}}$ is the canonical inclusion. Moreover, $\iota_{i}:Y'_{\Psi(f_{i})}\to Y_{f_{i}}$ is the identity map. Hence, the pp-divisors are equal $\D'_{f_{i}}=\D_{f_{i}}$ for every $f_{i}$. Given that $1\in (\Psi(f_{1}),\dots,\Psi(f_{r}))$, we have that 
\[\bigcap Z(f_{i})=\bigcap Z(\Psi(f_{i}))=\emptyset.\]
 Thus, $Y'=\bigcup \Loc(\D)_{f_{i}}$ and $\face(\Delta'_{D},m_{i})=\face(\Delta_{D},m_{i})$, when $D\not\subset Z(f_{i})$.

If there exist $m_{1},\dots,m_{r}\in\omega^{\vee}\cap M$ and $f_{i}\in \H^{0}(Y,\mathscr{O}_{Y}(\D(m_{i})))$ such that $\bigcup \Loc(\D)_{f_{i}}=\Loc(\D')$ and $\face(\Delta_{D}',m_{i})=\face(\Delta_{D},m_{i})$ for every $D\not\subset Z(f_{i})$ and $i\in\{1,\dots,r\}$, we have that $\D'_{f_{i}}=\D_{f_{i}}$ for every $f_{i}$. Hence, by \cref{Proposition 3.3 AHS08}, the maps $X(\D'_{f_{i}})\to X(\D)$ are open embeddings. Given that $\bigcup \Loc(\D)_{f_{i}}=\Loc(\D')$, we have that $X(\D')$ is covered by $X(\D_{f_{i}})$ and therefore that $X(\D')\to X(\D)$ is an open embedding.
\end{proof}

Let $(\psi,F,\mathfrak{f}):\D'\to \D$ be a morphism of pp-divisors such that the respective $T$-equivariant morphism $X(\psi,F,\mathfrak{f}):X(\D')\to X(\D)$ is an open embedding. It also induces a $T$-equivariant morphism $\tilde{X}(\D')\to \tilde{X}(\D)$. These morphisms fit into the following commutative diagram
\[
\xymatrixcolsep{5pc}\xymatrix{ Y' \ar[d]_{\psi} & \tilde{X}(\D') \ar[d]_{\tilde{\varphi}:=\tilde{X}(\psi,F,\mathfrak{f})} \ar[l]_{\pi'} \ar[r]^{r'} & X(\D') \ar[d]^{\varphi:=X(\psi,F,\mathfrak{f})} \\ Y & \tilde{X}(\D) \ar[l]^{\pi} \ar[r]_{r} & X(\D)\,. }
\]
The following result corresponds to \cite[Lemma 4.5]{AHS08}, whose proof adapts to this general context.

\begin{proposition}\label{Lemma 4.5 AHS08}
In the above notation, suppose that the morphism $\varphi: X(\D') \to X(\D)$ is an open embedding. Then the following statements hold:
\begin{enumerate}[(i)]
    \item We have $\tilde{\varphi}(\tilde{X}(\D')) = r^{-1}(\varphi(X(\D)))$, and the induced morphism 
    \[
    \tilde{\varphi}: \tilde{X}(\D') \to \tilde{\varphi}(\tilde{X}(\D'))
    \]
    is proper and birational.

    \item The image $\psi(Y') \subset Y$ is open and semi-projective, and $\psi: Y' \to \psi(Y')$ is a projective birational morphism.

    \item For every $y \in \psi(Y')$, the intersection
    \[
    \tilde{U}_y := \pi^{-1}(y) \cap \tilde{\varphi}(\tilde{X}(\D'))
    \]
    contains a unique $T$-orbit that is closed in $\tilde{U}_y$.
\end{enumerate}
\end{proposition}

Let $Y$ be a normal semi-projective variety over $k$. Let $\D=\sum\Delta_{D}\otimes D$ be a pp-divisor over $Y$. By \cite[Remark 4.18]{MN25}, the $T$-varieties $X:=X(\D)$ and $\tilde{X}:=\tilde{X}(\D)$ fit into the following commutative diagram 
\[\xymatrix{ \tilde{X} \ar[r]^{r} \ar[d]_{\pi} & X \ar[d] \\ Y \ar[r] & Y_{0},}\] 
where $\pi$ is affine and $r$ is a proper birational $T$-equivariant map. 

In the following we will see how to attach a pp-divisor for any $T$-stable affine open subvariety $\Spec(A)=X'\subset X(\D)$ defined over $Y$. Recall that the canonical inclusion $X'\to X$ is an open embedding if there exist $f_{1},\dots,f_{r}\in A[Y,\D]$, that can be considered homogeneous, such that $(f_{1},\dots,f_{r})=A$ and $A[Y,\D]_{f_{i}}=A_{f_{i}}$.

\begin{proposition}\label{Proposition 4.3 AHS08}
Let $Y$ be a normal semi-projective variety over $k$. Let $\D=\sum\Delta_{D}\otimes D$ be a pp-divisor over $Y$. Let $\tilde{X}:=\tilde{X}(\D)$ and $X:=X(\D)$. Let $X'\subset X$ be a $T$-stable affine open subvariety. Then, $Y':=\pi\left(r^{-1}(X')\right)$ is an open subvariety of $Y$ and it is semi-projective. Moreover, if $X'=\Spec (A)$ and the regular functions $f_{i}\in \H^{0}(Y,\mathscr{O}_{Y}(\D(m_{i})))$ are as above, then 
\[\D':=\bigcup\D_{f_{i}}:=\sum \Delta_{D}'\otimes D|_{Y'},\]
where $\Delta_{D}':=\bigcup_{D\cap Y_{f_{i}}'\neq\emptyset}\face\left(\Delta_{D},m_{i}\right)\subset \Delta_{D}$, is a pp-divisor on $Y'=\bigcup Y_{f_{i}}$ and the canonical map $\D'\to \D$ defines an open embedding of affine $T$-varieties and $X(\D')\cong X'$ under such a morphism.
\end{proposition}

\begin{proof}
First, the proofs of both \cite[Lemma 3.1]{Hau09} and \cite[Lemma 4.5]{AHS08} adapts to this context. Thus, $Y'\subset Y$ is an open subvariety and it is semi-projective.

On the one hand, given that $\face\left( \Delta_{D},m_{i} \right)\subset\Delta_{D}'\subset \Delta_{D}$, we have that $\face\left( \Delta_{D},m_{i} \right)\subset \face\left( \Delta_{D}',m_{i} \right)$, by \cref{lemma polyhedron face}~\eqref{part a lemma polyhedron face}. On the other hand, by \cref{lemma polyhedron face}~\eqref{part b lemma polyhedron face}, we have that 
\begin{equation*}
\resizebox{.98\hsize}{!}{
$\face\left( \Delta_{D}',m_{i} \right)\subset \bigcup\face\left( \face\left( \Delta_{D},m_{j} \right),m_{i} \right)= \bigcup\face\left( \face\left( \Delta_{D},m_{i} \right),m_{j} \right)=\bigcup\face\left(\Delta_{D},m_{i} \right)$
}
\end{equation*}
and, therefore, $\face\left( \Delta_{D}',m_{i} \right)\subset\face\left( \Delta_{D},m_{i} \right) $. Hence, $\face\left( \Delta_{D}',m_{i} \right)=\face\left( \Delta_{D},m_{i} \right)$. Since by construction $\Loc(\D')=\bigcup\Loc(\D_{f_{i}})$, in order to proof that the respective morphism of pp-divisors induces an open embedding it suffices to prove that $\D'$ is a pp-divisor and concludes by \cref{Proposition 3.4 AHS08}. It is worth mentioning that, as stated in the proof of \cite[Proposition 4.3]{AHS08}, the tail cone of the polyedron $\face\left( \Delta_{D},m_{i} \right)$ is the same for every $m_{i}$.

Assume that $\D$ is a minimal pp-divisor. By \cite[Proposition 4.30]{MN25}, there exists a minimal pp-divisor $\D''$ such that $X(\D'')\cong X'\subset X$. Hence, by \cite[Theorem 5.11]{MN25}, the open embedding $X'\subset X$ arises from a morphism $(\psi,\id_{N},\mathfrak{f}):\D''\to \D$. Such a morphism induces morphisms $(\psi_{i},\id_{N},\mathfrak{f}_{i}):\D''_{\psi^{*}(f_{i})}\to \D_{f_{i}}$ fitting in the following commutative diagram
\[
\xymatrixcolsep{5pc}\xymatrix{ \D'' \ar[r]^{(\psi,\id_{N},\mathfrak{f})}& \D \\ 
\D''_{\psi^{*}(f_{i})} \ar[u]^{(\iota''_{i},\id_{N},\mathfrak{1})} \ar[r]_{(\psi_{i},\id_{N},\mathfrak{f}_{i})} & \D_{f_{i}} \ar[u]_{(\iota_{i},\id_{N},\mathfrak{1})} \,.}
\]
We have then a morphism $(\psi,\id_{N},\mathfrak{f}):\D''\to \D$ factorizes through $(\id_{Y},\id_{N},\mathfrak{1}):\D'\to\D$ and $(\psi,\id_{N},\mathfrak{f}):\D''\to \D'$. This implies that $\D''$ and $\psi^{*}\D'$ are isomorphic and, therefore, $\psi^{*}\D'$ is a proper and, therefore, a pp-divisor. Consequently, we have that $\D'(m)+\D'(m')\leq \D'(m+m')$ for any pair $m$ and $m'$ in $M$. Hence, $A[Y',\D']$ is a $M$-graded $k$-algebra. Now, since $(\psi',\id_{N},\mathfrak{1}):\psi'^{*}\D'\to \D'$ induces an isomorphism of $M$-graded $k$-algebras, we conclude that $\D'$ is a pp-divisor.

In the case when $\D$ is not a minimal pp-divisor, you can take a minimal pp-divisor $\tilde{\D}$ encoding the same affine normal $T$-variety. Hence, proceed as before and then pullback through $(\vartheta,\id_{N},\tilde{\mathfrak{f}}):\D\to\tilde{\D}$ (see \cite[Lemma 5.10]{MN25}).  

\end{proof}

Then we have the following corollary that will help us later in the construction of a \emph{divisorial fan} given a certain normal $T$-variety (see \cref{Section Divisorial fan}).

\begin{corollary}\label{Corollary 4.4 AHS08}
Let $Y$ be a normal semi-projective over $k$. Let $\D$ be a pp-divisor on $Y$ and $X:=X(\D)$ be the respective $T$-variety. Let $X'$ and $X''$ be two $T$-stable affine open subvarieties of $X$ arising from pp-divisors 
\[\D'=\bigcup\D_{f_{i}}=\sum \Delta_{D}'\otimes D|_{Y'}\quad\textrm{ and }\quad\D''=\bigcup\D_{g_{i}}=\sum \Delta_{D}''\otimes D|_{Y''},\]
with loci $Y'$ and $Y''$, respectively. Then, 
\[\D'\cap\D'':=\bigcup\D_{f_{i}g_{j}}=\sum \left(\Delta_{D}'\cap \Delta_{D}''\right)\otimes D|_{Y'\cap Y''}\]
is a pp-divisor over $Y'\cap Y''$ and the canonical map $\D'\cap\D''\to\D$ induces a $T$-equivariant open embedding $X(\D'\cap\D'')\to X(\D)$ having $X'\cap X''$ as its image.
\end{corollary}

\subsection{Divisorial fans}\label{Section Divisorial fan}
\noindent

As stated in \cite{AH06} and \cite{MN25}, affine normal $T$-varieties arise from pp-divisors. In \cite{AHS08}, the authors introduced the notion of \emph{divisorial fan}, which is a generalization of a fan in the classical theory of toric varieties. In the following, we recall the definition of \emph{face of a pp-divisor} and the definition of a \emph{divisorial fan}.

\begin{definition}\label{definitionface}
   Let $k$ be a field. Let $N$ be a lattice, $\omega$ and $\omega'$ be cones in $N_{\Q}$, $Y$ be a normal variety over $k$, and consider two pp-divisors on $Y$: \[\D'=\sum \Delta_{D}'\otimes D\in \PPDiv_{\Q}(Y,\omega'),\quad \D=\sum \Delta_{D}\otimes D\in \PPDiv_{\Q}(Y,\omega). \] We call $\D'$ a \textit{face} of $\D$ (written $\D'\preceq\D$) if $\Delta_{D}'\subset \Delta_{D}$ holds for all $D$ and there exist $m_{1},\dots,m_{l}\in\omega^{\vee}\cap M$ and $D_{i}$ in the linear system $\lvert \D(m_{i}) \rvert$, for every $m_{i}$, such that 
	\begin{enumerate}[i)]
		\item $\Loc(\D')=\bigcup \left(\Loc(\D)\smallsetminus \mathrm{Supp}\left(D_{i}\right)\right)$ and
		\item $\face(\Delta_{D}',m_{i})=\face(\Delta_{D},m_{i})$ for every $D\not\subset \mathrm{Supp}(D_{i})$ and $i\in\{1,\dots,l\}$.
	\end{enumerate}   
\end{definition}

Recall that the relation $\D'\leq \D$ between two pp-divisors means that $\Delta_{D}\subset \Delta_{D}'$ for every prime divisor $D\in\cadiv(Y)$. Then, if $\D'\preceq\D$ we have $\Delta'_{D}\subset \Delta_{D}$ and, therefore, $\D\leq \D'$. Hence, for any pair of pp-divisors $\D'$ and $\D$ in $\mathfrak{PPDiv}(k)$ such that $\D'\preceq\D$ we have that the triple $(\id,\id,1)$ defines a morphism of pp-divisors $\D'\to\D$.

Let $\D'$ and $\D$ in $\mathfrak{PPDiv}(k)$ be such that $\D'\preceq\D$ and let $(\id,\id,1):\D'\to\D$ be its corresponding morphism of pp-divisors.  Given that $\D\leq\D'$, we have that $\D(m)\leq \D'(m)$ for every $m\in\omega^{\vee}\cap M$. Then, $\H^{0}(Y,\mathscr{O}(\D(m)))\subset \H^{0}(Y,\mathscr{O}(\D'(m)))$ for every $m\in\omega'^{\vee}\cap M$. These inclusions induce an inclusion of $M$-graded $k$-algebras
\[A[Y,\D]=\bigoplus_{m\in\omega^{\vee}\cap M} \H^{0}(Y,\mathscr{O}(\D(m)))\subset \bigoplus_{m\in\omega'^{\vee}\cap M} \H^{0}(Y,\mathscr{O}(\D'(m))) = A[Y,\D'], \]
which induces a $T$-equivariant dominant morphism 
\[ X(\D'):=\Spec(A[Y,\D'])\subset \Spec(A[Y,\D])=:X(\D) \]
between the respective affine normal $T$-varieties over $k$. Thus, from a morphism of pp-divisors of the form $(\id,\id,1):\D'\to\D$ we obtain a $T$-equivariant dominant morphism $\iota:X(\D')\to X(\D)$. We aim to determine when such a morphism of pp-divisors encodes an $T$-equivariant open immersions. Actually, the $T$-equivariant dominant morphism $\iota:X(\D')\to X(\D)$ that arise from a face relation are exactly those that are open immersions.

\begin{proposition}\label{corollaryfaceopenimmersion}
    Let $k$ be a field. Let $\D$ and $\D'$ be two pp-divisors such that $\D'\leq \D$. Then $\D'\preceq \D$ if and only if $\iota: X(\D')\to X(\D)$ is an open immersion.
\end{proposition}

\begin{proof}
    Let us suppose that $\D'\preceq \D$. By definition, there exist $m_{1},\dots,m_{l}\in\omega^{\vee}\cap M$ and respective $D_{i}\in \left| \D(m_{i})\right|$ such that $\Loc(\D')=\bigcup \left(\Loc(\D)\smallsetminus \mathrm{Supp}\left(D_{i}\right)\right)$ and $\face(\Delta_{D}',m_{i})=\face(\Delta_{D},m_{i})$ for every $D\not\subset \mathrm{Supp}(D_{i})$ and $i\in\{1,\dots,r\}$. Then, there exist $f_{i}\in \H^{0}(Y,\mathscr{O}_{Y}(\D(m_{i})))$ such that $D_{i}=\mathrm{div}\left(f_{i}\right)+\D(m_{i})$. Given that $\D'(m_{i})\leq \D(m_{i})$, we have that $f_{i}\in \H^{0}(Y,\mathscr{O}_{Y}(\D(m_{i})))\subset \H^{0}(Y,\mathscr{O}_{Y}(\D'(m_{i})))$. Thus, we have that $\Loc(\D)_{f_{i}}\subset \Loc(\D')$ and, moreover, $\Loc(\D')=\bigcup\Loc(\D)_{f_{i}}$. Also, we have that $\face(\Delta_{D}',m_{i})=\face(\Delta_{D},m_{i})$ for every $D$ satisfying
    \[  D\not\subset \mathrm{Supp}(D_{i}) = Z(f_{i}). \]
Then, by \cref{Proposition 3.4 AHS08}, we have that the canonical map $X(\D')\to X(\D)$ is an open embedding.

On the other direction, by \cref{Proposition 3.4 AHS08}, there exist $m_{1},\dots,m_{r}\in\omega^{\vee}\cap M$ and $f_{i}\in \H^{0}(Y,\mathscr{O}_{Y}(\D(m_{i})))$ such that $\bigcup \Loc(\D)_{f_{i}}=\Loc(\D')$ and $\face(\Delta_{D}',m_{i})=\face(\Delta_{D},m_{i})$ for every $D\not\subset Z(f_{i})$ and $i\in\{1,\dots,r\}$. Hence, the divisors $D_{i}:=\mathrm{div}\left(f_{i}\right)$ satisfy the conditions of \cref{definitionface}.
\end{proof}

\begin{remark}
Let $\D\in\PPDiv_{\Q}(Y,\omega)$ and $\D'\in\PPDiv_{\Q}(Y',\omega')$ be pp-divisors and denote by $T'$ the split algebraic torus acting on $X(\D')$ associated to $\D'$. Assume that $X(\D')\subset X(\D)$ is a $T'$-equivariant open embedding. If $\D'$ and $\D$ are minimal pp-divisors, it does not imply that there is an open immersion $Y'\to Y$. For example, the affine space $\A_{k}^{3}$ with the action of $\G_{\mathrm{m},k}$ given by 
\[\lambda\cdot (x,y,z)=(\lambda x,\lambda y,\lambda^{-1} z)\]
is encoded by the pp-divisor 
\[\D:=\{1\}\otimes D_{(1,0)}+\{0\}\otimes D_{(0,1)}+[0,1]\otimes D_{(1,1)} \in \PPDiv_{\Q}(\mathrm{Bl}_{0}(\A_{k}^{2}),\omega)\]
and the $\G_{\mathrm{m},k}$-stable affine open subvariety $\A_{k}^{2}\times\G_{\mathrm{m},k}$ of $\A_{k}^{3}$ is encoded by the following minimal pp-divisor 
\[\D':=\{1\}\otimes \tilde{D}_{(1,0)}+\{0\}\otimes \tilde{D}_{(0,1)} \in \PPDiv_{\Q}(\A_{k}^{2},\omega).\]

Let $\pi: \mathrm{Bl}_{0}(\A_{k}^{2}) \to \A_{k}^{2}$ be the canonical map. The pp-divisor $\D'$ is not a face of $\D$ because they do not \emph{live} over the same semi-projective variety ($\mathrm{Blow}_{0}(\A_{k}^{2})$ is not affine). However, after pulling back $\D'$ by $\pi$, the pp-divisor $\pi^{*}\D'$ is a face of $\D$ and $X(\D')\cong X(\pi^{*}(\D))$.
\end{remark}

\begin{remark}
    By \cite[Proposition 6.23]{MN25}, the face relation $\D'\preceq\D$ also induces a $T$-equivariant morphism $\tilde{\iota}:\tilde{X}(\D')\to\tilde{X}(\D)$ whose affinisation is $\iota$, that is, $\tilde{\iota}_{\mathrm{aff}}=\iota$. Moreover, $\tilde{\iota}$ is an open embedding by \cite[Remark 3.5(ii)]{AHS08}.
\end{remark} 

The two following results implicitly use the fact that for a pp-divisor $\D\in \PPDiv_{\Q}(Y,\omega)$ over $k$ there is a canonic isomorphism $X(\D)_{L}\cong X(\D_{L})$. This holds since for all $m\in \omega^{\vee}\cap M$ we have that $\H^{0}(Y,\mathscr{O}_{Y}(\D(m)))_{L}\cong \H^{0}(Y_{\sep{k}},\mathscr{O}_{Y_{L}}(\D_{L}(m)))$.

\begin{proposition}\label{Proposition face stable under base change}
Let $N$ be a lattice, $\omega$ and $\omega'$ in $N_{\Q}$ be pointed cones, and $Y$ be a normal semi-projective variety over $k$. Let $\D\in \PPDiv_{\Q}(Y,\omega)$ and $\D'\in \PPDiv_{\Q}(Y,\omega')$, then
    \[\D'\preceq\D \Leftrightarrow \D'_{L}\preceq \D_{L}.\]
\end{proposition}

\begin{proof}
   By \cref{Proposition 3.4 AHS08}, we have that $\D'\preceq \D$ if and only if $X(\D')\to X(\D)$ is an open immersion. Hence, by \cite[\href{https://stacks.math.columbia.edu/tag/02L3}{Tag 02L3}]{stacks-project}, $X(\D')\to X(\D)$ is an open immersion if and only if $X(\D_{L}')\to X(\D_{L})$ is an open immersion. Finally, again by \cref{Proposition 3.4 AHS08}, we have that $X(\D_{L}')\to X(\D_{L})$ is an open immersion if and only if $\D_{L}'\preceq \D_{L}$. Thus, $\D'\preceq \D$ if and only if $\D_{L}'\preceq \D_{L}$.
\end{proof}

\begin{definition}\label{definition divisorial fan}
    Let $N$ be a lattice and $Y$ a geometrically normal semi-projective $k$-variety. A \textit{divisorial fan} in $(Y,N)$ is a finite set $\S$ of pp-divisors $\D$ in $\mathfrak{PPDiv}_{N}(k)$ with base $Y$ such that, for any pair $\D,\D'\in \S$, the intersection $\D\cap\D'$ is in $\S$ and $\D\cap\D'\preceq\D,\D'$.
\end{definition}

\begin{proposition}\label{Proposition divisorial fan stable under base change}
Let $N$ be a lattice and $Y$ be a normal semi-projective variety over $k$. A finite set of pp-divisors $\S$ on $(Y,N)$ is a divisorial fan if and only if $\S_{L}$ is a divisorial fan on $(Y_{L},N)$.
\end{proposition}

\begin{proof}
This is a direct consequence of \cref{Proposition face stable under base change}.
\end{proof}

\begin{remark}
    Both \cref{Proposition face stable under base change,Proposition divisorial fan stable under base change} hold when $L=\sep{k}$ under the assumption of $Y$ being geometrically normal.
\end{remark}

\subsection{From divisorial fans to normal $T$-varieties}\label{section: from divisorial fans to normal T-varieties}
\noindent

We will see that \cref{theoremmainalthausus} holds for normal varieties endowed with an effective action of a split torus over any field.

Let $\S$ be a divisorial fan over a normal semi-projective $k$-variety $Y$ (we will use the pair $(\S,Y)$ to refer to this data). If $\D^{i}$ and $\D^{j}$ are in $\S$, then $\D^{i}\cap\D^{j}\in\S$ is a face of both. Therefore, by \cref{corollaryfaceopenimmersion}, we have $T$-equivariant open embeddings \[\xymatrix{X(\D^{i}) & X(\D^{i}\cap\D^{j}) \ar[l]_{\eta_{ij}} \ar[r]^{\eta_{ji}} & X(\D^{j}) }.\]

Let us set $X_{ij}:=\mathrm{im}(\eta_{ij})\subset X_{i}$ and $X_{ji}:=\mathrm{im}(\eta_{ji})\subset X_{j}$. Then, we have $T$-equivariant isomorphisms $\varphi_{ij}:=\eta_{ji}\circ\eta_{ij}^{-1}:X_{ij}\to X_{ji}$. When $k$ is algebraically closed, we construct the following quotient space by a relation over the points: 
\[X(\S):=\left(\bigsqcup_{\D\in\S} X(\D)\right)/\sim,\]
where the relation is given by $x\sim y$ if and only if for some $\varphi_{ij}$ we have $\varphi_{ij}(x)=y$.

Over nonalgebraically closed fields, we can also construct a scheme $X(\S)$. First, notice that a divisorial fan $(\S,Y)$ is a partially ordered set under the face relation, as in the toric case with the \emph{classical} fans. The $T$-equivariant open embeddings $\eta_{\E\D}:X(\E)\to X(\D)$, for $\E\preceq\D$ in $\S$, define a direct system satisfying the compatibility conditions: $\eta_{\mathfrak{F}\D}=\eta_{\mathfrak{F}\E}\circ\eta_{\E\D}$, for every $\mathfrak{F}\preceq\E\preceq\D$ in $\S$. Thus, the scheme $X(\S)$ can be defined as 
\[X(\S):=\varinjlim_{\D\in\S} X(\D).\]
Besides, when $Y$ is geometrically normal, over a separable closure $\sep{k}$, there is a canonic isomorphism between the two constructions.

Whereas in toric geometry the scheme $X(\Sigma)=\varinjlim_{\omega\in\Sigma} X_{\omega}$ is always separable, where $\Sigma$ is referred to a classical fan, the scheme $X(\S)$ might not be separable for a divisorial fan, but it is a prevariety (see for example: \cite[Example 7.9]{AHS08}).

\begin{proposition}
    Let $Y$ be a normal semi-projective variety over $k$. Let $(\S,Y)$ be a divisorial fan. The space $X(\S)$ is a prevariety over $k$ with affine diagonal $X(\S)\to X(\S)\times X(\S)$, and it comes with a (unique) $T$-action such that all the canonical maps $X(\D)\to X(\S)$ are $T$-equivariant.
\end{proposition}

\begin{proof}
By \cite[\href{https://stacks.math.columbia.edu/tag/01JB}{Tag 01JB}]{stacks-project} and \cite[\href{https://stacks.math.columbia.edu/tag/01JC}{Tag 01JC}]{stacks-project}, it is enough to prove that the maps $\eta_{ij}:X_{ij}\to X_{ji}$ define a gluing data and the proof follows straightforwardly as in the proof of \cite[Theorem 5.3]{AHS08}.
\end{proof}

Those divisorial fans defining separated schemes or varieties are called \emph{separated divisorial fans}.

\subsection{Separated and complete divisorial fans}\label{section: separated divisorial fans}
\noindent

Let $(\S, Y)$ be a divisorial fan over a field $k$, where $Y$ is a normal semi-projective $k$-variety. As said above, the scheme $X(\S)$ is not necessarily separated. To determine when $X(\S)$ is separated, the valuative criterion of separatedness can be applied, which can be translated in terms of \emph{valuations} on a variety. This translation has nice properties over the regular ones. 

Hironaka's Theorem, over fields of characteristic zero, guarantees the existence of a resolution of singularities 
\[ \pi : \tilde{Y} \to Y,\] 
with $\tilde{Y}$ regular. At the moment, it is possible to resolve singularities on positive characteristic for curves, surfaces (see: \cite{Abh56} and \cite{Lip78}) and threefolds (see: \cite{Abh66}, \cite{CP08} and \cite{CP09}). Thus, let us suppose that $Y$ admits a resolution of singularities $\pi : \tilde{Y} \to Y$ and let $\pi^{*}\S:=\{\pi^{*}\D \mid \D\in \S\}$ the pullback of $\S$, which is a divisorial fan over $\tilde{Y}$. The schemes $X(\S)$ and $X(\pi^{*} \S)$ are $T$-equivariantly isomorphic. Thus, if $Y$ admits a resolution of singularities, we can always construct a divisorial fan defining the same scheme but now over a regular semi-projective variety.

A valuation $\mu$ on a variety $Y$ over $k$ is a valuation on its field of rational functions $\mu:k(Y)^{*}\to \Q$ with $\mu=0$ along $k$. We say that $y\in Y$ is the (unique) center of $\mu$ if the valuation ring $(\mathcal{O}_{\mu},\mathfrak{m}_{\mu})$ dominates the local ring $(\mathscr{O}_{y},\mathfrak{m}_{y})$, which means that $\mathscr{O}_{y}\subset \mathcal{O}_{\mu}$ and $\mathfrak{m}_{y}=\mathfrak{m}_{\mu}\cap\mathscr{O}_{y}$ hold. For a valuation $\mu$ on $Y$ with center $y\in Y$, we have a well-defined group homomorphism
\begin{align*}
\mu:\cadiv(Y) &\to \Q ,\\
D &\mapsto \mu(f),
\end{align*}
where $D=\mathrm{div}(f)$ near $y$ with $f\in k(Y)$. Besides, for $Y$ regular, this provides a weight function $\mu:\{\text{prime divisors on }Y\}\to \Q$.

\begin{definition}
Let $N$ be a lattice, $\omega\subset N_{\Q}$ a cone and $\D\in\PPDiv_{\Q}(Y,\omega)$, with $Y$ a normal variety over $k$. If $\mu:\{\textrm{prime divisors of }Y\}\to\R$ is any map, then we define the associated \textit{weighted sum of polyhedral coefficients} as 
\[\mu(\D)=\sum \mu(D)\Delta_{D}\in\mathrm{Pol}_{\omega}(N).\]
\end{definition}

\begin{definition}\label{definition separated fans}
Let $Y$ be a regular semi-projective variety over $k$ and let $(\S,Y)$ be a divisorial fan with polyhedral coefficients defined over $N_{\Q}$, with $N$ a lattice. 
\begin{enumerate}[i)]
    \item We say that $(\S,Y)$ is \textit{separated} if for any two pp-divisors $\D$ and $\D'$ in $\S$ and any valuation $\mu$ on $Y$ we have that $\mu(\D\cap\D')=\mu(\D)\cap\mu(\D')$.
    \item We say that $(\S,Y)$ is \textit{complete} if $Y$ is complete and for any valuation $\mu$ on $Y$ we have that $\mu(\S):=\{\mu(\D)\mid \D\in\S\}$ form a complete subdivision of $N_{\Q}$.
\end{enumerate}
\end{definition}

Separated divisorial fans induce normal varieties endowed with an effective action of a split algebraic torus. The following result is proved by Altmann, Hausen and S\"u\ss\, which corresponds to \cite[Theorem 7.5]{AHS08}. The proof relies on the valuative criterion of separatness and the valuative criterion of properness, then all the arguments in the original proof can be performed over any field.

\begin{proposition}\label{proposition separated divisorial fans}
Let $Y$ be a regular semi-projective variety over $k$ and $(\S,Y)$ be a divisorial fan. Then, 
\begin{enumerate}
    \item $X(\S)$ is separated if and only if $(\S,Y)$ is separated.
    \item $X(\S)$ is complete if and only if $(\S,Y)$ is complete.
\end{enumerate}
\end{proposition}

\subsection{From normal $T$-varieties to divisorial fans}\label{section: from normal T-varieties to divisorial fans}
\noindent

A normal $k$-variety endowed with an effective action of a split algebraic $k$-torus $T$ has a $T$-stable affine open covering by Sumihiro's Theorem. Such a covering can be taken stable under intersections, since a variety is a separated scheme and a finite intersection of affine open subvarieties is still affine. By \cref{theoremmainMN25}, each $T$-stable affine open subvariety arises from a pp-divisor. However, all these pp-divisors could be defined over different normal semi-projective varieties, and a divisorial fan is a set of pp-divisors defined over the same semi-projective variety. Then, we need to construct a normal semi-projective that allows us to \textit{transfer} all the pp-divisors to compare them.

The following result is a generalization of \cite[Theorem 5.6]{AHS08}, whose proof adapts to this frame work.

\begin{theorem}\label{theoremahssplitcase2}
Let $T$ be a split torus over $k$. Up to isomorphism, every normal $T$-variety over $k$ arises from a divisorial fan $(\S, Y)$ over $k$.
\end{theorem}

\begin{remark}\label{remarkseveraldivfan}
    A normal $T$-variety does not have a unique $T$-stable affine open covering. Therefore, it can arise from several divisorial fans. For example, $\Gm$ acting on $\A^{1}\times\P^{1}$ by multiplication on the first coordinate  \begin{align*}
        \Gm\times\left( \A^{1}\times \P^{1}\right) &\to \left( \A^{1}\times \P^{1}\right) \\ (\lambda,(x,p)) &\mapsto (\lambda x,p).
    \end{align*} 
    By taking any two pairs of affine coverings of $\P^{1}$, for example 
    \[\{\P\smallsetminus\{0\},\P^{1}\smallsetminus\{\infty\}\}\textrm{ and }\{\P^{1}\smallsetminus\{p\},\P^{1}\smallsetminus\{q\}\}\] with $p,q\in\P^{1}\smallsetminus\{0,\infty\}$, we have two different $T$-stable affine open coverings of $\A^{1}\times\P^{1}$. Denote $\Delta:=[1,+\infty[$. The divisorial fans are 
    \[\S_{1}:=\{\emptyset\otimes \{0\}+\Delta \otimes \{\infty\},\Delta\otimes \{0\}+\emptyset\otimes \{\infty\},\emptyset\otimes \{0\}+\emptyset\otimes \{\infty\}\}\]
     and 
    \[\S_{2}:=\{\emptyset\otimes \{p\}+\Delta \otimes\{q\},\emptyset\otimes \{q\}+\Delta \otimes\{q\},\emptyset\otimes \{p\}+\emptyset\otimes \{q\}\},\]
which are different.
\end{remark}

A divisorial fan depends on the choice of a covering. Besides, the pp-divisors of the divisorial fan are in correspondence with the elements of the respective covering.

\begin{corollary}\label{remarkcovdivfancorrespondence}
   Let $T$ be a split torus over $k$. Let $X$ be a normal $T$-variety over $k$ and $\mathcal{U}$ be a finite $T$-stable affine open covering of $X$ that is stable under intersections, then there exists a divisorial fan $(\S, Y)$ over $k$ such that each pp-divisor in $\S$ corresponds to an element of $\mathcal{U}$.
\end{corollary}

\begin{remark}\label{remarkdivisorialfancomplexity1}
  If $X$ is a complexity-one normal $T$-variety, i.e. $\mathrm{tr. deg}(k(X)^{T})=1$, then all the pp-divisors $\D$ can be taken over the same regular projective curve. 
\end{remark}

For complexity-one normal $T$-varieties, we can ensure that every divisorial fan describing a given normal $T$-variety is defined over the same base.

\begin{corollary}\label{corollary cone isomorphic base}
Let $T$ be a split torus over $k$. Let $X$ be a complexity-one normal $T$-variety over $k$. If $(\S,Y)$ and $(\S',Y')$ are two divisorial fans over regular projective curves $Y$ and $Y'$, respectively, such that $X(\S)\cong X\cong X(\S')$ as $T$-varieties, then $Y\cong Y'$.
\end{corollary}

\begin{proof}
Identify the affine normal $T$-varieties $X(\D)$ and $X(\D')$, for $\D\in\S$ and $\D'\in \S'$, with open subvarieties of $X$. Their common intersection $X':=\bigcap_{\D\in \S\cup\S'}X(\D)$ is a $T$-stable affine open subvariety of $X$ and induces a birational map between $Y$ and $Y'$. Since $Y$ and $Y'$ are regular projective curves, this map extends to an isomorphism.
\end{proof}

\begin{example}\label{exampleFr}
 Let $r\in\N$ with $r\geq 1$. Denote by $\mathbb{F}_{r}$ the Hirzebruch surface and let us give the divisorial fan associated to the diagonal embedding $\G_{\mathrm{m},k}\to \G_{\mathrm{m},k}^{2}$. As a toric variety, $\mathbb{F}_{r}$ is encoded by the fan $\Sigma$:

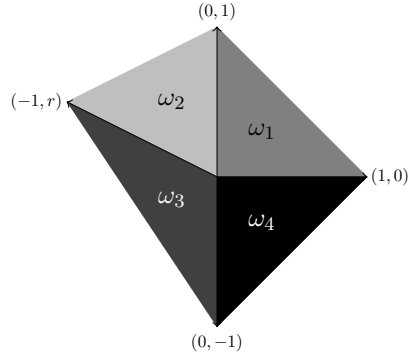
\begin{figure}
    \centering
    \begin{tikzpicture}
        \fill[gray] (0,0) -- (2,0) -- (0,2) -- (0,0);
        \fill[lightgray] (0,0) -- (0,2) -- (-2,1) -- (0,0);
        \fill[darkgray] (0,0) -- (-2,1) -- (0,-2) -- (0,0);
        \fill[black] (0,0) -- (2,0) -- (0,-2) -- (0,0);
        \node[scale=0.6] at (2.3,0) {$(1,0)$};
        \node[scale=0.6] at (0,2.2) {$(0,1)$};
        \node[scale=0.6] at (-2.4,1) {$(-1,r)$};
        \node[scale=0.6] at (0,-2.2) {$(0,-1)$};
        \draw[->] (0,0) -- (2,0);
        \draw[->] (0,0) -- (0,2);
        \draw[->] (0,0) -- (-2,1);
        \draw[->] (0,0) -- (0,-2);
        \node[scale=0.9] at (0.6,0.6) {$\omega_{1}$};
        \node[scale=0.9] at (-0.6,1) {$\omega_{2}$};
        \node[scale=0.9] at (-0.6,-0.3) {$\textcolor{white}{\omega_{3}}$};
        \node[scale=0.9] at (0.6,-0.6) {$\textcolor{white}{\omega_{4}}$};
    \end{tikzpicture}
    \caption{The fan of $\mathbb{F}_{r}$.}
    \label{figure: hirzebruch}
\end{figure}

The open subvarieties associated to the cones in $\Sigma$ are $\G_{\mathrm{m},k}$-stable. Thus, the set $\{X_{\omega}\mid \omega\in\Sigma\}$ is a $\G_{\mathrm{m},k}$-stable affine open covering of $\mathbb{F}_{r}$. The divisorial fan $(\S,\P_{k}^{1})$ associated to such a covering is generated by the pp-divisors:
\begin{itemize}
	\item $\D_{\omega_{1}}:=[1,\infty^{+}[\otimes\{\infty\},$
	\item $\D_{\omega_{2}}:=\left[-\frac{1}{r+1},0\right]\otimes\{\infty\}+\emptyset\otimes\{0\},$
	\item $\D_{\omega_{3}}:=\left]\infty^{-},-\frac{1}{r+1}\right[\otimes\{\infty\}$ and
	\item $\D_{\omega_{4}}:=\emptyset\otimes\{0\}+\left[0,1\right]\otimes\{\infty\},$
\end{itemize} 
 where $\omega_{1}':=\mathrm{Tail}(\D_{\omega_{1}})$ and $\omega_{3}':=\mathrm{Tail}(\D_{\omega_{3}})$. The divisorial fan can be put together into the following figure.

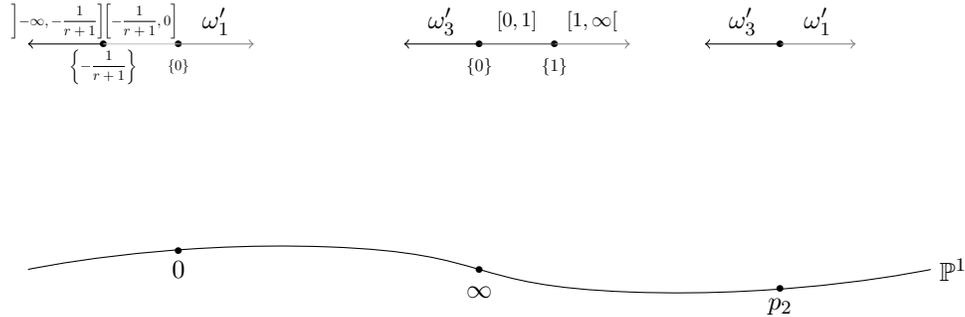
\begin{center}
	\begin{tikzpicture}
	\draw [black, xshift=1cm] plot [smooth, tension=1] coordinates { (0,0) (4,0.3) (8,-0.3) (12,0)};
	\node[scale=0.9] at (13.3,0) {$\P^{1}$};
	\node[scale=0.6] at (3,0.25) {$\bullet$};
	\node[scale=0.6] at (7,-0.01) {$\bullet$};
	\node[scale=0.6] at (11,-0.26) {$\bullet$};
	\node[scale=0.9] at (3,0) {$0$};
	\node[scale=0.9] at (7,-0.3) {$\infty$};
	\node[scale=0.9] at (11,-0.5) {$p_{2}$};
	\node[scale=0.6] at (3,3) {$\bullet$};
	\node[scale=0.6] at (2,3) {$\bullet$};
	\node[scale=0.6] at (7,3) {$\bullet$};
	\node[scale=0.6] at (8,3) {$\bullet$};
	\node[scale=0.6] at (11,3) {$\bullet$};
	
	
	
	
	
	
	\node[scale=0.5] at (2.5,3.3) {$\left[-\dfrac{1}{r+1},0\right]$};
	\node[scale=0.5] at (3,2.7) {$\left\{0\right\}$};
	\node[scale=0.5] at (2,2.7) {$\left\{-\dfrac{1}{r+1}\right\}$};
	\node[scale=0.9] at (6.5,3.3) {$\omega_{3}'$};
	\node[scale=0.9] at (10.5,3.3) {$\omega_{3}'$};
	
	\node[scale=0.5] at (1.4,3.3) {$\left]-\infty,-\dfrac{1}{r+1} \right]$};
	\node[scale=0.7] at (7.5,3.3) {$\left[0,1\right]$};
	\node[scale=0.6] at (7,2.7) {$\{0\}$};
	\node[scale=0.6] at (8,2.7) {$\{1\}$};
	
	\node[scale=0.9] at (3.5,3.3) {$\omega_{1}'$};
	\node[scale=0.7] at (8.5,3.3) {$\left[1,\infty\right[$};
	\node[scale=0.9] at (11.5,3.3) {$\omega_{1}'$};
	
	\draw[->,gray] (3,3) -- (4,3); 
	\draw[->,black] (2,3) -- (1,3); 
	\draw[-,lightgray] (3,3) -- (2,3); 
	
	\draw[->,black] (7,3) -- (6,3);
	\draw[->,gray] (8,3) -- (9,3);
	\draw[-,darkgray] (7,3) -- (8,3);
	
	\draw[->,gray] (11,3) -- (12,3);
	\draw[->,black] (11,3) -- (10,3);
	\end{tikzpicture}
\end{center}
The last diagram states that almost all the fibers over $\mathbb{P}_{k}^{1}$ are orbits whose closure is isomorphic to $\mathbb{P}_{k}^{1}$.
\end{example}

\begin{example}\label{exampleP3}
Let us consider $\P_{k}^{3}$ with the action of $\G_{\mathrm{m},k}$ given by 
\[(\lambda,\mu)\cdot [x_{0}:x_{1}:x_{2}:x_{3}]=[\lambda x_{0}:\mu x_{1}:\lambda\mu x_{2}:x_{3}].\]
The affine open covering given by $U_{i}:=\{x_{i}\neq 0\}$ is $\G_{\mathrm{m},k}$-stable. All the pp-divisors are defined over the quotient $\P^{3}_{k}\dasharrow \P_{k}^{1}$ given by 
\[[x_{0}:x_{1}:x_{2}:x_{3}]\mapsto [x_{0}x_{1}:x_{2}x_{3}].\]
Denote by $\D_{i}$ the pp-divisor corresponding to $U_{i}$. Notice that $U_{3}$, with the respective induced action of $\G_{\mathrm{m},k}$, corresponds to \cite[Example 4.34]{MN25}, but the nontrivial polyhedron is supported over $\{0\}$ instead of $\{\infty\}$ in $\P_{k}^{1}$. Thus, this $T$-variety is encoded by the following divisorial fan. 
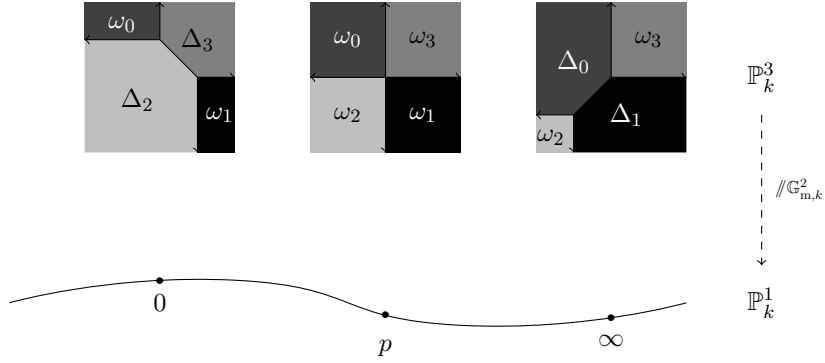
\begin{center}
	\begin{tikzpicture}
	\draw [black, xshift=1cm] plot [smooth, tension=1] coordinates { (0,0) (3,0.3) (6,-0.3) (9,0)};
	\node[scale=0.9] at (11,0) {$\P_{k}^{1}$};
	\node[scale=0.6] at (3,0.29) {$\bullet$};
	\node[scale=0.6] at (6,-0.17) {$\bullet$};
	\node[scale=0.6] at (9,-0.2) {$\bullet$};
	\node[scale=0.9] at (3,0) {$0$};
	\node[scale=0.9] at (6,-0.6) {$p$};
	\node[scale=0.9] at (9,-0.5) {$\infty$};
	
	\fill[gray] (3,3.5) -- (3,4) -- (4,4) -- (4,3) -- (3.5,3) ;
	\fill[gray] (6,3) -- (6,4) -- (7,4) -- (7,3);
	\fill[gray] (9,3) -- (9,4) -- (10,4) -- (10,3) -- (9,3) ;
	
	\fill[lightgray] (3,3.5) -- (3.5,3) -- (3.5,2) -- (2,2) -- (2,3.5) ;
	\fill[lightgray] (6,3) -- (6,2) -- (5,2) -- (5,3);
	\fill[lightgray] (8.5,2.5) -- (8,2.5) -- (8,2) -- (8.5,2) ;
	
	\fill[darkgray] (3,3.5) -- (3,4) -- (2,4) -- (2,3.5) -- (3,3.5);
	\fill[darkgray] (6,3) -- (5,3) -- (5,4) -- (6,4);
	\fill[darkgray] (9,3) -- (8.5,2.5) -- (8,2.5) -- (8,4) -- (9,4);
	
	\fill[black] (3.5,3) -- (4,3) -- (4,2) -- (3.5,2) ;
	\fill[black] (6,3) -- (7,3) -- (7,2) -- (6,2);
	\fill[black] (8.5,2.5) -- (9,3) -- (10,3) -- (10,2) -- (8.5,2);
	
	\node[scale=0.9] at (3.5,3.5) {$\Delta_{3}$};
	\node[scale=0.9] at (6.5,3.5) {$\omega_{3}$};
	\node[scale=0.9] at (9.5,3.5) {$\omega_{3}$};
	
	\node[scale=0.9] at (2.5,3.7) {$\textcolor{white}{\omega_{0}}$};
	\node[scale=0.9] at (5.5,3.5) {$\textcolor{white}{\omega_{0}}$};
	\node[scale=0.9] at (8.5,3.2) {$\textcolor{white}{\Delta_{0}}$};
	
	\node[scale=0.9] at (2.7,2.7) {$\Delta_{2}$};
	\node[scale=0.9] at (5.5,2.5) {$\omega_{2}$};
	\node[scale=0.9] at (8.2,2.2) {$\omega_{2}$};
	
	\node[scale=0.9] at (3.8,2.5) {$\textcolor{white}{\omega_{1}}$};
	\node[scale=0.9] at (6.5,2.5) {$\textcolor{white}{\omega_{1}}$};
	\node[scale=0.9] at (9.2,2.5) {$\textcolor{white}{\Delta_{1}}$};
	
	\draw[->] (3,3.5) -- (3,4);
	\draw[-] (3,3.5) -- (3.5,3);
	\draw[->] (3.5,3) -- (4,3);
	\draw[->] (3.5,3) -- (3.5,2);
	\draw[->] (3,3.5) -- (2,3.5);
	
	\draw[->] (6,3) -- (6,4);
	\draw[->] (6,3) -- (7,3);
	\draw[->] (6,3) -- (6,2);
	\draw[->] (6,3) -- (5,3);

	\draw[->] (9,3) -- (9,4); 
	\draw[->] (8.5,2.5) -- (8,2.5); 
	\draw[-] (8.5,2.5) -- (9,3);
	\draw[->] (9,3) -- (10,3); 
	\draw[->] (8.5,2.5) -- (8.5,2); 
	
	\node[scale=0.9] at (11,3) {$\P_{k}^{3}$};
	\draw[->,dashed] (11,2.5) -- (11,0.5);
	\node[scale=0.6] at (11.5,1.5) {$\sslash \G_{\mathrm{m},k}^{2}$};
	\end{tikzpicture} 
\end{center}

Such a diagram states that over almost all the points in $\P^{1}_{k}$ the fiber contains an orbit whose closure is isomorphic to $\P_{k}^{1}\times\P_{k}^{1}$.
\end{example}

\section{A localized category of pp-divisors and the general affine case}\label{Appendix localization}
\noindent

Then functor $X:\mathfrak{PPDiv}(L/k)\to\mathcal{E}(L/k)$ is faithful by \cref{Proposition X is faithful}, but it is not an equivalence of categories because it is not full. To address this obstruction, we \emph{localize} $\mathfrak{PPDiv}(L/k)$ (see \cite{KS06}), obtaining a category equivalent to $\mathcal{E}(L/k)$. In this section, we construct this category and prove that it satisfy the universal property of localization.


\subsection{A stuffed category}\label{section: a stuffed category}
\noindent

  The category $\mathfrak{PPDiv}(L/k)$ does not have enough morphisms. Our new category, denoted by $\mathfrak{PPDiv}_{S}(L/k)$, has the same objects as $\mathfrak{PPDiv}(L/k)$. Let us now define the morphisms. Let $\D$ and $\D'$ be two objects in $\mathfrak{PPDiv}_{S}(L/k)$, i.e. two pp-divisors over $L$. Consider the set $M_{\D,\D'}$ as the set of triangles 
\[\xymatrix{ & \kappa^{*}\D \ar[dl]_{(\kappa,\id_{N},\mathfrak{1})} \ar[rd]^{(\psi_{\gamma},F,\mathfrak{f})}   \\ \D & & \D' ,}\]
 where $\tilde{Y}$ is a normal semi-projective variety over $L$, $\kappa:\tilde{Y}\to Y$ is a projective and surjective morphism, and $(\psi_{\gamma},F,\mathfrak{f}):\kappa^{*}\D \to \D'$ is a dominating semilinear morphism of pp-divisors. In order to lighten the notation, the elements of $M_{\D,\D'}$ will be denoted by $(\kappa,\psi_{\gamma},F,\mathfrak{f})$. In $M_{\D,\D'}$ we define the following relation:
 \[ (\kappa,\psi_{\gamma},F,\mathfrak{f})\sim (\kappa',\psi'_{\gamma'},F',\mathfrak{f}') \]
 if and only if there exists a normal semi-projective variety $\hat{Y}$ and projective surjective morphisms $\hat{\kappa}:\hat{Y}\to \tilde{Y}$ and $\hat{\kappa}':\hat{Y}\to \tilde{Y}'$ such that
 \[(\kappa\circ \hat{\kappa},\psi_{\gamma}\circ\hat{\kappa},F,\hat{\kappa}^{*}\mathfrak{f})= (\kappa'\circ\hat{\kappa}',\psi'_{\gamma'}\circ\hat{\kappa}',F',\hat{\kappa}'^{*}\mathfrak{f}').\]
 This relation is clearly reflexive and symmetric.
 
 \begin{lemma}\label{lemma transitivity relation}
Let $(\kappa,\psi_{\gamma},F,\mathfrak{f})$ and $(\kappa',\psi'_{\gamma'},F',\mathfrak{f}')$ in $M_{\D,\D'}$. Then, \[(\kappa,\psi_{\gamma},F,\mathfrak{f})\sim (\kappa',\psi'_{\gamma'},F',\mathfrak{f}')\] if and only if $X(\psi_{\gamma},F,\mathfrak{f})=X(\psi'_{\gamma'},F',\mathfrak{f}')$.
 \end{lemma}
 
 \begin{proof}
 If $(\kappa,\psi_{\gamma},F,\mathfrak{f})\sim (\kappa',\psi'_{\gamma'},F',\mathfrak{f}')$ we have the following commutative diagram 
 \[ \xymatrix{
 & & & \kappa^{*} \D \ar[rd]^{(\psi_{\gamma},F,\mathfrak{f})} \ar[ld]_{(\kappa,\id_{N},\mathfrak{1})} & \\
 \hat{\kappa}^{*}\kappa^{*}\D = \hat{\kappa}'^{*}\kappa'^{*}\D \ar@/^2pc/[rrru]^{(\hat{\kappa},\id_{N},\mathfrak{1})} \ar@/_2pc/[rrrd]_{(\hat{\kappa}',\id_{N},\mathfrak{1})} & & \D   & & \D' \\
 & & & \kappa'^{*}\D \ar[ru]_{(\psi'_{\gamma'},F',\mathfrak{f}')} \ar[lu]^{(\kappa',\id_{N'},\mathfrak{1})} &.
  } \]
  This diagram induces the following in terms of affine normal varieties.
 \[ \xymatrix{
 & & & X(\kappa^{*} \D) \ar[rd]^{X(\psi_{\gamma},F,\mathfrak{f})} \ar[ld]_{X(\kappa,\id_{N},\mathfrak{1})} \ar[dd]^{\id} & \\
 X(\hat{\kappa}^{*}\kappa^{*}\D) = X(\hat{\kappa}'^{*}\kappa'^{*}\D) \ar@/^2pc/[rrru]^{X(\hat{\kappa},\id_{N},\mathfrak{1})} \ar@/_2pc/[rrrd]_{X(\hat{\kappa}',\id_{N},\mathfrak{1})} \ar[rr]^{\id} & & X(\D)   & & X(\D') \\
 & & & X(\kappa'^{*}\D) \ar[ru]_{X(\psi'_{\gamma'},F',\mathfrak{f}')} \ar[lu]^{X(\kappa',\id_{N'},\mathfrak{1})} &.
  } \]
 Given that, except for $X(\psi_{\gamma},F,\mathfrak{f})$ and $X(\psi'_{\gamma'},F',\mathfrak{f}')$, all the morphisms in the diagram are identities, we have that $X(\psi_{\gamma},F,\mathfrak{f})=X(\psi'_{\gamma'},F',\mathfrak{f}')$.
 
 In the other direction, if we assume $X(\psi_{\gamma},F,\mathfrak{f})=X(\psi'_{\gamma'},F',\mathfrak{f}')$, we have the following commutative diagram
 \[ \xymatrix{
 & & X(\kappa^{*} \D) \ar[rd]^{X(\psi_{\gamma},F,\mathfrak{f})} \ar[ld]_{X(\kappa,\id_{N},\mathfrak{1})} \ar[dd]^{\id} & \\
 & X(\D)   & & X(\D') \\
  & & X(\kappa'^{*}\D) \ar[ru]_{X(\psi'_{\gamma'},F',\mathfrak{f}')} \ar[lu]^{X(\kappa',\id_{N'},\mathfrak{1})} &.
  } \]
  By applying \cite[Theorem 5.11]{MN25} on the identity $\id:X(\kappa^{*}\D)\to X(\kappa'^{*}\D)$, we obtain the following commutative diagram
  \[ \xymatrix{
 & & & \kappa^{*} \D \ar[rd]^{(\psi_{\gamma},F,\mathfrak{f})} \ar[ld]_{(\kappa,\id_{N},\mathfrak{1})} & \\
 \hat{\kappa}^{*}\kappa^{*}\D = \hat{\kappa}'^{*}\kappa'^{*}\D \ar@/^2pc/[rrru]^{(\hat{\kappa},\id_{N},\mathfrak{1})} \ar@/_2pc/[rrrd]_{(\hat{\kappa}',\id_{N},\mathfrak{1})} & & \D   & & \D' \\
 & & & \kappa'^{*}\D \ar[ru]_{(\psi'_{\gamma'},F',\mathfrak{f}')} \ar[lu]^{(\kappa',\id_{N'},\mathfrak{1})} &.
  } \]
  Hence, we have that $(\kappa,\psi_{\gamma},F,\mathfrak{f})\sim (\kappa',\psi'_{\gamma'},F',\mathfrak{f}')$. Then, the assertion holds.
 \end{proof}
 
 By \cref{lemma transitivity relation}, the relation defined above is transitive. This implies that the relation is an equivalence relation. Thus, we define 
 \[\mathrm{Mor}_{\mathfrak{PPDiv}_{S}(L/k)}(\D,\D'):=M_{\D,\D'}/\sim.\]
 
 \begin{proposition}
The class $\mathfrak{PPDiv}_{S}(L/k)$ defines a category.
 \end{proposition}
 
 \begin{proof}

We need to prove that the composition is well defined. Let $\D$, $\D’$, and $\D''$ be objects in $\mathfrak{PPDiv}_{S}(L/k)$ and $(\kappa,\psi_{\gamma},F,\mathfrak{f}):\D\to\D'$ and $(\kappa',\psi'_{\gamma'},F',\mathfrak{f}'):\D'\to\D''$ be respectively in $\mathrm{Mor}_{\mathfrak{PPDiv}_{S}(L/k)}(\D,\D')$ and $\mathrm{Mor}_{\mathfrak{PPDiv}_{S}(L/k)}(\D',\D'')$. These morphisms form the following diagram
\[\xymatrix{ & \kappa^{*}\D \ar[dl]_{(\kappa,\id_{N},\mathfrak{1})} \ar[rd]^{(\psi_{\gamma},F,\mathfrak{f})} & & \kappa'^{*}\D' \ar[dl]^{(\kappa',\id_{N'},\mathfrak{1})} \ar[rd]^{(\psi'_{\gamma'},F',\mathfrak{f}')} & \\ \D & & \D' & & \D'' .}\]
The morphism $(\kappa,\psi_{\gamma},F,\mathfrak{f}):\D\to\D'$ induces a dominant equivariant semilinear morphism $X(\kappa,\psi_{\gamma},F,\mathfrak{f}):X(\D)\to(\D')$, which corresponds to 
\[X(\psi_{\gamma},F,\mathfrak{f}):X(\kappa^{*}\D)\to X(\D').\]
 Hence, given that $X(\D)=X(\kappa^{*}\D)$ and $X(\D')=X(\kappa'^{*}\D')$, we have the following morphism $X(\kappa,\psi_{\gamma},F,\mathfrak{f}):X(\kappa^{*}\D)\to X(\kappa'^{*}\D')$. Then, by \cite[Theorem 5.11]{MN25}, there exists a morphism $(\tilde{\kappa},\tilde{\psi}_{\gamma},\tilde{F},\tilde{\mathfrak{f}})$ in $\mathrm{Mor}_{\mathfrak{PPDiv}_{S}(L/k)}(\kappa^{*}\D,\kappa'^{*}\D')$ such that $X(\tilde{\kappa},\tilde{\psi}_{\gamma},\tilde{F},\tilde{\mathfrak{f}})=X(\kappa,\psi_{\gamma},F,\mathfrak{f})$. Notice that $\tilde{F}=F$ because both represent the same morphisms $\varphi_{\gamma}:T\to T'$. All these three localized dominating semilinear morphisms of pp-divisors fit into the following commutative diagram
\[\xymatrix{ 
 && \tilde{\kappa}^{*}\kappa^{*}\D \ar[dl]_{(\tilde{\kappa},\id_{N},\mathfrak{1})} \ar[rd]^{(\tilde{\psi}_{\gamma},F,\tilde{\mathfrak{f}})} && \\
& \kappa^{*}\D \ar[dl]_{(\kappa,\id_{N},\mathfrak{1})} \ar[rd]^{(\psi_{\gamma},F,\mathfrak{f})} & & \kappa'^{*}\D' \ar[dl]^{(\kappa',\id_{N},\mathfrak{1})} \ar[rd]^{(\psi'_{\gamma'},F',\mathfrak{f}')} & \\ 
\D & & \D' & & \D'' .}\]
Then, the composition gives 
\[(\kappa,\psi_{\gamma},F,\mathfrak{f})\circ(\kappa',\psi'_{\gamma'},F',\mathfrak{f}')=(\tilde{\kappa}\circ \kappa,\psi'_{\gamma'}\circ \tilde{\psi}_{\gamma},F'\circ F,\tilde{\psi}_{\gamma}^{*}(\mathfrak{f})\cdot F'_{*}(\tilde{\mathfrak{f}}))\]
In $\mathrm{Mor}_{\mathfrak{PPDiv}_{S}(L/k)}(\D,\D'')$. Then, the morphisms can be composed and, therefore, $\mathfrak{PPDiv}_{S}(L/k)$ defines a category.
\end{proof}

\begin{remark}
Notice that if $\kappa$ and $\kappa'$ are identities, we have that the composition turns to be 
\[(\id_{Y},\psi_{\gamma},F,\mathfrak{f})\circ(\id_{Y'},\psi'_{\gamma'},F',\mathfrak{f}')=(\id_{Y},\psi'_{\gamma'}\circ \psi_{\gamma},F'\circ F,\psi_{\gamma}^{*}(\mathfrak{f})\cdot F'_{*}(\mathfrak{f})), \]
since the semilinear morphisms of pp-divisors $(\kappa,\id_{N},\mathfrak{1})$ and $(\kappa',\id_{N'},\mathfrak{1})$ are isomorphisms.
\end{remark}

The category $\mathfrak{PPDiv}_{S}(L/k)$ has the same objects as $\mathfrak{PPDiv}(L/k)$. We prove now that the latter can be seen as a subcategory of the former.

The assignment $Q:\mathfrak{PPDiv}(L/k)\to \mathfrak{PPDiv}_{S}(L/k)$, defined as the identity on objects and by $Q(\psi_{\gamma},F,\mathfrak{f})=(\id_{Y},\psi_{\gamma},F,\mathfrak{f})$ in morphisms, satisfies the following
\begin{align*} 
Q(\psi_{\gamma},F,\mathfrak{f})\circ Q(\psi'_{\gamma'},F',\mathfrak{f}') &= (\id_{Y},\psi_{\gamma},F,\mathfrak{f})\circ (\id_{Y'},\psi'_{\gamma'},F',\mathfrak{f}') \\
&= (\id_{Y},\psi'_{\gamma'}\circ \psi_{\gamma},F'\circ F,\psi_{\gamma}^{*}(\mathfrak{f})\cdot F'_{*}(\mathfrak{f})) \\
&= Q(\psi'_{\gamma'}\circ \psi_{\gamma},F'\circ F,\psi_{\gamma}^{*}(\mathfrak{f})\cdot F'_{*}(\mathfrak{f})).
\end{align*} 
Then, $Q:\mathfrak{PPDiv}(L/k)\to \mathfrak{PPDiv}_{S}(L/k)$ defines a functor.

\begin{lemma}
The functor $Q:\mathfrak{PPDiv}(L/k)\to \mathfrak{PPDiv}_{S}(L/k)$ is faithful.
\end{lemma}

\begin{proof}
This is a consequence of \cref{lemma transitivity relation}. Let $\D$ and $\D'$ be objects in $\mathfrak{PPDiv}(L/k)$ and $(\psi_{\gamma},F,\mathfrak{f})$ and $(\psi'_{\gamma'},F',\mathfrak{f}')$ be two morphisms of pp-divisors such that 
\[Q(\psi_{\gamma},F,\mathfrak{f})= Q(\psi'_{\gamma'},F',\mathfrak{f}').\]
By \cref{lemma transitivity relation}, we have that this is equivalent to 
\[X(\psi_{\gamma},F,\mathfrak{f})= X(\psi'_{\gamma'},F',\mathfrak{f}'),\]
which implies that $(\psi_{\gamma},F,\mathfrak{f})=(\psi'_{\gamma'},F',\mathfrak{f}')$ by \cref{Proposition X is faithful}.
\end{proof}

According to \cite[Definition 7.1.1]{KS06}, given a collection of morphisms $S$ in a category $\mathcal{C}$, a localization $\mathcal{C}_{S}$ is a category satisfying a certain universal property that ``makes elements in $S$ invertible''. Under suitable conditions on $S$, one can ensure the existence of the category $\mathcal{C}_{S}$. However, the family of morphisms that we need to make invertible does not form a \textit{right multiplicative system}, in the sense of \cite[Definition 7.1.5]{KS06}. The axioms $S3$ and $S4$ in \cite[Definition 7.1.5]{KS06} both fail. Nevertheless, we prove here below that, if $S$ is the collection of morphisms in $\mathfrak{PPDiv}(L/k)$ of the form $(\kappa,\id,\mathfrak{1})$, then $\mathfrak{PPDiv}_{S}(L/k)$ is a localization with respect to $S$. We denote $s_{\kappa}:=(\kappa,\id,\mathfrak{1})$.

\begin{proposition}
The category $\mathfrak{PPDiv}_{S}(L/k)$ is a localization of the category $\mathfrak{PPDiv}(L/k)$ with respect to $S$.
\end{proposition}

\begin{proof}
In order to prove this, we verify the universal property described in \cite[Definition 7.1.1]{KS06}. Let $\mathcal{C}$ be a category and $G: \mathfrak{PPDiv}(L/k)\to \mathcal{C}$ be a functor such that $G(s)$ is invertible in $\mathrm{Mor}_{\mathcal{C}}$ for every $s\in S$. Let us define the functor $G_{S}:\mathfrak{PPDiv}_{S}(L/k)\to \mathcal{C}$ given by $G_{S}(\D):=\D$ and 
\[G_{S}(\kappa,\psi_{\gamma},F,\mathfrak{f}):=G(\psi_{\gamma},F,\mathfrak{f})\cdot G(\kappa,\id,\mathfrak{1})^{-1}.\]
This assignment does not depend on the representatives. Indeed, if $(\kappa,\psi_{\gamma},F,\mathfrak{f})\sim (\kappa',\psi'_{\gamma},F',\mathfrak{f}')$ we have the following commutative diagram 
 \[ \xymatrix{
 & & & \kappa^{*} \D \ar[rd]^{(\psi_{\gamma},F,\mathfrak{f})} \ar[ld]_{(\kappa,\id_{N},\mathfrak{1})} & \\
 \hat{\kappa}^{*}\kappa^{*}\D = \hat{\kappa}'^{*}\kappa'^{*}\D \ar@/^2pc/[rrru]^{(\hat{\kappa},\id_{N},\mathfrak{1})} \ar@/_2pc/[rrrd]_{(\hat{\kappa}',\id_{N},\mathfrak{1})} & & \D   & & \D' \\
 & & & \kappa'^{*}\D \ar[ru]_{(\psi'_{\gamma},F',\mathfrak{f}')} \ar[lu]^{(\kappa',\id_{N'},\mathfrak{1})} &.
  } \]
  Under the functor $G$, this diagram yields the following
 \[G(\psi_{\gamma},F,\mathfrak{f})\cdot G(s_{\hat{\kappa}}) \cdot G(s_{\hat{\kappa}})^{-1}\cdot G(s_{\kappa})^{-1}=G(\psi_{\gamma}',F',\mathfrak{f}')\cdot G(s_{\hat{\kappa}'}) \cdot G(s_{\hat{\kappa}'})^{-1}\cdot G(s_{\kappa}')^{-1},\]
and then
\begin{align*}
G_{S}(\kappa,\psi_{\gamma},F,\mathfrak{f}) &= G(\psi_{\gamma},F,\mathfrak{f}) \cdot G(s_{\kappa})^{-1} \\ &=G(\psi_{\gamma}',F',\mathfrak{f}')\cdot G(s_{\kappa'})^{-1} \\ &=G_{S}(\kappa',\psi_{\gamma}',F',\mathfrak{f}').
\end{align*}

Let us denote by $Q:\mathfrak{PPDiv}(L/k)\to\mathfrak{PPDiv}_{S}(L/k)$, given by $Q(\D):=\D$ and $Q(\psi_{\gamma},F,\mathfrak{f}):=(\id,\psi_{\gamma},F,\mathfrak{f})$, the localization functor. Let $\D$ and $\D'$ be objects of $\mathfrak{PPDiv}(L/k)$ and $(\psi_{\gamma},F,\mathfrak{f}):\D\to \D'$ be a morphism in $\mathfrak{PPDiv}(L/k)$. Notice that 
\[G_{S}\circ Q(\psi_{\gamma},F,\mathfrak{f})=G_{S}(Q(\psi_{\gamma},F,\mathfrak{f}))=G_{S}(\id,\psi_{\gamma},F,\mathfrak{f})=G(\psi_{\gamma},F,\mathfrak{f})\cdot G(s_{\id})^{-1}.\] 
Given that $G(s_{\id})=\id$, the following diagram commutes
\[\xymatrixcolsep{6pc}\xymatrix{ G(\D) \ar[r]^{G(\psi_{\gamma},F,\mathfrak{f})} & G(\D') \\ G_{S}\circ Q(\D) \ar[r]_{G_{S}\circ Q(\psi_{\gamma},F,\mathfrak{f})} \ar[u]^{\id} & G_{S}\circ Q(\D') \ar[u]_{\id} .}\]
Thus, this data gives rise to a natural transformation $\eta:G\Rightarrow G_{S}\circ Q$ that is an isomorphism in the category of functors.

Let $G_{1}$ and $G_{2}$ be elements of $\textbf{Fct}(\mathfrak{PPDiv}_{S}(L/k), \mathcal{C})$, the category of functors between $\mathfrak{PPDiv}_{S}(L/k)$ and $\mathcal{C}$. There is a natural map 
\[\varrho:\mathrm{Hom}_{\textbf{Fct}(\mathfrak{PPDiv}_{S}(L/k), \mathcal{C})}(G_{1},G_{2})\to \mathrm{Hom}_{\textbf{Fct}(\mathfrak{PPDiv}(L/k), \mathcal{C})}(G_{1}\circ Q,G_{2}\circ Q).\]
This map is injective since $Q$ is the identity on objects. Let us prove that it is also surjective. Let $\eta: G_{1}\circ Q\Rightarrow G_{2}\circ Q$ be a natural transformation. Given that $\mathfrak{PPDiv}(L/k)$ and $\mathfrak{PPDiv}_{S}(L/k)$ have the same objects, $\eta$ defines maps $\eta_{S,\D}:=\eta_{\D}:G_{1}(\D)\to G_{2}(\D)$, for each pp-divisor $\D$ in $\mathfrak{PPDiv}_{S}(L/k)$. In order to prove that $\eta_{S}$ defines a natural transformation between $\mathfrak{PPDiv}_{S}(L/k)$ and $\mathcal{C}$, it suffices to prove that the following diagram commutes.
\[\xymatrixcolsep{6pc}\xymatrix{ G_{1} (\D) \ar[r]^{G_{1}(\kappa,\psi_{\gamma},F,\mathfrak{f})} \ar[d]_{\eta_{S,\D}} & G_{1} (\D') \ar[d]^{\eta_{S,\D'}} \\ G_{2}(\D) \ar[r]_{G_{2}(\kappa,\psi_{\gamma},F,\mathfrak{f})} & G_{2}(\D') ,}\]
for every pair of pp-divisors $\D$ and $\D'$ in $\mathfrak{PPDiv}_{S}(L/k)$. In $\mathfrak{PPDiv}_{S}(L/k)$, we have that 
\[\xymatrix{ & \kappa^{*}\D \ar[ld]_{Q(\kappa,\id,\mathfrak{1})} \ar[rd]^{Q(\psi_{\gamma},F,\mathfrak{f})}& \\ \D \ar[rr]_{(\kappa,\psi_{\gamma},F,\mathfrak{f})} & & \D' .}\]
Otherwise stated, $(\kappa,\psi_{\gamma},F,\mathfrak{f})=Q(\psi_{\gamma},F,\mathfrak{f})\cdot Q(\kappa,\id,\mathfrak{1})^{-1}$. Then,
\begin{align*}
G_{i} (\kappa,\psi_{\gamma},F,\mathfrak{f}) &= G_{i}\left( Q(\psi_{\gamma},F,\mathfrak{f})\cdot Q(\kappa,\id,\mathfrak{1})^{-1}\right) \\
&= G_{i}\left(Q(\psi_{\gamma},F,\mathfrak{f})\right)\cdot G_{i}\left(Q(\kappa,\id,\mathfrak{1})^{-1}\right) \\
&= \left((G_{i}\circ Q)(\psi_{\gamma},F,\mathfrak{f})\right)\cdot \left((G_{i}\circ Q)(\kappa,\id,\mathfrak{1})\right)^{-1},
\end{align*}
for every $i\in\{1,2\}$. Besides, we have the following diagram, where the two rectangles in the back are commutative:
\[\xymatrixcolsep{5pc}\xymatrixrowsep{1pc}\xymatrix{ & (G_{1}\circ Q)(\kappa^{*}\D) \ar@{->}'[d][dd]_(0.3){\eta_{\kappa^{*}\D}} \ar[ld]_{(G_{1}\circ Q)(\kappa,\id,\mathfrak{1})} \ar[rd]^{(G_{1}\circ Q)(\psi_{\gamma},F,\mathfrak{f})}& \\ (G_{1}\circ Q)(\D) \ar[dd]_{\eta_{\D}} \ar[rr]_(0.7){G_{1}(\kappa,\psi_{\gamma},F,\mathfrak{f})} & & (G_{1}\circ Q)(\D') \ar[dd]^{\eta_{\D'}} \\ & (G_{2}\circ Q)(\kappa^{*}\D) \ar[ld]_{(G_{2}\circ Q)(\kappa,\id,\mathfrak{1})} \ar[rd]^{(G_{2}\circ Q)(\psi_{\gamma},F,\mathfrak{f})}& \\ (G_{2}\circ Q)(\D) \ar[rr]_{G_{2}(\kappa,\psi_{\gamma},F,\mathfrak{f})} & & (G_{2}\circ Q)(\D').  }\]
And since $(G_i\circ Q)(\kappa,\id,\mathfrak{1})$ is invertible in $\mathcal{C}$, we have that the rectangle in the front is commutative. Thus
\[\eta_{S,\D'}\cdot G_{1} (\kappa,\psi_{\gamma},F,\mathfrak{f}) = G_{2} (\kappa,\psi_{\gamma},F,\mathfrak{f}) \cdot \eta_{S,\D}.\]
Thus, $\eta_{S}$ is a natural transformation and $\varrho(\eta_{S})=\eta$. Hence, $\mathfrak{PPDiv}_{S}(L/k)$ is a localization of $\mathfrak{PPDiv}(L/k)$ with respect to $S$.
\end{proof}

\subsection{Equivalence of categories} 
\noindent

By the universal property of the localization, the functor $X:\mathfrak{PPDiv}(L/k)\to \mathcal{E}(L/k)$ fits into the following commutative diagram
\[ \xymatrix{ \mathfrak{PPDiv}_{S}(L/k) \ar[r]^<<{\mathcal{X}} & \mathcal{E}(L/k) \\
\mathfrak{PPDiv}(L/k) \ar[u]^{Q} \ar[ur]_{X} . } \]

\begin{proposition}\label{proposition equivalence of categories modified}
The functor $\mathcal{X}:\mathfrak{PPDiv}_{S}(L/k)\to\mathcal{E}(L/k)$ is an equivalence of categories.
\end{proposition}

\begin{proof}
The functor $\mathcal{X}$ is essentially surjective, since $X$ is essentially surjective and the diagram above is commutative. The fullness is a consequence of \cite[Theorem 5.11]{MN25}. The faithfulness follows from \cref{lemma transitivity relation} and the fact that 
\[\mathcal{X}(\kappa,\psi_{\gamma},F,\mathfrak{f})=X(\psi_{\gamma},F,\mathfrak{f})\circ X(\kappa,\id,\mathfrak{1})^{-1}=X(\psi_{\gamma},F,\mathfrak{f}).\] 
\end{proof}

\begin{remark}
Notice that \cref{proposition equivalence of categories modified} generalizes \cite[Corollary 8.14]{AH06}, whose proof is omitted in \textit{loc.~cit}. Indeed, it suffices to consider the semilinear morphisms with $\gamma$ the neutral element of the Galois group.
\end{remark}

\subsection{Affine case and general pp-divisors}\label{Section general pp-divisor}
\noindent

Let $T$ be an algebraic torus over $k$ that splits over $L$, and let $X$ be an affine normal $T$-variety over $k$. By \cite[Proposition 6.10]{MN25} and \cref{maintheoremofpaperaffine}, $X$ is encoded by a pair $(\D_{L},\varrho)$, where $\D_{L}$ is a minimal pp-divisor over $L$ and $\varrho$ is a Galois semilinear action. For a non-minimal pp-divisor $\D'$ with $X(\D)\cong X(\D')$, such an action may not exist, and we must consider a broader class of morphisms.

Indeed, by \cite[Theorem 5.11]{MN25}, every dominant semilinear equivariant morphism $(\varphi_{\gamma},f_{\gamma}):X(\D)\to X(\D')$ fits into a commutative diagram
\[
\xymatrix{
& X(\kappa^{*}\D) \ar[dl]_{X(\kappa,\id_{N},1)}^{\cong}
\ar[rd]^{X(\psi_{\gamma},F,\mathfrak{f})} & \\
X(\D) \ar[rr]_{(\varphi_{\gamma},f_{\gamma})} && X(\D'),
}
\]
where $\kappa:\tilde{Y}\to Y$ is projective and birational and $(\psi_{\gamma},F,\mathfrak{f}):\kappa^{*}\D\to\D'$ is a semilinear morphism of pp-divisors. Thus, as explained in \cref{section: a stuffed category}, the diagram
\[
\xymatrix{
& \kappa^{*}\D \ar[dl]_{(\kappa,\id_{N},1)}
\ar[rd]^{(\psi_{\gamma},F,\mathfrak{f})} & \\
\D && \D'
}
\]
defines a morphism in the localization $\mathfrak{PPDiv}_{S}(L/k)$.

For every pp-divisor $\D$ in $\mathfrak{PPDiv}_{S}(L/k)$, we denote by $\SAut_{\mathrm{loc}}(\D)$ its group of \textit{localized} semilinear automorphisms, i.e. its automorphisms in the category $\mathfrak{PPDiv}_{S}(L/k)$. The subgroup of localized automorphisms is denoted by $\Aut_{\mathrm{loc}}(\D)$.

\begin{definition}
Let $\D$ be an object in the localization $\mathfrak{PPDiv}_{S}(L/k)$. Let $H$ be an abstract group. A \textit{localized semilinear action of $H$ on $\D$}, or a \textit{$H$-localized semilinear action}, is a group homomorphism $\varrho:H\to \SAut_{\mathrm{loc}}(\D)$. A \textit{Galois localized semilinear action} is a localized semilinear action of $\Gamma_{L}$ on $\D$ and $\varrho$ is a section of the sequence
\[1\to\Aut_{\mathrm{loc}}(\D)\to\SAut_{\mathrm{loc}}(\D)\to\Gamma.\]
\end{definition}

Let $H$ be an abstract group. A $H$-localized semilinear action 
\[\varrho:H\to \SAut_{\mathrm{loc}}(\D)\] induces a $H$-semilinear equivariant action (see \cite[Section 6.1]{MN25})
\[\mathcal{X}(\varrho):H\to \SAut(T;X(\D))\]
via the functor $\mathcal{X}:\mathfrak{PPDiv}_{S}(L/k)\to \mathcal{E}(L/k)$. Given that $\mathcal{X}$ is an equivalence of categories by \cref{proposition equivalence of categories modified}, every $H$-semilinear equivariant action $\varrho:H\to\SAut(T;X(\D))$ arises from a $H$-localized semilinear action of pp-divisors. Actually, this defines a bijection between the set of localized semilinear actions on $\D$ and the set of semilinear equivariant actions on $X(\D)$. In other words, we have the following immediate result.

\begin{proposition}\label{proposition bijection localized semilinear actions}
Let $\D$ be an object in $\mathfrak{PPDiv}_{S}(L/k)$. Then, there exists a bijection between the set of localized semilinear actions on $\D$ and the set of semilinear equivariant actions on $X(\D)$.
\end{proposition}

Let $\mathfrak{PPDiv}_{S}(\Gamma_{L})$ be the category of pairs $(\D,\varrho)$, where $\D$ is a pp-divisor over $L$ and $\varrho:\Gamma_{L}\to \SAut_{\mathrm{loc}}(\D)$ is a Galois localized semilinear action. A morphism in this category is a morphism $(\psi,F,\mathfrak{f}):\D\to\D'$ in $\mathfrak{PPDiv}_{S}(L/k)$ such that 
\[\varrho'_{\gamma}\circ (\psi,F,\mathfrak{f})=(\psi,F,\mathfrak{f})\circ \varrho_{\gamma},\]
for every $\gamma \in \Gamma_{L}$. Let $(\D, \varrho)$ be an object in $\mathfrak{PPDiv}_{S}(\Gamma_{L})$. By \cref{theoremmainMN25}, $\mathcal{X}(\D)$ is a normal $T$-variety over $L$, where $T$ denotes its respective torus. Moreover, by \cref{proposition bijection localized semilinear actions}, $\mathcal{X}(\D)$ comes with a Galois semilinear equivariant action $\mathcal{X}(\varrho):G\to\SAut(T;X(\D))$. Then, by \cite[Proposition 6.10]{MN25}, there exists an affine normal $T$-variety $X:=\mathcal{X}(\D,\varrho)$ over $k$ such that $X_{L}\cong \mathcal{X}(\D)$ as $T$-varieties over $L$. This proves the first part of the following theorem.

\begin{theorem}\label{theorem main affine general}
Let $k$ be a field, $L/k$ be a finite Galois extension with Galois group $\Gamma_{L}$.
	\begin{enumerate}[a)]
            \item\label{theorem main affine general part i} Let $(\D_{L},\varrho)$ be an object in $\mathfrak{PPDiv}_{S}(\Gamma_{L})$. Then, $\mathcal{X}(\D_{L},\varrho)$ is an affine normal variety endowed with an effective action of an algebraic torus $T$ over $k$ such that $T$ splits over $L$ and $\mathcal{X}(\D_{L},\varrho)_{L}\cong X(\D_{L})$ as $T_{\D_{L}}$-varieties over $L$.
            \item\label{theorem main affine general part ii} Let $X$ be an affine normal variety over $k$ endowed with an effective $T$-action such that $T_{L}$ is split. Let $\D_L$ be a pp-divisor such that $X_{L}\cong X(\D_L)$. Then, there exists a Galois localized semilinear action $\varrho:\Gamma_{L}\to \SAut(\D_{L})$ such that $X \cong \mathcal{X}(\D_{L},\varrho)$ as $T$-varieties.
        \end{enumerate}
\end{theorem}

\begin{proof}
Let us prove part \eqref{theorem main affine general part ii}, the remaining part of the theorem. Let $X$ be a normal variety over $k$ endowed with an effective $T$-action and $\D_L$ be a pp-divisor such that $X_{L}\cong X(\D_{L})=\mathcal{X}(\D_{L})$ as $T_{L}$-varieties over $L$. By \cite[Proposition 6.10]{MN25}, as a $T$-variety over $k$, $X$ is equivalent to a pair $(X_{L},\varrho')$, where $X_{L}$ is a normal $T_{L}$-variety, with $T_{L}$ split over $L$, and a Galois semilinear equivariant action $\varrho'$. Now, by \cref{proposition bijection localized semilinear actions}, we have that the Galois semilinear equivariant action on $X(\D_{L})$ induces a unique Galois localized semilinear action $\varrho$ on $\D_{L}$. Then, the pair $(\D_{L},\varrho)$ encodes the pair $(X_{L},\varrho')$.
\end{proof}

By \cref{theorem main affine general}, every pair $(\D,\varrho)$ in $\mathfrak{PPDiv}_{S}(\Gamma_{L})$ corresponds to an affine normal variety $\mathcal{X}(\D,\varrho)$ endowed with a torus action over $k$ that is split over $L$. This construction induces a functor 
 \begin{align*}
 \mathcal{X}:\mathfrak{PPDiv}_{S}(\Gamma_{L}) &\to \mathcal{E}(k,L);\\
  (\D,\varrho) &\mapsto X(\D,\varrho),
 \end{align*}
where $\mathcal{E}(k,L)$ is the category of affine normal varieties over $k$ endowed with an effective action of an algebraic torus over $k$ that is split over $L$. This functor is the composition of the functor $(\D,\varrho)\mapsto (\mathcal{X}(\D),\mathcal{X}(\varrho))$, from the category $\mathfrak{PPDiv}_{S}(\Gamma_{L})$ to the category of normal varieties endowed with an effective action of a split algebraic torus over $L$ and a Galois semilinear equivariant action, and the equivalence of categories of \cite[Proposition 6.10]{MN25}. Given that the first functor is faithful, covariant, and essentially surjective, we have the following.

\begin{proposition}\label{theorem equivalence categories localized descent}
The functor $\mathcal{X}:\mathfrak{PPDiv}_{S}(\Gamma_{L}) \to \mathcal{E}(k,L)$ is an equivalence of categories. In particular, the functor $X:\mathfrak{PPDiv}(\Gamma_{L}) \to \mathcal{E}(k,L)$ is faithful, covariant, and essentially surjective.
\end{proposition}

\section{Normal varieties and non-split tori}\label{section: normal varieties and non-split tori}
\noindent

The main goal of this section is to prove \cref{maintheoremofpaper}, which we recall below for convenience. We first introduce semilinear morphisms and automorphisms of divisorial fans, and then study semilinear actions on them. After proving the main theorem, we conclude with applications to toric and equivariantly normal varieties \cite{Bri24a,Bri24b,Bri26}.

       \begin{theorem}\label{theorem: section normal non-split maintheoremofpaper}
        Let $k$ be a field and $L/k$ be a finite Galois extension with Galois group $\Gamma_{L}$.
	\begin{enumerate}[a)]
            \item Let $T$ be a split algebraic torus over $L$ and $X$ be a normal $T$-variety over $L$. If there exists a divisorial fan $(\S_{L},Y_{L})$ for $X$ admitting a Galois semilinear action such that 
\[\text{the subvariety }X(\S(\D,\Gamma_{L}))\text{ is quasi-projective for every }\D\in\S_{L},\]
then there exists an algebraic torus $T'$ over $k$ and a normal $T'$-variety $X'$ over $k$ such that $X'_{L}\cong X$ as $T$ varieties over $L$.
            \item Let $T$ be an algebraic torus over $k$ that splits over $L$. Let $X$ be a normal variety endowed with an effective $T$-action over $k$. Then, there exists a divisorial fan $(\S_{L},Y_{L})$ admitting a Galois semilinear action such that
\[\text{the subvariety }X(\S(\D,\Gamma_{L}))\text{ is quasi-projective for every }\D\in\S_{L},\]
            and $X_{L}\cong X(\S_{L})$ as $T_{L}$-varieties. 
        \end{enumerate} 
    \end{theorem}

\subsection{Semilinear morphisms of divisorial fans}\label{sectioneqaut}
\noindent

For minimal pp-divisors, descent data can be described by semilinear morphisms, whereas the general case requires localized semilinear morphisms. Since a divisorial fan need not contain any minimal pp-divisor, we work with localized morphisms. We now introduce semilinear morphisms of divisorial fans.

\begin{definition}
Let $\D$ be an object in $\mathfrak{PPDiv}(k)$ and $\mathfrak{f}=\sum v_{i}\otimes f_{i}$ be a plurifunction in $k(N,Y)^{*}$. We define \textit{the restriction of $\mathfrak{f}$ to $\D$} as
     \[\mathfrak{f}|_{\D}:=\sum_{\mathrm{div}(f_{i})\subset\mathrm{Supp}(\D)}v_{i}\otimes f_{i}.\] 
\end{definition}
 
\begin{definition}\label{Def semilinear morphism of div fan}
Let $\gamma\in\Gamma_{L}$. Let $(\S,Y)$ and $(\S',Y')$ be two divisorial fans over $L$. A \textit{semilinear morphism of divisorial fans with respect to $\gamma$} is a triple $(\psi_{\gamma},F,\mathfrak{f})$ such that

\begin{enumerate}[i)]
	\item\label{Def semilinear morphism of div fan part 1} for every $\D\in\S$ there exists $\D'\in\S'$ such that the triple $(\psi_{\gamma},F,\mathfrak{f}|_{\D}):\D\to\D'$ is a semilinear morphism of pp-divisors and
	\item\label{Def semilinear morphism of div fan part 2} For every $\D,\D'\in\S$ as in \eqref{Def semilinear morphism of div fan part 1}, if $\mathfrak{E}\in\S$ is a face of $\D$, then there exists $\mathfrak{E}'\in\S$ such that $(\psi_{\gamma},F,\mathfrak{f}|_{\mathfrak{E}}):\mathfrak{E}\to\mathfrak{E}'$ is a semilinear morphism of pp-divisors and $\mathfrak{E}'$ is a face of $\D'$.
\end{enumerate}
\end{definition}

We cannot ask for uniqueness of $\D'$ in \cref{Def semilinear morphism of div fan} part \eqref{Def semilinear morphism of div fan part 1}, because if $\D'$ is a proper face of some other pp-divisor in $\S'$, then we have a map from $\D$ to that pp-divisor under the composition. However, this is the only obstruction to the uniqueness. In the following, we write $g_{\gamma}:=(\psi_{\gamma},F,\mathfrak{f})$ and $g_{\D,\gamma}:=(\psi_{\gamma},F,\mathfrak{f}|_{\D})$ for simplicity.

From a semilinear morphism of divisorial fans $g_{\gamma}:(\S,Y)\to (\S',Y')$, we get equivariant semilinear morphisms $X(g_{\D,\gamma}):X(\D)\to X(\D')$ for each $\D\in\S$. From \cref{Def semilinear morphism of div fan} part \eqref{Def semilinear morphism of div fan part 2}, we have the following commutative diagram of equivariant semilinear morphisms of pp-divisors
\[\xymatrixcolsep{7pc}\xymatrix{ X(\D) \ar[r]^{X(g_{\D,\gamma})} & X(\D') \\ 
X(\D\cap\E) \ar[r]_{X(g_{\D\cap\E,\gamma})} \ar[u]^{X(\id_{Y},\id_{N},1)} \ar[d]_{X(\id_{Y},\id_{N},1)} & X(\D'\cap \E') \ar[u]_{X(\id_{Y'},\id_{N'},1)} \ar[d]^{X(\id_{Y'},\id_{N'},1)} \\ 
X(\E) \ar[r]_{X(g_{\E,\gamma})} & X(\E').}\]
Then, $X(g_{\D,\gamma})|_{X(\D\cap\E)}=X(g_{\E,\gamma})|_{X(\D\cap\E)}$. This implies that all the semilinear morphisms $X(g_{\D,\gamma})$ fit into a semilinear equivariant morphism 
\[
X(g_{\gamma}):X(\S)\to X(\S').
\]
Unfortunately, it is not always possible to associate a semilinear morphism of divisorial fans to every dominant equivariant semilinear morphism $(\varphi_{\gamma},f_{\gamma}):X(\S)\to X(\S')$. For example, when $\S:=\{\D\}$ and $\S':=\{\D'\}$ and none of the pp-divisors is minimal. In this manner, we need to consider a wider family of morphisms of divisorial fans. Recall that a localized semilinear morphism of pp-divisors $\D\in\PPDiv_{\Q}(Y,\omega_{\D})$ and $\D'\in\PPDiv_{\Q}(Y',\omega_{\D'})$ is represented by a pair of semilinear morphisms of pp-divisors 
\[\xymatrix{ & \kappa^{*}\D \ar[dl]_{(\kappa,\id_{N},\mathfrak{1})} \ar[rd]^{(\psi_{\gamma},F,\mathfrak{f})} & \\ \D & & \D',}\]
where $\kappa:\tilde{Y}\to Y$ is a projective morphism of $L$-varieties, $\psi_{\gamma}:\tilde{Y}\to Y'$ is a projective semilinear morphism, $\tilde{Y}$ is a normal semi-projective variety over $L$ and the morphism $(\psi_{\gamma},F,\mathfrak{f})$ is dominating. In order to simplify the notation, we will denote by $(\kappa,\psi_{\gamma},F,\mathfrak{f}):\D\to\D'$ the localized semilinear morphism of pp-divisors.

\begin{definition}\label{definition localized semilinear morphism of divisorial fans}
Let $\gamma\in\Gamma_{L}$. Let $(\S,Y)$ and $(\S',Y')$ be two divisorial fans over $L$. A \textit{localized semilinear morphism of divisorial fans with respect to $\gamma$} is a family of localized semilinear morphisms of pp-divisors $\mathcal{M}_{\gamma}$ such that 
\begin{enumerate}[i)]
	\item \label{definition localized semilinear morphism of divisorial fans part i} for every $\D\in\S$ and $\D'\in\S'$ there exists at most one localized semilinear morphism of pp-divisors $(\kappa,\psi_{\gamma},F,\mathfrak{f}):\D\to\D'$ in $\mathcal{M}_{\gamma}$, denoted by $\mathcal{M}_{\D,\gamma}$, and
	\item \label{definition localized semilinear morphism of divisorial fans part ii} let $\E,\D\in \S$ and $\E',\D'\in\S'$ such that $\E\preceq \D$ and $\E'\preceq \D'$. If there exist morphisms $\mathcal{M}_{\D,\gamma}:\D\to \D'$ and $\mathcal{M}_{\E,\gamma}:\E\to\E'$ in $\mathcal{M}_{\gamma}$, then the following diagram commutes
	\[ \xymatrix{ \D \ar[r]^{\mathcal{M}_{\D,\gamma}} & \D' \\ \E \ar[u]^{(\id_{Y},\id_{N},\mathfrak{1})} \ar[r]_{\mathcal{M}_{\E,\gamma}} & \E'. \ar[u]_ {(\id_{Y'},\id_{N'},\mathfrak{1}')} } \]
\end{enumerate}
\end{definition}
 
With every localized semilinear morphism of divisorial fans $\mathcal{M}_{\gamma}:(\S,Y)\to (\S,Y')$ we can associate a semilinear equivariant morphism $\mathcal{X}(\mathcal{M}_{\gamma}):\mathcal{X}(\S)\to \mathcal{X}(\S')$, constructed in the same way as for semilinear morphisms of divisorial fans, by \cref{lemma transitivity relation}.

\begin{theorem}\label{Theorem localised semilinear morphisms of fans}
Let $k$ be a field, $L/k$ be a Galois extension with Galois group $\Gamma_{L}$ and $\gamma\in\Gamma_{L}$. Let $(\S,Y)$ and $(\S',Y')$ be two divisorial fans over $L$. Let $(\varphi_{\gamma},f_{\gamma}):\mathcal{X}(\S)\to \mathcal{X}(\S')$ be a dominant semilinear equivariant morphism. If for every $\D\in\S$ there exists $\D'\in\S'$ such that $f_{\gamma}(\mathcal{X}(\D))\subset \mathcal{X}(\D')$, then there exists a localized semilinear morphism of divisorial fans $\mathcal{M}_{\gamma}:(\S,Y)\to (\S',Y')$ such that $\mathcal{X}(\mathcal{M}_{\gamma})=(\varphi_{\gamma},f_{\gamma})$. 
\end{theorem}

\begin{proof}
First notice that, for every $\D\in\S$ and $\D'\in\S'$ such that $f_{\gamma}(\mathcal{X}(\D))\subset \mathcal{X}(\D')$, we have that the induced semilinear equivariant morphism 
\[(\varphi_{\gamma},f_{\gamma})|_{\mathcal{X}(\D)}:\mathcal{X}(\D)\to \mathcal{X}(\D')\]
 is dominant. Let us prove this last assertion. Let $U\subset \mathcal{X}(\D')$ be an open subvariety, then $U\subset \mathcal{X}(\S')$ is an open subvariety. Given that $(\varphi_{\gamma},f_{\gamma})$ is dominant, $f_{\gamma}^{-1}(U)$ is a nonempty subvariety of $\mathcal{X}(\S)$. This implies that $\mathcal{X}(\D)\cap f_{\gamma}^{-1}(U)$ is not empty and, therefore, $f_{\gamma}(\mathcal{X}(\D))\cap U$ is not empty. Then, the restrictions $(\varphi_{\gamma},f_{\gamma})|_{\mathcal{X}(\D)}:\mathcal{X}(\D)\to \mathcal{X}(\D')$ are dominant for every $\D\in\S$.
 
 Given that, for every $\D\in\S$, the restriction map $(\varphi_{\gamma},f_{\gamma})|_{\mathcal{X}(\D)}:\mathcal{X}(\D)\to \mathcal{X}(\D')$ is a dominant semilinear equivariant morphism, there exists a localized semilinear morphism of pp-divisors $\mathcal{M}_{\D,\gamma}:\D\to\D'$ such that $\mathcal{X}(\mathcal{M}_{\D,\gamma})=(\varphi_{\gamma},f_{\gamma})|_{X(\D)}$ by \cref{proposition equivalence of categories modified}. Denote $\mathcal{M}_{\gamma}:=\{\mathcal{M}_{\D,\gamma}\}$. By construction, $\mathcal{M}_{\gamma}$ satisfies part \eqref{definition localized semilinear morphism of divisorial fans part i} of \cref{definition localized semilinear morphism of divisorial fans}. In order to prove part \eqref{definition localized semilinear morphism of divisorial fans part ii} of \cref{definition localized semilinear morphism of divisorial fans}, it suffices to do it on the respective varieties by \cref{proposition equivalence of categories modified} and this is obvious from the commutative diagram
\[ \xymatrixcolsep{5pc}\xymatrix{ \mathcal{X}(\D) \ar[r]^{\mathcal{X}(\mathcal{M}_{\D,\gamma})} & \mathcal{X}(\D') \\ \mathcal{X}(\E) \ar[u]^{\mathcal{X}(\id_{Y},\id_{N},\mathfrak{1})} \ar[r]_{\mathcal{X}(\mathcal{M}_{\E,\gamma})} & \mathcal{X}(\E') \ar[u]_ {\mathcal{X}(\id_{Y'},\id_{N'},\mathfrak{1}')} .} \]
Finally, we have that $\mathcal{X}(\mathcal{M}_{\gamma})=(\varphi_{\gamma},f_{\gamma})$, since $\mathcal{X}(\mathcal{M}_{\gamma})$ is the gluing of the restriction maps $(\varphi_{\gamma},f_{\gamma})|_{\mathcal{X}(\D)}$. Then, the assertion holds.
\end{proof}

 A semilinear equivariant morphism between two normal varieties endowed with an effective torus action does not necessarily satisfy the hypothesis of \cref{Theorem localised semilinear morphisms of fans}. However, given that such varieties have several divisorial fans, we can always find a pair of divisorial fans satisfying the condition of the theorem.
 
 The following is a tool that allows us to construct suitable divisorial fans.

\begin{lemma}\label{Lemma equivariant divisorial fans}
Let $\gamma\in\Gamma_{L}$. Let $T$ and $T'$ be split algebraic tori over $L$. Let $X$ be a normal $T$-variety and $X'$ be a normal $T'$-variety, both over $L$. If $(\varphi_{\gamma},f_{\gamma}):X\to X'$ is a semilinear equivariant morphism, then there exist divisorial fans $(\S,Y)$ and $(\S',Y')$ over $L$ such that $X\cong \mathcal{X}(\S)$, $X'\cong \mathcal{X}(\S')$; and for every $\D\in\S$ there exists $\D'\in\S'$ satisfying $f_{\gamma}(\mathcal{X}(\D))\subset \mathcal{X}(\D')$.
\end{lemma}

\begin{proof}
By Sumihiro's Theorem, $X'$ has a $T'$-stable affine open covering $\mathcal{U}'$. Again by Sumihiro's Theorem, for each $U'\in\mathcal{U}'$, $f_{\gamma}^{-1}(U')$ has a $T$-stable affine open covering $\mathcal{U}_{U'}$. Since the $f^{-1}(U')$ form $T$-stable open covering, we have that $\mathcal{U}:=\bigcup_{U'\in\mathcal{U}'}\mathcal{U}_{U'}$ is a $T$-stable affine open covering of $X$. Moreover, this covering is such that for every $U\in\mathcal{U}$ there exists $U'\in\mathcal{U}'$ satisfying $f_{\gamma}(U)\subset U'$. Then, applying \cref{remarkcovdivfancorrespondence} to these affine open coverings, we get the divisorial fans satisfying the conditions of the lemma.
\end{proof}

\begin{corollary}
Let $\gamma\in\Gamma_{L}$. Let $T$ and $T'$ be split algebraic tori over $L$. Let $X$ be a normal $T$-variety and $X'$ be a normal $T'$-variety, both over $L$. Let $(\varphi_{\gamma},f_{\gamma}):X\to X'$ be a dominant semilinear equivariant morphism. Then, there exist divisorial fans $(\S,Y)$ and $(\S',Y')$ such that $X\cong \mathcal{X}(\S)$, $X'\cong \mathcal{X}(\S')$ and a localized semilinear morphism of divisorial fans $\mathcal{M}_{\gamma}$ such that $\mathcal{X}(\mathcal{M}_{\gamma})=(\varphi_{\gamma},f_{\gamma})$.
\end{corollary}

\begin{proof}
By \cref{Lemma equivariant divisorial fans}, there exist two divisorial fans $(\S,Y)$ and $(\S',Y')$ over $L$ satisfying the hypothesis of \cref{Theorem localised semilinear morphisms of fans}. Then, there exists a localized semilinear morphism of divisorial fans $\mathcal{M}_{\gamma}:(\S,Y)\to(\S',Y')$ such that $\mathcal{X}(\mathcal{M}_{\gamma})=(\varphi_{\gamma},f_{\gamma})$. Then, the assertion holds.
\end{proof}

For normal toric varieties, from a semilinear equivariant isomorphism, we get a semilinear isomorphism of fans. 

\begin{corollary} \label{corollary semilinear automorphisms toric}
Let $\gamma\in\Gamma_{L}$. Let $X_{\Delta}$ and $X_{\Delta'}$ be two normal split toric varieties over $L$ and $(\varphi_{\gamma},f_{\gamma}):X_{\Delta}\to X_{\Delta'}$ be a semilinear equivariant isomorphism. Then, there exists a triple $(\psi_{\gamma},F,\mathfrak{f})$, where $\psi_{\gamma}=\gamma^{\natural}:\Spec(L)\to \Spec(L)$, $F:N\to N'$ is an isomorphism of lattices such that $F(\omega)=\omega'$ for $\omega\in\Delta$ and $\omega'\in\Delta'$ satisfying $f_{\gamma}(X_{\omega})=X_{\omega'}$ and $\mathfrak{f}\in N\otimes L^{*}$ such that $(\varphi_{\gamma},f_{\gamma})=X(\psi_{\gamma},F,\mathfrak{f})$.
\end{corollary}

\begin{proof}
This result follows from \cref{Theorem localised semilinear morphisms of fans} applied to the coverings given by the fans. The last part comes from the fact that $\mathfrak{f}|_{\omega}=\mathfrak{f}$ for any cone, where the cone is seen as a pp-divisor.
\end{proof}

\subsection{The group of semilinear automorphisms of a divisorial fan}\label{Section semilinear aut div fan}
\noindent

In the following, we will focus on the case where $(\S,Y)$ and $(\S',Y')$ are the same divisorial fan.

\begin{definition}
Let $\gamma\in\Gamma_{L}$. Let $(\S,Y)$ be a divisorial fan over $L$. A \textit{semilinear automorphism with respect to $\gamma$ of} $(\S,Y)$ is a semilinear morphism of divisorial fans $g_{\gamma}:(\S,Y)\to(\S,Y)$ such that, for every $\D\in\S$, there exists a unique $g_{\gamma}(\D)\in\S$ such that $g_{\D,\gamma}:\D\to g_{\gamma}(\D)$ is a semilinear isomorphism of pp-divisors. If $\gamma$ is the neutral element of $\Gamma_{L}$, then we say that $g:=g_{\gamma}$ is an \textit{automorphism of divisorial fans}.
\end{definition}

Let $g_{\gamma}$ and $h_{\beta}$ be semilinear automorphisms of $(\S,Y)$. Then the composition $h_{\beta}\circ g_{\gamma}$ is a semilinear automorphism $r_{\beta\gamma}$ satisfying: 
\begin{enumerate}[a)]
	\item For every $\D\in\S$, $r_{\beta\gamma}(\D)=h_{\beta}(g_{\gamma}(\D))$ and 
	\item $r_{\D,\beta\gamma}=h_{g_{\gamma}(\D),\beta}\circ g_{\D,\gamma}$. 
\end{enumerate}
Thus, the set of all semilinear automorphisms $\SAut(\S,Y)$ comes with a semi-group structure. Even in the affine case, $\SAut(\D)$ does not necessarily have a group structure. This problem is solved by considering the localized category $\mathfrak{PPDiv}_{S}(L/k)$.

\begin{definition}\label{definition localized semilinear automorphism of divisorial fans}
Let $\gamma\in\Gamma_{L}$. Let $(\S,Y)$ be a divisorial fan over $L$. A \textit{localized semilinear automorphism of divisorial fans with respect to $\gamma$} is a localized semilinear morphism of divisorial fans $\mathcal{M}_{\gamma}:(\S,Y)\to (\S,Y)$ with respect to $\gamma$, such that each $\mathcal{M}_{\D,\gamma}$ is a localized semilinear isomorphism.
\end{definition}

The set of localized semilinear automorphisms $\SAut_{\mathrm{loc}}(\S,Y)$ has a group structure and the functor $\mathcal{X}:\mathfrak{PPDiv}_{S}(L/k)$ induces a map
\begin{align*}
\mathcal{X}:\SAut_{\mathrm{loc}}(\S,Y) &\to \SAut(T_{\S};\mathcal{X}(\S)), \\
\mathcal{M}_{\gamma} &\mapsto \mathcal{X}(\mathcal{M}_{\gamma}),
\end{align*}
 where $T_{\S}$ denotes the split algebraic torus over $L$ associated with $(\S,Y)$.
\begin{proposition}\label{propositiontoricactiontoequivariantaut}
  Let $(\S,Y)$ be a divisorial fan over $L$. The map $\SAut_{\mathrm{loc}}(\S,Y)\to\SAut(T_{\S};\mathcal{X}(\S))$ is a monomorphism of abstract groups. 
\end{proposition}

\begin{proof}
    For each $\mathcal{M}_{\gamma}\in \SAut_{\mathrm{loc}}(\S,Y)$ there exists an equivariant semilinear automorphism $\mathcal{X}(\mathcal{M}_{\gamma})$. Thus, we need to prove that 
    \[\mathcal{X}(\mathcal{N}_{\beta}\circ \mathcal{M}_{\gamma})=\mathcal{X}(\mathcal{N}_{\beta})\circ \mathcal{X}(\mathcal{M}_{\gamma})\] 
    for all $\mathcal{N}_{\beta},\mathcal{M}_{\gamma}\in\SAut_{\mathrm{loc}}(\S,Y)$. 
    
    Let $\D$ be in $\S$. By restricting $\mathcal{X}(\mathcal{N}_{\beta}\circ \mathcal{M}_{\gamma})$ to $\mathcal{X}(\D)$ we have
\[(\varphi_{g},g)|_{X(\D^{i})}=X(\psi_{ij}(g),F_{ij}(g),\mathfrak{f}_{ij}(g))\,\mathrm{ and }\,(\varphi_{h},h)|_{X(\D^{j})}=X(\psi_{jk}(h),F_{jk}(h),\mathfrak{f}_{jk}(h)).\]
\begin{align*}
        \mathcal{X}(\mathcal{N}_{\beta}\circ \mathcal{M}_{\gamma})|_{\mathcal{X}(\D)} &= \mathcal{X}((\mathcal{N}\circ \mathcal{M})_{\D,\beta\cdot \gamma})\\ 
        &=\mathcal{X}(\mathcal{N}_{\mathcal{M}_{\gamma}(\D),\beta}\circ \mathcal{M}_{\D,\gamma}) \\ 
        &= \mathcal{X}(\mathcal{N}_{\mathcal{M}_{\gamma}(\D),\beta})\circ \mathcal{X}(\mathcal{M}_{\D,\gamma}) \\ 
        &= \mathcal{X}(\mathcal{N}_{\beta})|_{\mathcal{X}(\mathcal{M}_{\gamma}(\D))}\circ \mathcal{X}(\mathcal{M}_{\gamma})|_{\mathcal{X}(\D)} \\ 
        &= (\mathcal{X}(\mathcal{N}_{\beta})\circ \mathcal{X}(\mathcal{M}_{\gamma}))|_{\mathcal{X}(\D)}.
    \end{align*}
    Given that $\mathcal{X}((\mathcal{N}\circ \mathcal{M})_{\beta\cdot\gamma})|_{\mathcal{X}(\D)}=(\mathcal{X}(\mathcal{N}_{\beta})\circ \mathcal{X}(\mathcal{M}_{\gamma}))|_{\mathcal{X}(\D)}$ for every $\D$ and the $\mathcal{X}(\D)$ cover $\mathcal{X}(\S)$, we have that $\mathcal{X}((\mathcal{N}\circ \mathcal{M})_{\beta\cdot\gamma})=\mathcal{X}(\mathcal{N}_{\beta})\circ \mathcal{X}(\mathcal{M}_{\gamma})$. Then, the map $\SAut_{\mathrm{loc}}(\S,Y)\to \SAut(T_{\S};(\mathcal{X}(\S))$ is a group homomorphism. The injectivity follows from \cref{proposition equivalence of categories modified}. Then, the proposition holds. 
\end{proof}

\begin{remark}
Notice that $\Aut_{\mathrm{loc}}(\S,Y)$ is a subgroup of $\SAut_{\mathrm{loc}}(\S,Y)$.
\end{remark}

In the following, if it yields no confusion, the group of localized semilinear automorphisms of a divisorial fan will be referred to as $\SAut(\S,Y)$.

\subsection{Semilinear actions of finite groups on divisorial fans}\label{Section semilinear actions of finite groups}
\noindent

A semilinear equivariant automorphism in the image of $\mathcal{X}:\SAut(\S,Y)\to \SAut(\mathcal{X}(\S))$ must preserve the set $\{\mathcal{X}(\D)\mid \D\in\S\}$. Thus this homomorphism is far from surjective, even when $k$ is algebraically closed. For example, let $X:=\A^{1}\times E$, where $E$ is an elliptic curve and $\G_{\mathrm{m},k}$ acts naturally on $\A^{1}_{k}$ and trivially on $E$. Let $a\in E$ be a non-torsion point and consider
\begin{align*}
	(\id_{\G_{\mathrm{m},k}},f_{a}):X &\to X,\\
	(x,y) &\mapsto (x,y+a).
\end{align*}
This automorphism admits no finite $T$-stable affine open cover $\mathcal{U}$ that is also $f_{a}$-stable. Indeed, since $\{0\}\times E$ is $\G_{\mathrm{m},k}$-stable, such a cover must contain an open subset $U=\A_k^1\times V$, with $V\subset E$ affine. By finiteness and $f_a$-stability, $f_a^n(U)=U$ for some $n>0$, hence $f_a^{mn}(U^c)=U^c$ for every $m\in\N$. Since $U^c=\A_k^1\times V^c$, with $V^c$ finite, this contradicts the fact that $a$ is non-torsion. Therefore, there is no divisorial fan $(\S,Y)$ for $X$ such that $(\id_{\G_{\mathrm{m},k}},f_a)$ lies in the image of $\mathcal{X}$.

Even semilinear equivariant automorphisms of finite order need not lie in the image of $\mathcal{X}$. Nevertheless, for every finite subgroup of $\SAut(T_{\S};\mathcal{X}(\S))$, we will show that there exists a divisorial fan $(\S',Y')$ with $\mathcal{X}(\S')\cong\mathcal{X}(\S)$ for which the subgroup lies in the image of
\[
\mathcal{X}:\SAut(\S',Y')\to\SAut(T_{\S'};\mathcal{X}(\S')).
\]

 \begin{proposition}\label{propositioninducetoricactioneveryorder}
  Let $\gamma\in\Gamma_{L}$. Let $T$ be a split algebraic torus over $L$. Let $X$ be a normal $T$-variety over $L$. Let $(\varphi_{\gamma},f_{\gamma}):X\to X$ be a semilinear equivariant automorphism. If $X$ has a $(\varphi_{\gamma},f_{\gamma})$-stable finite affine open covering, then there exists a divisorial fan $(\S,Y)$ and $\mathcal{M}_{\gamma}\in\SAut(\S,Y)$ such that $X\cong X(\S)$ and $X(\mathcal{M}_{\gamma})=(\varphi_{\gamma},f_{\gamma})$.
 \end{proposition}
 
 \begin{proof}
 Let $\mathcal{U}$ be a $(\varphi_{\gamma},f_{\gamma})$-stable affine open covering of $X$. Notice that $\mathcal{U}$ is also $T$-stable. Let $(\S,Y)$ be the divisorial fan associated to $\mathcal{U}$ (cf. \cref{remarkcovdivfancorrespondence}). The existence of $\mathcal{M}_{\gamma}\in\SAut(\S,Y)$ is ensured by \cref{Theorem localised semilinear morphisms of fans}, therefore, the assertion holds.
 \end{proof}

  Not every semilinear equivariant automorphism $(\varphi_{\gamma},f_{\gamma}):X\to X$ has such a covering. However, semilinear equivariant automorphisms of finite order, i.e. $(\varphi_{\gamma},f_{\gamma})^{n}=(\varphi_{\gamma}^{n},f_{\gamma}^{n})=(\id_{T},\id_{X})$, do have one.

\begin{lemma}\label{lemmafstability}
Let $\gamma\in\Gamma_{L}$. Let $T$ be a split algebraic torus over $L$. Let $X$ be a normal $T$-variety and $(\varphi_{\gamma},f_{\gamma}):X\to X$ be a semilinear equivariant automorphism of order $n\in \N$. There exists a $T$-stable affine open covering $\mathcal{U}$ of $X$ that is $(\varphi_{\gamma},f_{\gamma})$-stable.
\end{lemma}

\begin{proof}
    By Sumihiro's Theorem, $X$ has a $T$-stable affine open covering $\mathcal{U'}$. For every $U\in\mathcal{U'}$, $f_{\gamma}(U)$ is an affine open subvariety and also $T$-stable. Indeed, let $x\in f_{\gamma}(U)$ and $t\in T$. There exist $x'\in U$ and $t'\in T$ such that $x=f_{\gamma}(x')$ and $t=\varphi_{\gamma}(t')$. Thus, $t\cdot x=\varphi_{\gamma}(t')\cdot f_{\gamma}(t')=f_{\gamma}(t'\cdot x')\in f_{\gamma}(U)$. Therefore, $f_{\gamma}(U)$ is $T$-stable. Finally, the set $\mathcal{U}$ consisting of all the (finitely many) intersections of the elements of $\{f_{\gamma}^{k}(U) \mid U\in\mathcal{U'}\textrm{ and }k\in\N \}$ is $(\varphi_{\gamma},f_{\gamma})$-stable.
\end{proof}

Then, we have the following result.

\begin{corollary}\label{lemmainducedtoricaction}
Let $\gamma\in\Gamma_{L}$. Let $T$ be a split algebraic torus over $L$. Let $X$ be a normal $T$-variety and $(\varphi_{\gamma},f_{\gamma})$ be a semilinear equivariant automorphism of order $n\in\N$. Then, there exists a divisorial fan $(\S,Y)$ and $\mathcal{M}_{\gamma}\in\SAut(\S,Y)$ such that $X\cong \mathcal{X}(\S)$ and $(\varphi_{\gamma},f_{\gamma})=\mathcal{X}(\mathcal{M}_{\gamma})$.
\end{corollary}

\begin{proof}
By \cref{lemmafstability}, $X$ has a $(\varphi_{\gamma},f_{\gamma})$-stable finite affine open covering. Then, by \cref{propositioninducetoricactioneveryorder}, the assertion holds.
\end{proof}

A normal $T$-variety $X$ can have many divisorial fans and not every divisorial fan admits an action induced by a semilinear equivariant automorphism of $X$. However, \cref{lemmainducedtoricaction} states that for every finite cyclic subgroup of $\SAut(T;X)$, there exists a divisorial fan $(\S,Y)$ such that $X\cong \mathcal{X}(\S)$ and the cyclic group is in the image of $\mathcal{X}:\SAut(\S,Y)\to \SAut(T;X)$. Before we prove such an assertion for any finite subgroup, we present the notion of \textit{semilinear action} for divisorial fans.

\begin{definition}\label{definitiontoricactiononppdiv}
Let $H$ be an abstract group and $(\S,Y)$ be a divisorial fan over $L$. A \textit{semilinear action of $H$ on $(\S,Y)$} is a group homomorphism $\varrho:H\to\SAut(\S,Y)$. If the image of $H$ lies on $\Aut(\S,Y)$, it is said to be an \textit{action of $H$ on $(\S,Y)$}. A \textit{Galois semilinear action over} $(\S,Y)$ is a semilinear action of $\Gamma_{L}$ on $(\S,Y)$ with $\varrho$ a section of the sequence
\[1\to\Aut_{}(\S,Y)\to\SAut_{}(\S,Y)\to\Gamma_{L}.\]
\end{definition}

\begin{lemma}\label{propositionGstablecovering}
Let $T$ be a split algebraic torus over $L$. Let $X$ be a normal $T$-variety and $H$ be a finite group of semilinear equivariant automorphisms of $X$. Then, there exists a $T$-stable affine open covering $\mathcal{U}$ of $X$ that is $H$-stable.
\end{lemma}

\begin{proof}
 By Sumihiro's Theorem, $X$ has a $T$-stable affine open covering $\mathcal{U}'$. As in \cref{lemmafstability}, each $g(U)$ is a $T$-stable affine subvariety. Let us define $R:=\{g(U) \mid U\in\mathcal{U}'\textrm{ and }g\in H\}$. Then the set $\mathcal{U}$, consisting of all possible (finitely many) intersections of elements in $R$, is $H$-stable.
\end{proof}

\begin{proposition}\label{propositionexistencetoricaction}
Let $X$ be a normal variety over $L$ endowed with an effective action of a split algebraic torus over $L$ and $H$ be a finite group of semilinear equivariant automorphisms. Then, there exists a divisorial fan $(\S,Y)$, such that $X\cong \mathcal{X}(\S)$ and that admits a semilinear action of $H$ encoding the one on $X$.
\end{proposition}

\begin{proof}
    By \cref{propositionGstablecovering}, there exists a $T$-stable affine open covering, which is also $G$-stable. Let $(\S,Y)$ be the divisorial fan associated to such an open covering by \cref{remarkcovdivfancorrespondence}. By \cref{Theorem localised semilinear morphisms of fans}, for every $(\varphi_{\gamma},f_{\gamma})\in H$ there exists $\mathcal{M}_{\gamma}\in\SAut(\S,Y)$ such that $\mathcal{X}(\mathcal{M}_{\gamma})=(\varphi_{\gamma},f_{\gamma})$. Hence, given that $\SAut(\S,Y)\to \SAut(T;X)$ is a group homomorphism and that $H$ is in the image of $\SAut(\S,Y)\to \SAut(T;X)$, we have a group homomorphism of $H\to\SAut(\S,Y)$. Then, the assertion holds.
\end{proof}

Both results, \cref{propositiontoricactiontoequivariantaut} and \cref{propositionexistencetoricaction}, can be summarized in the following theorem.

\begin{theorem}\label{thereom equivalence semilinear equivariant action and existence divisorial fans}
Let $H$ be an abstract finite group and $X$ be a normal variety over $L$ endowed with an effective action of a split torus $T$ over $L$. There exists a semilinear equivariant algebraic action of $H$ over $X$ if and only if there exists a divisorial fan $(\S,Y)$ with a semilinear action of $H$ such that $X\cong \mathcal{X}(\S)$.
\end{theorem}

\begin{definition}\label{definitionstability}
Let $N$ be a lattice and $H\leq \Aut(N)$ be a finite group. We say that a divisorial fan $(\S,Y,N)$ is \textit{$H$-stable} if there exists a semilinear action of $H$ over $(\S,Y)$. In this case, for every $g\in H$ and $\D$, we denote by $g(\D)$ the corresponding pp-divisor in part \eqref{Def semilinear morphism of div fan part 1} of \cref{Def semilinear morphism of div fan}. If the semilinear action factors through $\Aut(\S,Y)$, then it is said to be \textit{$H$-stable}.
\end{definition}

\begin{definition}\label{definition semilinear orbit divisorial sub fan}
Let $H$ be an abstract group. Let $(\S,Y)$ be a $H$-stable divisorial fan over $L$. Let $\D\in\S$, we denote by $\S(\D,H)$ the \textit{sub divisorial fan generated by $\D$ and $H$}, which is defined as the smallest divisorial fan containing the set $\{g(\D)\mid g\in H\}$ in $(\S,Y)$.
\end{definition}

\subsection{Proof of the main theorem}\label{Section proof main theorem}
\noindent

We now have nearly all the tools needed to prove \cref{maintheoremofpaper}. The remaining ingredient is a generalization of Sumihiro's Theorem, since a variety with an action of a non-split torus need not admit a $T$-stable affine open cover. For example, consider the natural action of $\mathbb{S}_{\R}^{1}$ on the conic $C:=V(x^{2}+y^{2}+z^{2})\subset \P_{\R}^{2}$. The curve is the disjoint union of the open subset $C\cap U_{z\neq 0}$, a non-trivial $\mathbb{S}_{\R}^{1}$-torsor, and the closed point $C\cap\{z=0\}$ od degree two over $\R$. Every nonempty $\mathbb{S}_{\R}^{1}$-stable open subset contains this torsor. Thus, any such open subset containing the ramaining closed point must equal $C$, which is projective and hence not affine. Consequently, no $\mathbb{S}_{\R}^{1}$-stable affine open cover of $C$ exists.


\begin{theorem}\cite[Lemma 8]{Sum74}\label{Theorem Sumihiro general}
Let $T$ be an algebraic $k$-torus. Let $X$ be a normal variety over $k$ endowed with an action of $T$, then $X$ has a $T$-stable quasi-projective open covering.
\end{theorem}

\begin{remark}
In \cite{Sum74}, Sumihiro states that the hypothesis of normality cannot be relaxed. As an example, consider the cuspidal curve $C := V(y^{2}z - x^{3}) \subset \mathbb{P}_{k}^{2}$. Such a curve has an effective action of $\mathbb{G}_{\mathrm{m},k}$ given by $\lambda \cdot [x:y:z] = [\lambda^{2}x: \lambda^{3}y: z]$. If $C$ were covered by $\mathbb{G}_{\mathrm{m},k}$-stable affine open subvarieties, then the open containing the singular point should intersect the open orbit isomorphic to $\mathbb{G}_{\mathrm{m},k}$. However, this would imply that $C$ is such an affine variety, and that is a contradiction. (This example can be found in \cite{CLS11}).
\end{remark}

We are ready to prove the main theorem of this work.

\begin{proof}[Proof of \cref{maintheoremofpaper}]

Let us prove first \cref{maintheoremofpaper}~\eqref{maintheoremofpaper part a}. Let $(\S_{L},Y_{L})$ be a divisorial fan such that $X_{L}\cong X(\S_{L})$ and that admits a Galois semilinear action such that $X(\S(\D,\Gamma_{L}))$ is quasi-projective for every $\D\in\S_{L}$. By \cref{thereom equivalence semilinear equivariant action and existence divisorial fans}, it is equivalent to a Galois semilinear equivariant action on $X(\S_{L})$ with a $T_{L}$-stable and $\Gamma_{L}$-stable quasi-projective covering. Hence, by \cref{proposition equivalence normal G-varieties}, $X_{L}$ has a $k$-form as a $T$-variety.

Now, let us prove \cref{maintheoremofpaper}~\eqref{maintheoremofpaper part b}. By \cref{Theorem Sumihiro general}, $X$ can be covered by $T$-stable quasi-projective open subvarieties $\{X_{i}\}_{i\in I}$. Given that $T_{L}:=T\times_{\Spec(k)}L$ is split, each $X_{i,L}:=X_{i}\times_{\Spec(k)}L$ can be covered by $T_{L}$-stable open subvarieties $\{U_{ij}\}_{j\in J_{i}}$, by Sumihiro's Theorem. Notice that each $X_{i,L}$ is $\Gamma_{L}$-stable and quasi-projective. Denote by $(\varphi_{\gamma},f_{\gamma})$ the semilinear equivariant morphism induced by the base change, where $f_{\gamma}:=\id_{X}\times \gamma^{\natural}$. For each $\gamma\in\Gamma_{L}$, the subvariety $f_{\gamma}(U_{ij})\subset X_{i,L}\subset X_{L}$ is open and $T$-stable. Therefore, the set 
\[\mathcal{U}:=\{f_{\gamma}(U_{ij}) \mid \gamma\in\Gamma_{L}\,\textrm{ and } U_{ij}\subset X_{i,L}\}\]
 form a $T_{L}$-stable affine open covering such that is $\Gamma_{L}$-stable and for any $U\in \mathcal{U}$ the union $\cup_{\gamma\in\Gamma_{L}}(f_{\gamma}(U))$ is quasi-projective. Moreover, the covering can be considered stable under intersections. Let $(\S_{L},Y_{L})$ be the divisorial fan associated to $\mathcal{U}$ (cf. \cref{remarkcovdivfancorrespondence}). By \cref{proposition equivalence of categories modified}, for every $\gamma\in\Gamma_{L}$, there exists a localized semilinear automorphism of divisorial fans $\mathcal{M}_{\gamma}\in\SAut(\S_{L},Y_{L})$ such that $X(\mathcal{M}_{\gamma})=(\varphi_{\gamma},f_{\gamma})$. This gives the following commutative diagram of abstract groups 
 \[\xymatrix{ & \Gamma_{L} \ar[rd]^{} \ar[ld] & \\ \SAut(\S_{L},Y_{L}) \ar[rr] & & \SAut(T_{L}; X(\S_{L})). }\]
 Hence, there exists a semilinear action of $\Gamma_{L}$ on $(\S_{L},Y_{L})$. Moreover, for every $\D\in\S_{L}$ we have that $X(\S(D,\Gamma_{L}))$ is a quasi-projective variety. Therefore, $(\S_{L},Y_{L})$ is a $\Gamma_{L}$-stable divisorial fan and $X_{L}\cong X(\S_{L})$ as $T_{L}$-varieties.
\end{proof}

\begin{remark}
    The quasi-projectivity condition in \cref{maintheoremofpaper}, by \cite[Corollary 3.28]{PS11}, can be tested on the divisorial fan when $X$ is a complexity-one normal $T$-variety.
\end{remark}

\begin{example}
Let $k$ be a field and $L/k$ be a quadratic extension with Galois group $\Gamma_{L}$. Denote by $\gamma \in \Gamma_{L}$ the nontrivial element of the Galois group. Let us consider $\P_{L}^{3}$ with the action of $\G_{\mathrm{m},L}$ given by 
\[(\lambda, \mu) \cdot [x_{0}:x_{1}:x_{2}:x_{3}] = [\lambda x_{0}:\mu x_{1}:\lambda\mu x_{2}:x_{3}],\]
as \cref{exampleP3} where all the pp-divisors are defined over the quotient $\P^{3}_{L} \dasharrow \P_{L}^{1}$ given by $[x_{0}:x_{1}:x_{2}:x_{3}]\mapsto [x_{0}x_{1}:x_{2}x_{3}]$. Recall the divisorial fan is as follows.

\begin{center}
	\begin{tikzpicture}
	\draw [black, xshift=1cm] plot [smooth, tension=1] coordinates { (0,0) (3,0.3) (6,-0.3) (9,0)};
	\node[scale=0.9] at (11,0) {$\P_{L}^{1}$};
	\node[scale=0.6] at (3,0.29) {$\bullet$};
	\node[scale=0.6] at (6,-0.17) {$\bullet$};
	\node[scale=0.6] at (9,-0.2) {$\bullet$};
	\node[scale=0.9] at (3,0) {$0$};
	\node[scale=0.9] at (6,-0.6) {$p$};
	\node[scale=0.9] at (9,-0.5) {$\infty$};
	
	\fill[gray] (3,3.5) -- (3,4) -- (4,4) -- (4,3) -- (3.5,3) ;
	\fill[gray] (6,3) -- (6,4) -- (7,4) -- (7,3);
	\fill[gray] (9,3) -- (9,4) -- (10,4) -- (10,3) -- (9,3) ;
	
	\fill[lightgray] (3,3.5) -- (3.5,3) -- (3.5,2) -- (2,2) -- (2,3.5) ;
	\fill[lightgray] (6,3) -- (6,2) -- (5,2) -- (5,3);
	\fill[lightgray] (8.5,2.5) -- (8,2.5) -- (8,2) -- (8.5,2) ;
	
	\fill[darkgray] (3,3.5) -- (3,4) -- (2,4) -- (2,3.5) -- (3,3.5);
	\fill[darkgray] (6,3) -- (5,3) -- (5,4) -- (6,4);
	\fill[darkgray] (9,3) -- (8.5,2.5) -- (8,2.5) -- (8,4) -- (9,4);
	
	\fill[black] (3.5,3) -- (4,3) -- (4,2) -- (3.5,2) ;
	\fill[black] (6,3) -- (7,3) -- (7,2) -- (6,2);
	\fill[black] (8.5,2.5) -- (9,3) -- (10,3) -- (10,2) -- (8.5,2);
	
	\node[scale=0.9] at (3.5,3.5) {$\Delta_{3}$};
	\node[scale=0.9] at (6.5,3.5) {$\omega_{3}$};
	\node[scale=0.9] at (9.5,3.5) {$\omega_{3}$};
	
	\node[scale=0.9] at (2.5,3.7) {$\textcolor{white}{\omega_{0}}$};
	\node[scale=0.9] at (5.5,3.5) {$\textcolor{white}{\omega_{0}}$};
	\node[scale=0.9] at (8.5,3.2) {$\textcolor{white}{\Delta_{0}}$};
	
	\node[scale=0.9] at (2.7,2.7) {$\Delta_{2}$};
	\node[scale=0.9] at (5.5,2.5) {$\omega_{2}$};
	\node[scale=0.9] at (8.2,2.2) {$\omega_{2}$};
	
	\node[scale=0.9] at (3.8,2.5) {$\textcolor{white}{\omega_{1}}$};
	\node[scale=0.9] at (6.5,2.5) {$\textcolor{white}{\omega_{1}}$};
	\node[scale=0.9] at (9.2,2.5) {$\textcolor{white}{\Delta_{1}}$};
	
	\draw[->] (3,3.5) -- (3,4);
	\draw[-] (3,3.5) -- (3.5,3);
	\draw[->] (3.5,3) -- (4,3);
	\draw[->] (3.5,3) -- (3.5,2);
	\draw[->] (3,3.5) -- (2,3.5);
	
	\draw[->] (6,3) -- (6,4);
	\draw[->] (6,3) -- (7,3);
	\draw[->] (6,3) -- (6,2);
	\draw[->] (6,3) -- (5,3);

	\draw[->] (9,3) -- (9,4); 
	\draw[->] (8.5,2.5) -- (8,2.5); 
	\draw[-] (8.5,2.5) -- (9,3);
	\draw[->] (9,3) -- (10,3); 
	\draw[->] (8.5,2.5) -- (8.5,2); 
	
	\node[scale=0.9] at (11,3) {$\P_{L}^{3}$};
	\draw[->,dashed] (11,2.5) -- (11,0.5);
	\node[scale=0.6] at (11.5,1.5) {$ \sslash \G_{\mathrm{m},L}^{2}$};
	\end{tikzpicture} 
\end{center}

Let us consider the following Galois semilinear equivariant action given by 
\[\varphi_{\gamma}(\lambda,\mu)=(\gamma(\mu),\gamma(\lambda))\, \mathrm{ and }\, f_{\gamma}([x_{0}:x_{1}:x_{2}:x_{3}])=[\gamma(x_{1}):\gamma(x_{0}):\gamma(x_{2}):\gamma(x_{3})].\]
In terms of divisorial fans, this action is given by an automorphism of the divisorial fan given $g_{\gamma}:=(\psi_{\gamma},F,\mathfrak{f})$, where $\psi_{\gamma}([z:w])=[\gamma(z):\gamma(w)]$, $F(a,b)=(-b,-a)$ and $\mathfrak{f}=\mathfrak{1}$. Notice that $g_{\gamma}\D_{0}=\D_{1}$ and $g_{\gamma}\D_{2}=\D_{3}$. This combinatorial datum encodes an action of $\mathrm{Res}_{L/k}(\G_{\mathrm{m},k})$ over $\P_{k}^{3}$.
\end{example}

\section{Complexity one and Applications}\label{Section Complexity one and applications}
\noindent

This section is divided into two parts. First, we restate \cref{maintheoremofpaper} for complexity-one normal $T$-varieties, where the statement takes a simpler form. The second part is devoted to applications.

\subsection{Complexity one}
\noindent

Let $T$ be a split algebraic torus over $k$. Let $X$ be a complexity one affine normal $T$-variety over $k$. In this setting, every complexity-one normal $T$-variety over $k$ arises from a divisorial fan over a regular complete curve (see \cref{remarkdivisorialfancomplexity1}). From now on, every divisorial fan will be considered over a regular complete curve. 

\paragraph{Semilinear morphisms in complexity one.} Let $L/k$ be a Galois extension with Galois group $\Gamma_{L}$ and $\gamma\in\Gamma_{L}$. Let $T$ be a split algebraic torus over $L$. Let $X(\D)$ and $X(\D')$ be two complexity one affine normal $T$-varieties over $L$ and $(\varphi_{\gamma},f_{\gamma}):X(\D)\to X(\D')$ be a semilinear equivariant dominant morphism. By \cite[Theorem 5.11]{MN25}, this morphism arises from a triangle of semilinear morphisms of pp-divisors as follows
\[\xymatrixcolsep{5pc}\xymatrix{  \D & \kappa^{*}\D \ar[l]_-{(\kappa,\id_{N},1)} \ar[r]^-{(\psi_{\gamma},F,\mathfrak{f})} & \D‘,}\]
 where $\kappa:\tilde{Y}\to Y$ is a projective birational map and $\tilde{Y}$ is a regular semi-projective variety that is actually projective. However, given that $Y$ is a regular complete curve, $\kappa:\tilde{Y}\to Y$ can be considered as the identity. Then, in this setting, \cref{Appendix localization} is not needed. Thus, the complexity one version of \cite[Theorem 5.11]{MN25} is the following.

\begin{theorem}\label{theorem88complexity1}
Let $\D$ and $\D'$ be two pp-divisors over complete curves and $(\varphi_{\gamma},f_{\gamma}):X(\D)\to X(\D')$ be a dominant semilinear equivariant morphism. Then, there exists a dominating semilinear morphism of pp-divisors $(\psi_{\gamma}, F,\mathfrak{f}):\D\to\D'$ such that $(\varphi_{\gamma},f_{\gamma})=X(\psi,F,\mathfrak{f})$. Moreover, if $(\varphi_{\gamma},f_{\gamma})$ is a semilinear equivariant isomorphism, then $\psi_{\gamma}$ is a semilinear isomorphism of varieties over $L$, $F$ is a lattice isomorphism such that $F(\omega_{\D})=\omega_{\D'}$ and $\psi_{\gamma}^{*}\D'=F_{*}\D+\mathrm{div}(\mathfrak{f})$.
\end{theorem}

If we restrict the category $\mathfrak{PPDiv}(L/k)$ to the full subcategory of pp-divisors over regular complete curves, the functor $X:\mathfrak{PPDiv}(L/k)\to\mathcal{E}(L/k)$ becomes an equivalence of categories with the category of complexity one affine normal varieties endowed with a split algebraic torus action and dominant semilinear equivariant morphisms. 

\begin{proposition}
The functor $\D\to X(\D)$ is an equivalence of categories between the category of pp-divisors over complete curves and whose morphisms are dominant semilinear morphisms of pp-divisors, with the category of complexity-one normal varieties endowed with a split algebraic torus action over $L$ and dominant semilinear equivariant morphisms.
\end{proposition}

 As we saw in \cref{sectioneqaut}, for a semilinear equivariant automorphism of finite order $(\varphi_{\gamma},f_{\gamma}):X \to X$ there exists a $\langle (\varphi_{\gamma},f_{\gamma}) \rangle$-stable divisorial fan $(\S_{L},Y_{L})$ (cf. \cref{lemmainducedtoricaction}). For higher complexity, the data of $\langle (\varphi_{\gamma},f_{\gamma}) \rangle$ over $(\S_{L},Y_{L})$ could not be given by a semilinear morphism of divisorial fans. However, we have the following version of \cref{Theorem localised semilinear morphisms of fans}.
 
 \begin{proposition}\label{proposition theorem 8.8 split semilinear}
 Let $X$ be a normal $L$-variety endowed with an effective action of a split algebraic torus $T$ over $L$. Let $(\varphi_{\gamma},f_{\gamma}):X\to X$ be a semilinear equivariant automorphism. If there exists a finite $T$-stable affine open covering $\mathcal{U}$ of $X$ such that $f_{\gamma}(\mathcal{U})\subset \mathcal{U}$, then there exists a divisorial fan $(\S_{L},Y_{L})$ over a regular complete curve $Y_{L}$ such that $X\cong X(\S_{L})$ and semilinear automorphism of divisorial fans $(\psi_{\gamma},F,\mathfrak{f}):(\S,Y_{L})\to(\S_{L},Y_{L})$ such that $X(\psi_{\gamma},F,\mathfrak{f})=(\varphi_{\gamma},f_{\gamma})$.
 \end{proposition}
 
 \begin{proof}
 Let $(\S,Y)$ be the divisorial fan corresponding to $\mathcal{U}$ (see \cref{remarkcovdivfancorrespondence}). For every $\D\in\S$, the restricted semilinear equivariant isomorphism \[(\varphi_{\gamma},f_{\gamma})|_{X(\D)}:X(\D)\to X(f_{\gamma}(\D))\] arises from a semilinear isomorphism of pp-divisor $(\psi_{\gamma,\D},F_{\gamma,\D},\mathfrak{f}_{\gamma,\D}):\D\to f_{\gamma}(\D)$ by \cref{theorem88complexity1}.
 
 The isomorphisms $\psi_{\gamma,\D}:Y\to Y$ are equal over a Zariski open subset of $Y$, then they are the same semilinear isomorphism $\psi_{\gamma}:Y\to Y$. Let us prove that $F_{\gamma,\D}=F_{\gamma,\E}$ and $\mathfrak{f}_{\gamma,\D}=\mathfrak{f}_{\gamma,\E}$ for every $\D,\E\in\S$. Let us suppose first that $\E\preceq\D$, then we have the following commutative diagram
 \[\xymatrixcolsep{6pc}\xymatrix{X(\D) \ar[r]^{(\varphi_{\gamma},f_{\gamma})|_{X(\D)}}& X(f_{\gamma}(\D)) \\ X(\E) \ar[u] \ar[r]_{(\varphi_{\gamma},f_{\gamma})|_{X(\E)}} & X(f_{\gamma}(\E)) \ar[u], }\]
 where the vertical arrows are the inclusions. This diagram, in the category of pp-divisors turns to be
 \[\xymatrixcolsep{6pc}\xymatrix{\D \ar[r]^{(\psi_{\gamma},F_{\gamma,\D},\mathfrak{f}_{\gamma,\D})} & f_{\gamma}(\D) \\ \E \ar[u]^{(\id_{Y},\id_{N},\mathfrak{1})} \ar[r]_{(\psi_{\gamma},F_{\gamma,\E},\mathfrak{f}_{\gamma,\E})} & f_{\gamma}(\E) \ar[u]_{(\id_{Y},\id_{N},\mathfrak{1})}. }\]
 Then, by the commutativity of the latter diagram we have that \[(\psi_{\gamma},F_{\gamma,\D},\mathfrak{f}_{\gamma,\D})=(\psi_{\gamma},F_{\gamma,\E},\mathfrak{f}_{\gamma,\E}).\]
 Then, we have that $F_{\gamma,\D}=F_{\gamma,\E}$ and $\mathfrak{f}_{\gamma,\D}=\mathfrak{f}_{\gamma,\E}$ for every $\D,\E\in\S$ and $\E\preceq \D$. Finally, given that for any pair $\D$ and $\E$ in $\S$ we have that $\D\cap \E$ is a face of both, we conclude that the assertion holds.
 \end{proof}
 
 As stated in \cref{remarkseveraldivfan}, a normal $T$-variety of higher complexity can be encoded by infinitely many divisorial fans. However, in complexity one, any two of them can always be compared through a third one.
 
 \begin{definition}
Let $Y$ be a normal variety and $(\S,Y)$ a divisorial fan. We define
\begin{enumerate}
\item the \emph{locus} of $(\S,Y)$ by
\[
\Loc(\S,Y):=\bigcup_{\D\in\S}\Loc(\D),
\]
\item the \emph{support} of $(\S,Y)$ by
\[
\mathrm{Supp}(\S,Y):=
\{D\in\Div(Y)\mid \exists\, \D=\textstyle\sum\Delta_{E}\otimes E\in\S
\text{ such that }\Delta_{D}\neq\omega_{\D}\}.
\]
\end{enumerate}
\end{definition}
 
 \begin{proposition}\label{proposition: comparing div fans}
Let $C$ be a regular projective curve and let $\S,\S'$ be divisorial fans over $C$. Then $X(\S)\simeq X(\S')$ as $T$-surfaces if and only if there exist a divisorial fan $\S''$ over $C$ and monomorphisms of divisorial fans
\[
(\psi,\mathrm{id}_{N},\mathfrak f):(\S,C)\to(\S'',C),\qquad
(\psi',\mathrm{id}_{N},\mathfrak f'):(\S',C)\to(\S'',C)
\]
such that
\[
\psi(\Loc(\S))=\Loc(\S'')=\psi'(\Loc(\S')).
\]
\end{proposition}

\begin{proof}
Assume that $X(\S)\simeq X(\S')$ as $T$-varieties and let $f:X(\S)\to X(\S')$ be a $T$-equivariant isomorphism. 
Let $\mathcal U=\{X(\D)\mid\D\in\S\}$
be the $T$-stable affine open covering of $X(\S)$ defined by $\S$, and 
$\mathcal{U}'=\{f^{-1}(X(\D'))\mid \D'\in\S'\}$
the covering induced by $f$. 
Let $\mathcal{U}''$ be their common refinement obtained by taking all finite intersections. 
This refinement corresponds to a divisorial fan $(\S'',C)$ containing $(\S,C)$ and $(\S',C)$, and the inclusion of coverings induces the required monomorphisms of divisorial fans.

Conversely, the monomorphisms $(\psi,\mathrm{id}_{N},\mathfrak f):(\S,C)\to(\S'',C)$ and $(\psi',\mathrm{id}_{N},\mathfrak f'):(\S',C)\to(\S'',C)$ induce $T$-equivariant open embeddings $f_1:X(\S)\to X(\S'')$ and $f_2:X(\S')\to X(\S'')$.
The condition $\psi(\Loc(\S))=\Loc(\S'')=\psi'(\Loc(\S'))$
implies that both $f_1$ and $f_2$ are surjective, hence isomorphisms. 
Therefore $f_2^{-1}\circ f_1:X(\S)\to X(\S')$ is a $T$-equivariant isomorphism.
\end{proof}

\begin{remark}
This result does not hold in higher complexity; it may already fail in complexity two. 
The main obstruction is that the third divisorial fan cannot in general be chosen over the same semi-projective variety.
\end{remark}

\paragraph{Descent in complexity one.}In \cite{PS11}, Petersen and S\"{u}\ss\ provided a caracterization for a divisorial fan to be \emph{quasi-projective}: if there exists a $T$-invariant divisor as in \cite[Corollary 3.28]{PS11}. The following result is a direct consequence of \cite[Corollary 3.28]{PS11}.

\begin{proposition}
    Let $(\S,Y)$ be a divisorial fan over a complete regular curve. Then, $X(\S)$ is a quasi-projective variety if and only if $(\S,Y)$ is a quasi-projective divisorial fan.
\end{proposition}

In the complexity one case, \cref{maintheoremofpaper} can be stated in terms of quasi-projective divisorial fans:
\begin{theorem}\label{theoremdescentcomplexityone}
Let $T$ be an algebraic $k$-torus, $k\subset L\subset \overline{k}$ a finite Galois extension that splits the torus and $\Gamma:=\Gal(L/k)$.

        \begin{enumerate}[a)]
            \item Let $X_L$ be a complexity one $T_{L}$-variety over $L$. If there exists a divisorial fan $(\S_{L},Y_{L})$ over a complete curve admitting a Galois semilinear action such that $X_{L}\cong X(\S_{L})$ and for every $\D\in\S_{L}$ the divisorial subfan $\S(\D,\Gamma_{L})$ is quasi-projective, then $X_L$ has a $k$-form $X$ as a $T$-variety.
            \item Let $X$ be a complexity-one normal $T$-variety over $k$. Then, there exists a divisorial fan $(\S_{L},Y_{L})$ over a complete curve admitting a Galois semilinear action such that $X_{L}\cong X(\S_{L})$ as $T_{L}$-varieties over $L$ and for every $\D\in\S_{L}$ the divisorial subfan $\S(\D,\Gamma_{L})$ is quasi-projective.
        \end{enumerate}
\end{theorem}

\subsection{Applications}
\noindent

Among the applications of \cref{maintheoremofpaper} we can find the classification of smooth real Fano varieties with torus action. For example, Hendrik S\"{u}\ss\ in \cite{S14} studies complex Fano threefolds with two dimensional torus action that admit a Kähler-Einstein metric. 

\subsubsection{Non-split toric varieties}
\noindent

  Let $k$ be a field, $T$ an algebraic torus over $k$, and $L/k$ a finite Galois extension with Galois group $\Gamma_{L}$ that splits $T$. Let $X$ be a normal $T$-toric variety over $k$. Since $T_{L}$ is split, there exists a fan $\Sigma$ such that $X_{L}\cong X_{\Sigma}$. By \cref{maintheoremofpaper}, there exists a divisorial fan $(\tilde{\S}_{L},\Spec(L))$ endowed with a Galois semilinear action such that $X_{L}\cong X(\tilde{\S}_{L})$ and $X(\S(\D,\Gamma_{L}))$ is quasi-projective for every $\D\in\tilde{\S}_{L}$.

  Consider the divisorial fan
  \[\S:=\{\D_{\omega}:=\omega\otimes\Spec(L) \mid \omega\preceq \omega_{\D}\textrm{ for some $\D\in \tilde{\S}_{L}$}\}.\]
  The Galois semilinear action extends to $(\S,\Spec(L))$, and $X(\S)\cong X(\Sigma)$ with $\mathrm{Tail}(\S)=\Sigma$. Moreover,
  \[
  X(\S(\D,\Gamma_{L}))=X_{\Sigma(\omega_{\D},\Gamma_{L})}
  \]
  is quasi-projective for every $\D\in \S$. Hence, by \cite[Proposition 1.9]{Hur11}, $\Sigma(\omega_{\D},\Gamma_{L})$ is a quasi-projective fan.

  Finally, the Galois semilinear action induces a group homomorphism $\Gamma_{L}\to \Aut(N)$, given by $\gamma\mapsto F_{\gamma}$, with $F_{\gamma}(\Sigma)=\Sigma$, and therefore a homomorphism $\Gamma_{L}\to\Aut(\Sigma)$. Thus, $\Sigma$ is $\Gamma_{L}$-stable and $\Sigma(\omega_{\D},\Gamma_{L})$ is quasi-projective for every $\omega\in\Sigma$.
  
\begin{corollary}\cite[Theorem 1.22]{Hur11}
Let $N$ be a lattice and $M\cong \Hom(N,\Z)$ its dual lattice. Let $T:=\Spec(L[M])$ be a split algebraic torus over $L$ and $X(\Sigma)$ be the normal $T$-toric variety associated to a fan $\Sigma$ in $N_{\Q}$. If there exists a group homomorphism $F:\Gamma_{L}\to \Aut(N)$ such that $\Delta$ is $\Gamma_{L}$-stable and the fans $\Sigma(\omega,\Gamma_{L})$ are quasi-projective for all $\omega\in\Sigma$, then $X(\Sigma)$ has a $k$-form as a toric variety.
\end{corollary}

As a consequence of \cref{maintheoremofpaper}, we have the following corollary.

\begin{corollary}
Let $N$ be a lattice and $M\cong \Hom(N,\Z)$ its dual lattice. Let $T:=\Spec(L[M])$ be a split algebraic torus over $L$ and $X(\Sigma)$ be the normal $T$-toric variety associated to a fan $\Sigma$ in $N_{\Q}$. There is a bijection between the set of $k$-forms and the isomorphism classes of pairs $(\Sigma,\varrho)$, where $\varrho:\Gamma_{L}\to \SAut(\Sigma,L)$ is a Galois semilinear equivariant action such that $\Sigma(\omega,\Gamma_{L})$ is a quasi-projective fan for every $\omega\in \Sigma$.
\end{corollary}

\begin{example}
Let $\gamma\in\Gal(\C/\R)$ be the non-trivial element of the group. The toric variety $\P_{\C}^{1}$ is given by the fan $\{\Q_{\leq 0},\{0\},\Q_{\geq 0}\}$. The Galois semilinear actions of fans $(\gamma,\id_{\Z},1)$, $(\gamma,-\id_{\Z},1)$ and $(\gamma,-\id_{\Z},-1)$ correspond, respectively, to $\P_{\R}^{1}$ with an action of $\G_{\mathrm{m},\R}$, $\P_{\R}^{1}$ with an action of $\mathbb{S}_{\R}^{1}$ and $V(x^{2}+y^{2})\subset\P_{\R}^{2}$ with an action of $\mathbb{S}_{\R}^{1}$.
\end{example}

\begin{remark}
A Galois semilinear action on a fan is given by triples $(\gamma,F_{\gamma},\mathfrak{f}_{\gamma})$. The map $\mathfrak{f}:\Gamma_{L}\to L(Y,N)^{*}=N\otimes L^{*}$, satisfies the relation
\[f_{\gamma\gamma'}=F_{\gamma *}(f_{\gamma'})\cdot \gamma'^{*}(\mathfrak{f}_{\gamma}).\]
This relation encodes the information of a torsor of the acting torus.
\end{remark}

In the following we study the particular case of the Hirzebruch surface.

\paragraph{Hirzebruch surface}Let $L/k$ be a finite Galois extension with Galois group $\Gamma_{L}$. Let $r\in\N$ and $\mathbb{F}_{r}$ be the toric $L$-variety generated by the fan $\Sigma$ (see \cref{figure: hirzebruch}). Such a variety is known as Hirzebruch surface. If $L$ is a quadratic extension, by \cite[Theorem 1.25]{Hur11}, $\mathbb{F}_{r}$ has a toric $k$-form if and only if the fan $\Sigma$ is $\Gamma_{L}$-stable.

A direct computation gives the following.

\begin{proposition}\label{proposition automorphisms of Hirzebruch surface's fan}
    The group of automorphisms of $\Sigma$ is isomorphic to $C_{2}$ and is generated by \[F=\begin{pmatrix}
        -1 & 0 \\ r & 1
    \end{pmatrix}.\]
\end{proposition}

According to \cite{Dun16}, a toric variety is said to be neutral if it has a $k$-rational point. The following provides an example of a toric variety having just neutral $k$-forms.

\begin{proposition}
    Assume that $L/k$ is a quadratic extension over a finite field. There are two non-isomorphic toric $k$-forms of the Hirzebruch surface $\mathbb{F}_{r}$ and both of them are neutral.
\end{proposition}

\begin{proof}
    By \cite[Theorem 1.25]{Hur11}, the fan $\Sigma$ has to be $\Gamma_{L}$-stable. By \cref{proposition automorphisms of Hirzebruch surface's fan}, the group of automorphisms of $\Sigma$ is generated by 
    \[F=\begin{pmatrix} -1 & 0 \\ r & 1 \end{pmatrix}.\]
    Then, there are two possibilities. On the one hand, $\Gamma_{L}$ acts trivially on $\Sigma$ and the corresponding $k$-form is the Hirzebruch surface over $k$. Hence, it is neutral because a split algebraic torus has no non-trivial torsors. On the other hand, if $\gamma\in\Gamma_{L}$ is the non-trivial element, then $F_{\gamma}=F$. Let $\mathbb{F}_{r}'$ be the respective $k$-form and $T'$ be the respective $k$-form of $\G_{\mathrm{m},L}^{2}$. The respective $\Gamma_{L}$-action on the torus is the following 
    \begin{align*}
        \G_{\mathrm{m},L}^{2} &\to \G_{\mathrm{m},L}^{2}, \\ (t_{1},t_{2}) &\mapsto (\gamma(t_{1})^{-1}\gamma(t_{2})^{r},\gamma(t_{2})).
    \end{align*}
    Then, 
    \[(\G_{\mathrm{m},L}^{2})^{\Gamma_{L}}:=\{(t_{1},t_{2})\in \G_{\mathrm{m},L}^{2} \mid t_{1}\gamma(t_{1})=t_{2}^{r}\quad\mathrm{ and }\quad t_{2}=\gamma(t_{2}) \}.\] Notice that the norm-one torus 
    \[\mathrm{R}_{L/k}^{1}(\G_{\mathrm{m},L})\cong\{(t_{1},1)\in \G_{\mathrm{m},L}^{2} \mid t_{1}\gamma(t_{1})=1\}\leq (\G_{\mathrm{m},L}^{2})^{\Gamma_{L}}.\] This implies that $\mathrm{R}_{L/k}^{1}(\G_{\mathrm{m},L})\leq T'.$ Hence, we have that $T'\cong \mathrm{R}_{L/k}^{1}(\G_{\mathrm{m},L})\times \G_{\mathrm{m},L}$. Since this torus has no non-trivial torsors over a finite field, we have that this $k$-form is also neutral.
\end{proof}

\begin{proposition}\label{proposition Hirzebruch no non trivial forms}
    Let $k$ be a field and $L/k$ be a finite Galois extension with Galois group $\Gamma_{L}$ such that $|\Gamma_{L}|$ is odd, then the Hirzebruch surface $\mathbb{F}_{r}$ over $L$ has a unique toric $k$-form, the trivial one. 
\end{proposition}

\begin{proof}
    By \cite[Theorem 1.25]{Hur11}, the fan $\Sigma$ has to be $\Gamma_{L}$-stable. By \cref{proposition automorphisms of Hirzebruch surface's fan}, the group of automorphisms of $\Sigma$ is $C_{2}$. Then, there exists a morphism of groups $\Gamma_{L}\to C_{2}$. Given that $|\Gamma_{L}|$ is odd, the only group homomorphism is the trivial one. Finally, $\mathbb{F}_{r}$ has a unique $k$-form, the trivial one.
\end{proof}

Let $\mathbb{F}_{r}$ be the Hirzebruch surface over $L$ with torus $T'$ and $T\leq T'$ be a one-dimensional subtorus. Then, $T$ acts effectively over $\mathbb{F}_{r}$. By \cref{theoremahssplitcase}, $\mathbb{F}_{r}$ arises from a divisorial fan $(\S,Y)$ over a complete curve.

\begin{proposition}\label{proposition complexity 1 toric base}
    Let $T'$ be a split algebraic torus and $X$ be a $T'$-toric variety. Let $T\leq T'$ be a codimension one subtorus and $(\S,Y)$ be a divisorial fan such that $X\cong X(\S)$ as $T$-varieties and $Y$ is a complete curve, then $Y\cong \P^{1}$. 
\end{proposition}

\begin{proof}
	By \cref{corollary cone isomorphic base}, we can assume that every $X(\D)$ is also $T'$-stable. 
    Let $t\in T'$ be an element such that $[t]\in T'/T$ is not trivial. This element defines a $T$-equivariant automorphism $t:X\to X$ and, by \cref{proposition theorem 8.8 split semilinear}, there exists an isomorphism $\psi_{t}:Y\to Y$. Such morphisms satisfy $\psi_{tt'}=\psi_{t}\psi_{t'}$. Then we have an action of $T'/T$ over $Y$. This action is effective, because if $\psi_{t}=\id_{Y}$ we have that $t\in T$ by \cref{theorem88complexity1}. Given that $Y$ is a regular complete curve with an effective action of $T'/T$, we have that $Y\cong \P^{1}$.
\end{proof}

From \cref{proposition complexity 1 toric base}, there exists a divisorial fan $(\S,\P_{L}^{1})$ for the Hirzebruch surface $\mathbb{F}_{r}$. By \cref{proposition theorem 8.8 split semilinear}, the semilinear equivariant action of $\Gamma_{L}$ over $\mathbb{F}_{r}$ induces a semilinear equivariant action of $\Gamma_{L}$ over $\P_{L}^{1}$.

By \cite[Proposition 1.1]{Bea10}, we know that $\PGL_{2}(L)$ contains $S_{4}$ if and only if $-1$ is a sum of two squares in $L$.

\begin{proposition}\label{proposition non toric k form P1}
    Let $k$ be a field of characteristic zero and $L/k$ be a finite Galois extension with Galois group $\Gamma_{L}\cong S_{4}$. If $-1$ is a sum of two squares in $k$, then $\P_{L}^{1}$ has a non-toric $k$-form.
\end{proposition}

\begin{proof}
Given that $-1$ is a sum of two squares in $k$, by \cite[Proposition 1.1]{Bea10}, we have that $S_{4}\subset\mathrm{PGL}_{2}(k)$. Then, $\Hom(\Gamma_{L},\mathrm{PGL}_{2}(k))$ has an injective morphism $\alpha:\Gamma_{L}\to\mathrm{PGL}_{2}(k)$. There is an injective map
\[\xymatrix{\Hom(\Gamma_{L},\mathrm{PGL}_{2}(k))/\sim \, \ar@{^{(}->}[r] & H^{1}(\Gamma_{L},\mathrm{PGL}_{2}(L))},\]
where $\sim$ denotes the relation given by conjugacy. Let $[\alpha]$ be the corresponding class of $\alpha$ under this map. This implies that there exists a semilinear action of $\Gamma_{L}$ over $\P_{L}^{1}$ and, therefore, a group homomorphism $\Gamma_{L}\to\Aut(\Sigma)$, where $\Sigma$ is the fan of $\P_{L}^{1}$. However, given that $\Aut(\Sigma)\cong\Z/2\Z$, this implies that there is no semilinear action of $\Gamma_{L}$ over $\P_{L}^{1}$ that is compatible with its toric structure. Thus, the assertion holds.
\end{proof}

\subsubsection{Equivariantly normal varieties}

Let’s assume $k$ is a perfect field of characteristic $p > 0$. Consider a diagonalized finite group scheme $G \subset T$, where $T$ is a split algebraic torus over $k$. Define $X$ as a $G$-normal variety over $k$ (see for further details \cite[Definition 4.1]{Bri24a}). Brion \cite[Remark 4.3]{Bri24b} offers an alternative perspective on $X$ by considering the normal $T$-variety $X^{\#}:=(T \times X)/G$, where the action of $G$ on $T$ is defined as $g \cdot t = g^{-1} t$, and $G$ acts diagonally on the product. Normality of $X^{\#}$ is established by \cite[Proposition 2.1.3]{Bri26}, accompanied by a $T$-equivariant morphism $\Phi:X^{\#}\to T/G$ induced by the first projection $\mathrm{pr}_{1}:T\times X\to T$. The fiber of $\Phi$ over the base point of $T/G$ is geometrically integral, which corresponds precisely to the $G$-normal curve $X$. This establishes a bijection between the set of $G$-normal varieties $X$ and the set of normal $T$-varieties $Z$ equipped with a $T$-equivariant map $\Phi:Z\to T/G$ such that the fiber over the base point is geometrically integral. Consequently, $G$-normal varieties provide a vast collection of examples of normal $T$-varieties.

Both $X$ and $X^{\#}$ allow for the construction of good quotients $X/G$ and $X^{\#}/T$, respectively, which are canonically isomorphic. Since a $G$-normal variety $X$ gives rise to a normal $T$-variety $X^{\#}$, we can express a $G$-normal variety $X$ in terms of a divisorial fan $(\S,X/G)$ by \cref{maintheoremofpaperspli}. Let us now follow the example provided in \cite[Example 2.1.12]{Bri26}.

With respect to the $T$-equivariant morphism $\Phi:X^{\#}\cong X(\S)\to T/G$ into a combinatorial form, this problem is explored in \cite{MNT26} in the specific scenario where $X$ is a curve and $\mu_{d}=G\subset\G_{\mathrm{m}}$ with $d\geq 2$. In this context, it is demonstrated that the $T$-equivariant morphism $\Phi:X^{\#}\to T/G$ is equivalent to selecting an element $\beta \in k(X/G)^{*}$ such that $\mathrm{div}(\beta|_{\Loc(\D)})=-\D(d)$ for every $\D\in\S$, which is referred to as a $\Phi$-\emph{structure}. Consequently, a $G$-normal curve $X$ corresponds to a $\Phi$-\emph{pair} $(\S,\beta)$ over $X/G$. This pair consists of a divisorial fans $\S$ over $X/G$ and a $\Phi$-structure $\beta$ over $\S$. Nevertheless, not every $\Phi$-pair $(\D,\alpha)$ over a smooth curve $Y$ results in a $G$-normal curve. In fact, such pairs are classified as \emph{geometrically integral} $\Phi$-pairs (as defined in \cite[Definition 5.10]{MNT26}). The following proposition corresponds to \cite[Proposition 5.2]{MNT26} and \cite[Proposition 5.11]{MNT26}.

\begin{proposition}\label{lemma invertible element}
\begin{enumerate}[a)]
	Let $G\cong\mu_{d}\subset \G_{\mathrm{m}}$, where $d\geq 2$.
	\item Let $X$ be a $G$-normal curve. Then there exists a geometrically integral $\Phi$-pair $(\S,\beta)$ over $X/G$ such that $\Phi_{\beta}^{-1}(\{1\})\cong X$.
	\item Conversely, given a geometrically integral $\Phi$-pair $(\S,\beta)$ over a smooth curve $Y$, there exists a $G$-normal curve $X$ such that $X/G\cong \mathrm{Loc}(\S)$, $X^{\#}\cong X(\S)$ and $\Phi=\Phi_{\beta}$.
\end{enumerate}
\end{proposition}

Thus, this theory might be use to study $G$-normal varieties for diagonalizable finite group schemes. For example, Brion proves that there are finitely many isomorphism classes of $G$-normal curves with a prescribed quotient and branch locus (see \cite[Theorem 1.2]{Bri26}). In \cite{MNT26}, the authors improves this result by giving an upper bound for the number of classes by applying this theory (see \cite[Proposition 5.18]{MNT26}).





\Addresses

\end{document}